\documentclass[letterpaper, 12pt, oneside]{book}

    \usepackage[left=1.4in, right=1in, top=1in, bottom=0.8in, dvips, pdftex]{geometry}
    
    \usepackage{fancyhdr}
    \renewcommand{\bibname}{Bibliography}

    \usepackage{amsmath,amsthm, amsfonts,amssymb}
    \usepackage{setspace}
    \usepackage{appendix}
    \usepackage{graphicx}
    \usepackage{pstricks, pst-plot}
    \usepackage[stdtext,vcentermath]{youngtab}
\usepackage{times,amssymb,amscd,verbatim,graphics,axodraw}
\usepackage[all]{xy}

\numberwithin{equation}{section}
                                                                                
\newtheorem{prop}{Proposition}
\newtheorem{theorem}[prop]{Theorem}
\newtheorem{corollary}[prop]{Corollary}
\newtheorem{lemma}[prop]{Lemma}

\theoremstyle{definition}
\newtheorem{definition}[prop]{Definition}
\newtheorem{example}[prop]{Example}
\newtheorem{remark}[prop]{Remark}
\numberwithin{prop}{section} 
    
    \newcommand{\nn}{\vspace{12pt}}
      
    \newcommand{\newchap}[1]
	{
        \chapter{#1}
        \thispagestyle{myheadings}
	}

\newcommand{\A}{\mathcal{A}}
\newcommand{\cc}{cc}
\newcommand{\Conf}{\mathrm{C}}
\newcommand{\cb}{\overline{c}}
\newcommand{\D}{\mathcal{D}}
\newcommand{\Dp}{\Delta p}
\newcommand{\Dt}{\overleftarrow{D}}
\newcommand{\HH}{\mathcal{H}}
\newcommand{\Hom}{\mathrm{Hom}}
\newcommand{\s}{\overline{s}}
\newcommand{\Jb}{\overline{J}}
\newcommand{\la}{\lambda}
\newcommand{\lab}{\overline{\la}}
\newcommand{\La}{\Lambda}
\newcommand{\lb}{\mathrm{lb}}
\newcommand{\lh}{\mathrm{lh}}
\newcommand{\lm}{\la^-}
\newcommand{\ls}{\mathrm{ls}}
\newcommand{\nub}{\overline{\nu}}
\newcommand{\Path}{\mathcal{P}}
\newcommand{\qbin}[2]{\genfrac{[}{]}{0pt}{}{#1}{#2}}
\newcommand{\RC}{\mathrm{RC}}
\newcommand{\rcls}{\ls_{rc}}
\newcommand{\rclb}{\lb_{rc}}
\newcommand{\rk}{\mathrm{rk}}
\newcommand{\SA}{\mathcal{SA}}
\newcommand{\tb}{\overline{t}}
\newcommand{\word}{\mathrm{word}}
\newcommand{\wt}{\mathrm{wt}}
\newcommand{\Z}{\mathbb{Z}}
\newcommand{\ellb}{\overline{\ell}}
\newcommand{\ft}{\tilde{f}}
\newcommand{\et}{\tilde{e}}
\newcommand{\nut}{\tilde{\nu}}
\newcommand{\Jt}{\tilde{J}}
\newcommand{\ellt}{\tilde{\ell}}
\newcommand{\m}{{\bf m}}
\newcommand{\tm}{{\bf \tilde{m}}}
\newcommand{\n}{{\bf n}}
\newcommand{\tn}{{\bf \tilde{n}}}
\newcommand{\uu}{{\bf u}}
\newcommand{\vv}{{\bf v}}
\newcommand{\e}{{\bf e}}
\newcommand{\I}{\mathcal{I}} 
\newcommand{\p}{{\bf p}}
\newcommand{\yb}{\overline{y}}
\newcommand{\R}{\mathbb{R}}
\newcommand{\eb}{{\overline{\bf e}}}
\newcommand{\uub}{{\overline{\bf u}}}

\newcommand{\CC}{{\bf C}}
\newcommand{\Cb}{{\overline{{\bf C}}}}

\begin{document}
    

    \pagenumbering{roman}
    \pagestyle{plain}

    %
    %

    \singlespacing

    \begin{center}

        \begin{Large}
         Fermionic Formulas For Unrestricted Kostka Polynomials And 
        Superconformal Characters
        \end{Large}\\\nn
        By\\\nn
        {\sc Lipika Deka}\\
        B.Sc. (St. Stephen's College, Delhi, India) 1997\\
        B.A. (University of Cambridge, UK ) 1999 \\
        M.S. (University of Wisconsin, Madison, USA) 2001\\\nn
        DISSERTATION\\\nn
        Submitted in partial satisfaction of the requirements for the degree of\\\nn
        DOCTOR OF PHILOSOPHY\\\nn
        in\\\nn
        MATHEMATICS\\\nn
        in the\\\nn
        OFFICE OF GRADUATE STUDIES\\\nn
        of the\\\nn
        UNIVERSITY OF CALIFORNIA\\\nn
        DAVIS\\\nn
        Approved:\\\nn\nn
        
        Anne Schilling (Chair)\\
        \rule{4in}{1pt}\\\nn\nn
        
        Jes\'{u}s A. De Loera\\
        \rule{4in}{1pt}\\\nn\nn
   
       Greg Kuperberg\\
       \rule{4in}{1pt}\\\nn
        
        Committee in Charge\\[8pt]\nn
        
       2005\\\nn
       
    \copyright\ Lipika Deka. All rights reserved.

    \end{center}
       
    \newpage
    \centerline{To my parents}
    \newpage

     \doublespacing

    %
    %
    \tableofcontents
    
    \newpage

    %
    %
    

    \begin{flushright}
        \singlespacing
        Lipika Deka\\
        December 2005\\
        Mathematics

    \end{flushright}
    \begin{center}
    {\Large Fermionic Formulas For Unrestricted Kostka Polynomials And 
        Superconformal Characters}
     \end{center}   
    \centerline{\textbf{\underline{Abstract}}}

The problem of finding fermionic formulas for the many generalizations of 
Kostka polynomials and for the characters of conformal field theories has been
a very exciting research topic for the last few decades. In this dissertation we present new
fermionic formulas for the unrestricted Kostka polynomials extending the work of
Kirillov and Reshetikhin. We  also present new
fermionic formulas  for the characters of $N=1$ and $N=2$ superconformal
algebras which extend the work of Berkovich, McCoy and Schilling. 

Fermionic formulas for the unrestricted Kostka polynomials of type
$A_{n-1}^{(1)}$ in the case of symmetric and anti-symmetric crystal paths  
were given by Hatayama et~al. We present fermionic formulas for the 
unrestricted Kostka polynomials of type $A_{n-1}^{(1)}$ for all crystal paths 
based on Kirillov-Reshetihkin modules. Our formulas and method of proof even 
in the symmetric and anti-symmetric cases  are different from the work
of Hatayama et~al. We interpret the fermionic formulas in terms of a 
new set of unrestricted rigged configurations. For the proof we give a 
statistics preserving  bijection from this new set of unrestricted rigged
configurations to the set of unrestricted crystal paths 
 which generalizes a bijection of Kirillov and Reshetikhin.

We present fermionic formulas for the characters of $N=1$
superconformal models \\$SM(p',2p+p')$ and $SM(p',3p'-2p)$,
and the $N=2$ superconformal model with central charge $c=3(1-\frac{2p}{p'})$. 
The method used to derive these formulas is known as the Bailey flow.  
We  show Bailey 
flows  from the nonunitary minimal model $M(p,p')$ with $p,p'$ coprime
positive integers to $N=1$ and $N=2$
superconformal algebras. We derive a new Ramond sector character
formula for the $N=2$ superconformal algebra with central charge
$c=3(1-\frac{2p}{p'})$ and calculate its fermionic formula.

    \newpage

    %
    %

    \chapter*{Acknowledgments}

First of all I take this opportunity to thank my advisor Prof. Anne Schilling for working with me and guiding me through the years of my graduate school at Davis. I thank her for everything she has taught me starting from the very beginning. Without her help and direction this dissertation would not have been possible.  I am grateful to her for being patient with me when the progress 
was slow, for always reminding me not to lose focus, for reading every word I wrote and for correcting my mathematical and English mistakes.   

I would like to thank my committee members, Prof. Jes\'{u}s De Loera  and Prof. Greg 
Kuperberg, for all their help and especially for reading the dissertation in a short period of time. I am thankful to Prof. Monica Vazirani for all the help and encouragement. A special thanks to Prof. Arun Ram for inspiring  me to work on combinatorial representation theory and for all his support during my days at University of Wisconsin, Madison. I take this opportunity to  express my gratitude for all my professors at St. StephenÕs College, Delhi, University of Cambridge, UW-Madison and UC-Davis for all the amazing classes that inspired me to do a PhD in mathematics. 
I am grateful to my high school math teacher Mr. Mahendra Kalita for being the first person to show me the beauty of mathematics and my undergraduate professor Prof. Geetha Venkataraman for her inspiring classes on group theory.  I specially thank  Prof. Georgia Benkart and  Prof. Martin Isaacs from UW-Madison for the most amazing series of algebra and Lie algebra classes. Their lecture notes have been one of the most useful resources for me.     

I am grateful  to Professors Gaberdiel and Hanno Klemm for the email correspondence which helped me navigate through the jungle of literature regarding $N=2$ character formulas. I also acknowledge their help regarding the spectral flow of $N=2$ superconformal algebras. 
 
I greatly appreciate the financial support from the U.S. National Science Foundation (NSF). This work was supported in part by the NSF Grant DMS-0200774.

I am grateful to Celia Davis for all her help, understanding and support during my days at Davis. I thank my office-mate Isaiah Lankham for providing me with the template for the dissertation, for all the help with \LaTeX\ and computers, and for proof-reading part of this dissertation. 
  
I thank my husband Deepankar for being my constant source of inspiration, for re-living his graduate school days with me, and for being the reason for me  to keep going when the going got tough. Without having him by my side this journey would have been  impossible.  I also thank him for teaching me C++ which is an invaluable part of this dissertation and for  proof-reading my dissertation. 

I am grateful to my parents for always believing in me and encouraging me to pursue
a career in mathematics. I would like to thank my sisters Bini and Simi for their love, and my
in-laws  for their  support and patience during the last few years. I thank my friends Sagarika, Paban and  Lakhima for always being there for me no matter what. I thank Geetika, Pawan and Abhilasha for making the first few years of my graduate school a lot of fun. I thank my friend Maya Ahmed for all the inspiring conversations.

    \newpage

    %
    %

    \pagestyle{fancy}
    \pagenumbering{arabic}
       
    %
    %

    \newchap{Introduction}
    \label{chap:intro}
       \section{Summary of the main results} 
Fermionic formulas have been widely researched in mathematics and physics. 
In this dissertation we consider two problems regarding fermionic formulas 
that arise in the context of combinatorial representation theory and conformal 
field theory (CFT). This dissertation is based on  two papers that resulted from 
research performed with Prof.~ Anne Schilling during my  years of graduate 
school at the University of California, Davis. 
 
In Chapter~\ref{chap:kostka},  we present a new fermionic formula
for the unrestricted Kostka polynomials.  This work is based upon the paper 
``New Fermionic formula for the unrestricted Kostka polynomials"   with Anne 
Schilling. An extended abstract of this paper has appeared in the proceedings 
of 17th International conference, \textit{Formal Power Series and Algebraic
Combinatorics} 2005, held at the University of Messina, Italy, in 
June 2005. The full version of the  paper is available as a preprint at 
{\tt http://front.math.ucdavis.edu/math.CO /0509194}. We have submitted 
this paper for publication to \textit{The Journal of Combinatorial Theory, Series A}. 
Our results extend the work of Kerov, Kirillov and Reshetikhin 
\cite{KKR:1986,KR:1988} who used the Bethe Ansatz to find a 
fermionic formula for the Kostka polynomials. This was first extended 
in \cite{KSS:2002} to generalized Kostka polynomials by establishing 
a bijection between the highest weight paths in the tensor products of 
Kirillov-Reshetikhin crystals of type $A_n$ and rigged configurations. 
We prove our new formula for the unrestricted Kostka polynomial case by 
giving a statistics preserving  algorithmic bijection between all crystal paths 
in the tensor products of Kirillov-Reshetikhin crystals of type $A_n$ 
and the corresponding set of rigged configurations.  An explicit 
description of the new set of rigged configurations is presented, which is 
called the set of unrestricted rigged configurations. Our formula when restricted to 
symmetric and anti-symmetric crystals is different from the results of
Hatayama et~al.~\cite{HKKOTY:1999} where fermionic formulas are given 
for the unrestricted Kostka polynomials in these special cases.

In Chapter~\ref{chap:bailey},  we present new fermionic formulas for the 
characters of $N=1$ and $N=2$ superconformal algebras using the method of 
Bailey construction. The work in Chapter~\ref{chap:bailey} is based on  
the  paper ``Non Unitary minimal models, Bailey's lemma and $N=1,2$ 
superconformal algebras" with Anne Schilling. This paper is published 
in \textit{Communications in Mathematical Physics}, Volume 260, number 3  (2005) 711-725.    
We show that there are Bailey flows from the nonunitary minimal 
models $M(p,p')$ for arbitrary coprime positive integers $p,p'$ to $N=1$ and $N=2$ 
superconformal models. The superconformal models are also indexed by a pair of
coprime positive integers $(p,p')$. Denote the $N=1$ superconformal algebras
by $SM(p,p')$ and $N=2$ superconformal algebras by $\mathcal{A}(p,p')$. We  
find Bailey flows specifically from the model $M(p,p')$ to $SM(2p+p'p')$, $SM(p',3p,-2p)$ and 
to $\mathcal{A}(p,p')$ with central charge given by $3(1-\frac{2p}{p'})$. Using the
known fermionic formulas for the minimal models $M(p,p')$ \cite{BMS:1997},
we explicitly calculate the fermionic formulas for the characters of $SM(2p+p'p')$, 
$SM(p',3p,-2p)$ and $\mathcal{A}(p,p')$. Moreover, we
derive a new Ramond sector character for $N=2$ superconformal algebras and 
calculate its fermionic formula. 

The new bijection given in Chapter~\ref{chap:kostka} as well as its inverse 
have been implemented as C++ programs and are included in Chapter~\ref{chap:program}. 
In early stages of the project on unrestricted Kostka polynomials, these programs were 
used extensively to produce data and to 
verify conjectures regarding the unrestricted rigged configurations.
The programs have also been incorporated into MuPAD-Combinat~\cite{MuPAD:2005} 
as a dynamic module by Francois Descouens. 
In Chapter~\ref{chap:program}, we describe three programs 
 which were used to verify different parts of our conjectures
for our results presented in Chapter~\ref{chap:kostka}.  The programs in 
Sections~\ref{prog:allpaths}, \ref{prog:onepath} and \ref{prog:inverse_bijection}
are designed for use by anyone who would like to do calculations using the 
bijection or the inverse
bijection.  Working out the bijection for even a small example 
is time-consuming and very tedious. Therefore, we believe that these programs 
are very helpful to anyone studying unrestricted Kostka polynomials. The program in
Section~\ref{prog:allpaths} also calculates the unrestricted Kostka polynomials.  
  
Having stated the main results, it is worth mentioning that the bridge between the two 
chapters is the fermionic formula. 
A fermionic formula is a $q$-polynomial or a $q$-series that is a specific
sum of products of the $q$-binomial coefficients
\begin{equation}\label{eq:qbinomial}
\left[\begin{array}{c}m\\ n\end{array}\right]_q=\frac{(q)_m}{(q)_n(q)_{m-n}},
\end{equation}
where 
\begin{equation}
(q)_m :=(q;q)_m = \prod_{k=1}^{m}(1-q^k) \quad \text{for} \quad m\in\Z_{>0} 
\quad \text{and} \quad (q)_0=1.
\end{equation}

In Chapter~\ref{chap:kostka}, the fermionic formulas appear in the context of  
combinatorial representation theory and  in Chapter~\ref{chap:bailey}, they 
appear in the context of conformal field theories (CFTs). The following section 
provides a brief background on fermionic formulas, Kostka polynomials and 
CFTs. More details are given in the subsequent chapters.

\section{Background and motivation}
A partition $\la=(\la_1,\cdots,\la_k)$ is a $k$-tuple of positive integers satisfying
$\la_1\ge \la_2\ge \cdots \ge \la_k\ge 0$. Let $l(\la)$ be the length of the partition $\la$
which is the number of nonzero parts.  In symmetric functions theory, 
the ring of symmetric functions have various bases including the monomial 
symmetric functions, Schur functions, and Hall-Littlewood symmetric functions 
\cite{Macdonald}. The Kostka polynomial 
$K_{\la \mu}(q)$, indexed by two partitions $\la$ and $\mu$ is defined as 
the matrix elements of the transition matrix between the  Schur 
functions $s_{\la}(X)$ and the  Hall-Littlewood symmetric functions $P_{\mu}(X;q)$. 
That is:
\begin{equation}
s_{\la}(X)=\sum_{\mu}K_{\la \mu}(q) P_{\mu}(X;q).
\end{equation}
In representation theory, the Kostka polynomials $K_{\la \mu}(q)$  are 
a $q$-analog of the multiplicity of the 
irreducible $sl_n$ representation $V_{\la}$, indexed by the highest weight $\la$, in the
expansion  of the L-fold tensor 
product $V_{(\mu_1)} \otimes \cdots \otimes V_{(\mu_L)}$. Here $\mu=(\mu_1,\cdots,\mu_L)$
is a partition and $V_{(\mu_i)}$ is the symmetric tensor representation of $sl_n$ with 
weight $\mu_i$. 
These polynomials  have been generalized in many ways in algebraic combinatorics. In some 
generalizations, for example \cite{KSS:2002,LLT:1997,LT:2000,SW:1999,Sh:1998a,Sh:1998b,
Sh:1998c,Sh:2002}, the components of the tensor product are replaced by tensor
representations which are not always symmetric. In some other 
generalizations \cite{HKOTT:1999,HKOTY:1999,OSS:2003a,OSS:2003b,
Sch:2005a, SS:2005}, the representations of $sl_n$ are replaced by representations of other 
Kac-Moody algebras \cite{Kac}. There are many combinatorial descriptions of Kostka 
polynomials.  
Lascoux and Sch{\"u}tzenberger~\cite{LS:1978} gave the description  
\begin{equation}
K_{\la \mu}(q)=\sum_{t \in \mathcal{T}(\la,\mu)} q^{c(t)},
\end{equation}
where $\mathcal{T}(\la,\mu)$ is the set of semi-standard Young tableaux 
\cite{Fulton, Macdonald} of shape 
$\la$ and content $\mu$ and where $c(t)$ \cite{Macdonald} is the  charge statistic 
of the tableau $t\in \mathcal{T}(\la,\mu)$. This expression proved the non-negativity of 
the  coefficients of the Kostka polynomials as conjectured by H.O. Foulkes 
\cite{Foulkes:1974}.  

In the mid 1980's, Kirillov and Reshetikhin \cite{KR:1988}  
obtained a new expression for the  Kostka polynomials which is a fermionic formula. 
The formula was a consequence of studying Bethe Ansatz in statistical mechanics. Bethe 
Ansatz is a method to construct the eigenvectors of the Hamiltonian of 
an integrable quantum system.  When $\mu=(1^L)$, the fermionic formula for the 
Kostka polynomial looks like
\begin{equation}\label{eq:KR_form}
K_{\la,\mu}(q)=q^{\frac{L(L-1)}{2}}M(\la,  \mu;q^{-1})
\end{equation}
where
\begin{equation*}
\begin{split}
 M(\la,\mu;q)&=\sum_{\{m\}} q^{c(\{m\})}\prod_{1\le a \le n,i\ge 1} 
\left[\begin{array}{c}p_i^{(a)}
                           +m_i^{(a)}\\ m_i^{(a)} \end{array}\right]_q, \\
 c{(\{m\})}&=\frac{1}{2}\sum_{1\le a,b\le n} C_{a,b} \sum_{i,j\ge 1}\min (i,j)m_i^{(a)} 
m_j^{(b)},\\
 p_i^{(a)}&=L\delta_{a1}-\sum_{1\le b\le n}C_{ab}\sum_{j\ge 1}\min (i,j)m_j^{(b)}.
\end{split}
\end{equation*}     
The  $\{m\}$-indexed sum is over the set $\{m_i^{(a)}\in \Z_{\ge 0} | 1\le a\le n,i \ge 1\}$ 
such that
 for $1\le a \le n, i\ge 1$, 
 $$p_i^{(a)}\ge 0,\quad \sum_{i\ge 1}im_i^{(a)}=\sum_{j>a}\la_j.$$ 
 Here $(C_{ab})_{1\le a,b\le n}$ is the Cartan matrix for $sl_{n+1}$  and the partition $\la$ 
 has at most $n+1$  nonzero parts.
 The importance of the fermionic formula lies in the fact that there are no minus signs. 
 These formulas can be used to study the limiting behavior, which is a key ingredient
 in finding different formulas for the characters related to affine Lie algebras
 and Virasoro algebras. Some examples of such applications can be found in 
\cite{HKKOTY:1999,Kir:1995, W:2002}.
   
To prove that the fermionic formula of the Kostka polynomial is given by equation
\eqref{eq:KR_form}, Kirillov and Reshetikhin  \cite{KR:1988} gave a bijection
between the set $T(\la,\mu)$ and a new combinatorial object called rigged 
configurations. Rigged configurations index the solutions of the Bethe Ansatz 
equations and they are sequences of partitions satisfying certain size
restrictions along with some labellings called riggings for the parts of the partitions. 
This connection between fermionic formulas and Kostka polynomials is the beginning of
a whole new era of research in combinatorial representation theory.        

In 1997 Nakayashiki and Yamada \cite{NY:1997} gave a different representation 
of the Kostka polynomials in terms of  paths:
\begin{equation}\label{eq:path_form}
K_{\la \mu}(q)=\sum_{p\in \overline{\Path}(\la,\mu)} q^{E(p)},
\end{equation}
where a path $p \in \overline{\Path}(\la,\mu) $ (section 2.2.2, Chapter~\ref{chap:kostka}) is
a highest weight  element of weight $\la$  in Kashiwara's crystal base 
\cite{Kash:1990} corresponding to 
the tensor product representation $V_{(\mu_1)} \otimes V_{(\mu_2)}\otimes 
\cdots \otimes V_{(\mu_L)}$ of $sl_n$. The statistic $E(p)$ \cite{NY:1997} 
associated with a path $p$ is called energy. This energy statistic can be defined 
as $E(p)=n(\mu)-\Dt(p)$ where $\Dt(p)$ is defined in section 2.2.2, Chapter~\ref{chap:kostka}
and $n(\mu)=\sum_{1\le i<j\le L} \min(\mu_i,\mu_j)$.
This new representation was derived by realizing that paths are in bijection with the 
set of rigged configurations. This bijection is done by sending a path (which can be 
viewed as a word in the $sl_n$ case) to its Robinson-Schensted \cite{Fulton} 
recording $Q$-tableau, which 
is then sent to the rigged configuration using Kirillov-Reshetikhin bijection. The path 
form of the Kostka polynomials is particularly  important because this definition 
can be generalized to  any Kac-Moody Lie algebras using the crystal base theory.
 Therefore, the Kostka polynomial for any  Kac-Moody Lie algebra
is defined as the generating function of paths when graded by the energy statistic and  
is called the \textit{one dimensional sum}  $X$. The fermionic formula $M$ for the 
\textit{one dimensional sum} was conjectured in full generality by 
Hatayama et~al.~ in \cite{ HKOTT:1999,HKKOTY:1999}.
This is known as the famous $X=M$ \textit{conjecture}. 
Although this conjecture in full generality is 
still open, many special cases have been proved in a series of papers \cite{OSS:2003a,
OSS:2003b, Sch:2005a,SS:2005}. 

Similar to the Kostka polynomials, the unrestricted Kostka
polynomials  $X_{\la \mu}(q)$,
indexed by two partitions $\la$ and $\mu$, can be defined as the matrix elements of
the transition matrix between the monomial symmetric functions
and the modified Hall-Littlewood symmetric functions \cite{Kir:1998,Macdonald}. 
Let $\la$ be a partition with $l(\la)\le n$, and let $P_{\la}(X_n:q)$ and $Q_{\la}(X_n;q)$
be the Hall-Littlewood polynomials \cite{Macdonald}. A modified Hall-Littlewood
polynomial $Q'_{\la}(X_n;q)$ is defined to be
\begin{equation}
Q'_{\la}(X_n;q)=Q_{\la}(X_n/(1-q);q),
\end{equation}  
where the variables $X_n/(1-q)$ are the products $xq^{j-1}, j\ge 1$ for $x\in X_n:=
(x_1,\cdots, x_n)$. Note that $Q'_{\la}(X;0)=s_{\la}(X)$ and $Q'_{\la}(X;1)=h_{\la}(X)$
where $h_{\la}(X)$ is the complete homogeneous symmetric function \cite{Macdonald}.
With this notation the Kostka polynomial can also be defined as
$$Q'_{\la}(X;q)=\sum_{\mu}s_{\mu}(X)K_{\mu \la}(q).$$
The unrestricted Kostka polynomial, $X_{\la \mu}(q)$ is then defined as
$$Q'_{\la}(X;q)=\sum_{\mu}X_{\la \mu}(q)m_\mu(X).$$
Combinatorially ~\cite{HKOTT:1999,HKKOTY:1999,HKOTY:1999,SW:1999,Sh:2002},  
\begin{equation}\label{eq:unpath_form}
X_{\la \mu}(q)=\sum_{p\in \Path(\la,\mu)} q^{E(p)}
\end{equation}
where  $ \Path(\la,\mu) $ (section 2.2.2, Chapter~\ref{chap:kostka}) is the  set of 
all unrestricted paths of weight $\la$. Unrestricted 
paths are elements in the crystal base of the tensor product representation $V_{(\mu_1)} 
\otimes \cdots \otimes V_{(\mu_L)}$ of $sl_n$. The unrestricted Kostka polynomials
we described above correspond to type $A_{n-1}$  Lie algebras.  One should note that
the set of unrestricted paths of weight $\la$ contains the set of highest weight
paths with the same weight vector. 

A fermionic formula for the $A_{n-1}$ unrestricted Kostka polynomials, when 
$\mu$ is a sequence of row partitions or a sequence of  column partitions,  
was proved in~\cite{HKKOTY:1999,Kir:1998}. 
The existence of crystal bases have been conjectured in \cite{HKOTT:1999,HKOTY:1999} 
for all Kirillov-Reshetikhin modules. A Kirillov-Reshetikhin module, $W^{r,s}$ 
is a finite dimensional module over an affine Lie algebra, which corresponds to
the weight vector $s\Lambda_r$, where $\Lambda_r$ is the fundamental weight of
the affine Lie algebra. The corresponding crystal is denoted by $B^{r,s}$.
For $A_{n-1}^{(1)}$, the affine Lie algebra of type $A_{n-1}$ \cite{Kac}, the 
existence of the Kirillov-Reshetikhin crystals are known \cite{KKMMNN:1992,Sh:2002}. 
In the type $A_{n-1}^{(1)}$ case, the weight vector $s\Lambda_r$ is 
a rectangular partition of height $r$ and width $s$. Having the crystal basis, 
it is natural to extend the definition of unrestricted Kostka polynomials 
to tensor products of Kirillov-Reshetikhin 
modules using the path definition \eqref{eq:unpath_form}.  
The fermionic formula for the unrestricted Kostka polynomials of
type $A_{n-1}^{(1)}$ in this general setup has not been studied 
until now. In Chapter~\ref{chap:kostka},
we study these unrestricted Kostka polynomials for tensor products of all 
Kirillov-Reshetikhin modules of type $A_{n-1}^{(1)}$ and present 
new fermionic formulas.

Fermionic expressions for generating functions of unrestricted paths
for type $A_1^{(1)}$ have also recently surfaced in connection with box-ball systems.
The box-ball systems invented by Takahashi and Satsuma \cite{TS:1990} is an 
important soliton cellular automata. It is a discrete dynamical system in which finitely many
balls move along the one dimensional array of boxes under certain rules. The 
interesting relation between the box-ball systems and the crystal base theory has been 
widely studied, see for example \cite{FOY:2000, HHIKTT:2001, HKT:2000,
HKOTY:2002, T:2004}.
Takagi~\cite{T:2004} established a bijection between box-ball systems and a new set of 
rigged configurations to prove a fermionic formula for the $q$-binomial coefficient.
His set of rigged configurations coincides with our set in the type $A_1^{(1)}$ case.
There is also a generalization of Takagi's bijection to the type $A_{n-1}^{(1)}$ 
case~\cite{KOSTY:2005}. Hence, our bijection composed with the  generalized Takagi's 
bijection establishes  a new  connection between box-ball systems and 
the unrestricted Kostka polynomials.

One motivation for seeking an explicit expression for unrestricted Kostka
polynomials is their appearance in generalizations of Bailey's 
lemma~\cite{Bailey:1949}. Bailey's lemma is a powerful method to
prove Rogers-Ramanujan-type identities \cite{R:1894,R:1917,RR:1919}. 
The Bailey transform of \cite{A:1984} starts with a seed identity and then produces an
infinite family of identities. The original Bailey's lemma corresponds to
type $A_1$.  In~\cite{SW:1999} a type $A_n$ generalization 
of Bailey's lemma was conjectured which was subsequently
proven in~\cite{Wa:2002}. A type $A_2$ Bailey chain, which yields an infinite
family of identities, was given in~\cite{ASW:1999}. In these generalizations a key
ingredient was an explicit fermionic formula for the unrestricted Kostka polynomial. 
If the method we used in this dissertation to find the new fermionic formula  for the 
unrestricted Kostka polynomial can be generalized to other Kac-Moody algebras, it
 might trigger further progress towards generalizations of the Bailey lemma to 
 Kac-Moody algebras other than type $A_n$.  Unrestricted rigged configurations for
 the simply-laced type Lie algebras have already been studied in \cite{Sch:2005b}.         

In the physics context, finding explicit formulas for the characters of
the the solvable lattice models has been a fundamental problem.  The minimal models
denoted by $M(p,p')$ are conformal field theories (CFTs)  invented by Belvin, 
Polyakov and Zamolodchikov \cite{BPZ:1984a,BPZ:1984b}. These are conformally 
invariant two dimensional field theories that describe second order phase 
transitions.  The symmetry algebras of these theories  are the  infinite-dimensional 
algebras known as the Virasoro algebras. The Virasoro algebra is generated by 
generators $L_m$ satisfying
\begin{equation}
[L_n,L_m]=(n-m)L_{m+n}+\frac{c}{12}(n^3-n)\delta_{n+m,0},\quad \text{for} \quad m,n \in \Z,
\end{equation}
where $c$ is the central charge given by 
\begin{equation}
c(p,p')=1-\frac{6(p'-p)^2}{pp'},
\end{equation}
where $1\le p< p'$ with $p,p'$ coprime.
Hence the minimal models are indexed by pairs $p,p'$. The conformal dimension for 
this  model is given by
\begin{equation}
 \quad \Delta_{r,s}^{p,p'}=\frac{(p'r-ps)^2-(p'-p)^2}{4pp'},
\end{equation}
where $$1\le r\le p-1, 1\le s\le p'-1.$$
The characters of these models are given in \cite{CIZ:1987,Dob:1987,Rc:1984} as 
\begin{equation}
\chi_{r,s}^{(p'p')}(q)=q^{\Delta_{r,s}^{p,p'}-c/24}\frac{1}{(q)_{\infty}}
\sum_{j=-\infty}^{\infty}(q^{j(jpp'+rp'-sp)}-q^{(jp'+s)(jp+r)}),
\end{equation}      
where $(q)_\infty=\prod_{i=1}^{\infty}(1-q^i)$. This expression is derived via 
 the Feigin and Fuchs construction \cite{FF:1984} of a basis using bosonic generators and hence 
 known as a  bosonic formula. In 1993, Kedem et~al.~\cite{DKMM:1994,KM:1993,KKMM:1993a,
 KKMM:1993b} found  a new expression for such characters in their study of 
 the three state Potts models. The new formula had no minus signs like the bosonic form.  
 They interpreted this new formula as the partition
function of quasi-particles satisfying fermionic exclusion principles, and so called the new 
expression  fermionic.
 
The fermionic expression  for the minimal models $M(p,p')$ are calculated in 
\cite{ BMS:1997,W:2002}, which have the following form:
\begin{equation}\label{eq:fermionic}
\sum_{\m \quad \text{restriction}}q^{\frac{1}{4}\m^tB\m+A\m}\prod_{j=1}^n\left[ \begin{array}{c} 
((1-B)\m+\frac{\uu}{2})_a\\m_a \end{array} \right]_{q},
\end{equation} 
where $\m$ is an $n$ component vector of non-negative integers that may be subject to
restrictions in the sum, $B$ is an $n\times n$ matrix, and $A$ and $\uu$ are $n$ component 
vectors. Kedem et~al.~ showed in \cite{DKMM:1994,KM:1993,KKMM:1993a,
 KKMM:1993b}  that any expression of this form can be interpreted physically. Consequently,
 finding fermionic formulas is a very important problem in physics.   
 The character identities obtained by equating the bosonic and fermionic expression
 for the CFT characters are known as Bose-Fermi identities. 

As mentioned above, the fermionic formulas for minimal models 
are very well-studied, but the fermionic formulas for other CFTs are not yet known 
in full generality. The
$N=1$ and $N=2$ superconformal algebras are two classes of CFTs in which the 
symmetry algebras are  extended Virasoro algebras.  Berkovich, McCoy and Schilling
demonstrated in \cite{BMS:1995} that some of the characters of $N=1$ and $N=2$ 
superconformal algebras can be obtained from the minimal models $M(p-1,p)$ by 
means of a construction known as the Bailey lemma. Bailey's Lemma first 
appeared in the 1949 paper \cite{Bailey:1949}. Bailey 
observed this important result while trying to clarify Rogers 
second proof of Rogers-Ramanujan (RR) identities (1917).  The 
first and second RR identities are
 \begin{equation}\label{eq:RR_id1}
\sum_{n=0}^{\infty}\frac{q^{n^2}}{(q)_n}=\frac{1}{(q)_{\infty}}\sum_{n=-\infty}^{\infty}
               (q^{n(10n+1)}-q^{(5n+2)(2n+1)})
                         =\prod_{n=1}^{\infty} \frac{1}{(1-q^{5n-1})(1-q^{5n-4})}
\end{equation}    
\begin{equation}\label{eq:RR_id2}
\sum_{n=0}^{\infty}\frac{q^{n(n+1)}}{(q)_n}=\frac{1}{(q)_{\infty}}\sum_{n=-\infty}^{\infty}
               (q^{n(10n+3)}-q^{(5n+1)(2n+1)})
             =\prod_{n=1}^{\infty} \frac{1}{(1-q^{5n-2})(1-q^{5n-3})}
\end{equation} 
There are many different proofs of these identities for example  
\cite{Mac:1916,R:1894,R:1917,
RR:1919, Schur:1917}  and there are also many generalizations 
\cite{A:1974,B:1980,G:1961,Sl:1951,Sl:1952}
in the theory of partitions. It is worth mentioning that the equality of the first 
two expressions 
of \eqref{eq:RR_id1} is the Bose-Fermi identity  for the minimal model $M(2,5)$. Therefore,
Bose-Fermi identities can be interpreted  as a generalization of the RR-identities.
 
Bailey's lemma has been a very useful method for proving RR-type $q$-identities.  We
use Bailey's lemma to prove our results in Chapter~\ref{chap:bailey}. Let us 
state the lemma here. 
A pair $(\alpha_n,\beta_n)$ of sequences $\{\alpha_n\}_{n\ge 0}$ and 
$\{\beta_n\}_{n\ge 0}$ is called a \textbf{Bailey pair} with respect to $a$ if 
\begin{equation}\label{eq:baileypair}
\beta_n=\sum_{j=0}^n \frac{\alpha_j}{(q)_{n-j}(aq)_{n+j}}
\end{equation}
where 
\begin{equation*}
\begin{split}
(a)_n:=(a;q)_n &= \prod_{k=0}^{n-1}(1-aq^k),\\
\end{split}
\end{equation*}
\begin{theorem}(\textbf{Bailey's lemma}) If $(\alpha_n,\beta_n)$ is a Bailey pair
with respect to $a$ then for two parameters $\rho_1,\rho_2$,
\begin{equation}\label{eq:baileylemma}
\begin{split}
\sum_{n=0}^{\infty}(\rho_1)_n(\rho_2)_n & (aq/\rho_1\rho_2)^n\beta_n\\
= & \frac{(aq/\rho_1)_{\infty}(aq/\rho_2)_{\infty}}{(aq)_{\infty}(aq/\rho_1 
\rho_2)_{\infty}}\sum_{n=0}^{\infty} \frac{(\rho_1)_n(\rho_2)_n
(aq/\rho_1\rho_2)^n\alpha_n}{(aq/\rho_1)_n(aq/\rho_2)_n}.
\end{split}
\end{equation}
\end{theorem}
Putting different Bailey pairs in this lemma many  RR type identities
were proved by Rogers ~\cite{R:1894,R:1917}, Bailey~\cite{Bailey:1949}
and Slater~\cite{Sl:1951,Sl:1952}. In their proofs they considered the following two 
specializations of the parameters:
\begin{equation}\label{eq:infinity_spec}
\rho_1\longrightarrow \infty,\quad \rho_2\longrightarrow \infty
\end{equation}
\begin{equation}\label{eq:half_finite_spec}
\rho_1\longrightarrow  \infty,\quad \rho_2=\text{finite}. 
\end{equation}   

To get a feel of how Bailey's lemma can be used to prove RR-type identities let 
us discuss here briefly the proof of RR identities \eqref{eq:RR_id1} and 
\eqref{eq:RR_id2}. We use the specialization \eqref{eq:infinity_spec}. Let 
$\rho_1\longrightarrow \infty,\quad \rho_2\longrightarrow \infty$
in the Bailey lemma. We obtain
\begin{equation} \label{eq:sp_baileylemma}
\sum_{n=0}^{\infty} a^n q^{n^2} \beta_n = \frac{1}{(aq)_{\infty} }\sum_{n=0}^{\infty} 
a^nq^{n^2}\alpha_n
\end{equation}

The Bailey pair used to prove \eqref{eq:RR_id1} is given by
\begin{equation*}
\begin{split}
\alpha_0 &=1\\
\alpha_n &=(-1)^nq^{n(3n-1)/2}(1+q^n), \quad n\ge 1\\
\beta_n &= \frac{1}{(q)_n} \quad n\ge 0
\end{split}
\end{equation*}
Inputing this Bailey pair in \eqref{eq:sp_baileylemma} with $a=1$ and 
simplifying we find,
\begin{equation} \label{eq:sp_RR}
\sum_{n=0}^{\infty} \frac{q^{n^2}} {(q)_n} = \frac{1}{(q)_{\infty} }
       \sum_{n=-\infty}^{\infty} (-1)^n q^{n(5n-1)/2}
\end{equation}
\noindent
Note the left-hand side is exactly the left-hand side of \eqref{eq:RR_id1}. To prove
that the right hand side of \eqref{eq:sp_RR} equals the right hand 
side of \eqref{eq:RR_id1} we use Jacobi's triple product identity,
\begin{equation}\label{eq:jacobi}
\sum_{n=-\infty}^{n=\infty}q^{k^2}z^n=(q^2,-qz, -q/z;q^2)_{\infty},
\end{equation}
where $(a_1,a_2,a_3;q^2)_{\infty}=(a_1;q^2)_{\infty}(a_2;q^2)_{\infty}
(a_3;q^2)_{\infty}$.

Using  \eqref{eq:jacobi} we can rewite right-hand side of \eqref{eq:sp_RR} as
\begin{equation*}
\begin{split}
\frac{1}{(q)_{\infty} }\sum_{n=-\infty}^{\infty} (-1)^n q^{n(5n-1)/2}&=\frac{1}{(q)_{\infty} }
\sum_{n=-\infty}^{\infty}  (-q^{-1/2})^n(q^{5/2})^{n^2}\\
&=\frac{1}{(q)_{\infty} }\times (q^5;q^5)_{\infty} \times (q^2;q^5)_{\infty}\times (q^3;q^5)_{\infty}\\
&=\prod_{n=1}^{\infty} \frac{1}{(1-q^{5n-1})(1-q^{5n-4})}
\end{split}
\end{equation*}
\noindent
which is the right hand side of \eqref{eq:RR_id1}.  Similarly, the second 
RR identity \eqref{eq:RR_id2} can be 
proved using the Bailey pair
\begin{equation*}
\begin{split}
\alpha_n &=(-1)^n\frac{q^{n(3n+1)/2} (1-q^{2n+1})}{1-q}, \quad n\ge 0,\\
\beta_n &= \frac{1}{(q)_n} \quad n\ge 0.
\end{split}
\end{equation*}

Slater \cite{Sl:1951,Sl:1952} used this lemma extensively in order to prove 130 different
RR-type identities.  In connecting the Bailey construction to physics, the most 
remarkable step was achieved when Foda and Quano~\cite{FQ:1995,FQ:1996}  derived
identities for the Virasoro characters using Bailey's lemma. Their method is
a constructive procedure that starts from a polynomial generalization of
a Bose-Fermi  identity for one  CFT and then produces a  Bose-Fermi identity
for the character of another CFT. This is known as the Bailey flow.  Hence, new fermionic 
formulas for CFTs can be calculated via Bailey flow from known fermionic formulas
of another CFT. The Bailey flow from $M(p-1,p)$ to $M(p,p+1)$ is
presented in \cite{BMS:1995,FQ:1996}, and further flows to some 
special $N=1$ and $N=2$ supersymmetric models are given in \cite{BMS:1995}.
This led us to investigate Bailey flows from $M(p,p')$ with $p,p'$
arbitrary coprime positive integers to other CFTs. In \cite{BMS:1995}
it is conjectured that their methods, which was applied to the unitary case when $p=p'-1$
can also be applied to the general case. This is the problem we study  in 
Chapter~\ref{chap:bailey}. We demonstrate new Bailey flows from $M(p,p')$ to $N=1$
and $N=2$ superconformal algebras and prove the conjectures of \cite{BMS:1995}.
We present new Bose-Fermi identities for the characters of $N=1$ and $N=2$
superconformal algebras. These new identities can be thought of as the generalized RR-type 
identities for the $N=1$ and $N=2$ superconformal characters.

    %
    %

    \newchap{Fermionic formulas  for unrestricted Kostka polynomials}
    \label{chap:kostka}
        \section{Introduction}
The Kostka numbers $K_{\lambda\mu}$, indexed by the two partitions $\lambda$ and
$\mu$, play an important role in symmetric function theory, representation theory, 
combinatorics, invariant theory and mathematical physics. The Kostka polynomials 
$K_{\lambda\mu}(q)$ are $q$-analogs of the Kostka numbers. There are several
combinatorial definitions of the Kostka polynomials. For example Lascoux and 
Sch{\"u}tzenberger~\cite{LS:1978} proved that the Kostka polynomials are generating
functions of semi-standard tableaux of shape $\lambda$ and content $\mu$
with charge statistic. In~\cite{NY:1997} the Kostka polynomials are expressed
as generating function over highest-weight crystal paths with energy statistics. Crystal
paths are elements in tensor products of finite-dimensional crystals.
Dropping the highest-weight condition yields unrestricted Kostka 
polynomials~\cite{HKOTT:1999,HKKOTY:1999,HKOTY:1999,SW:1999}. In the
$A_1^{(1)}$ setting, unrestricted Kostka polynomials or $q$-supernomial
coefficients were introduced in~\cite{SW:1997} as $q$-analogs of the coefficient of 
$x^a$ in the expansion of $\prod_{j=1}^N (1+x+x^2+\cdots+x^j)^{L_j}$. 
An explicit formula for the $A_{n-1}^{(1)}$ unrestricted Kostka polynomials for 
completely symmetric and completely antisymmetric crystals was proved 
in~\cite{HKKOTY:1999,Kir:1998}. This formula is called fermionic as it is a 
manifestly positive expression.  

In this chapter we give a new explicit fermionic formula for the unrestricted Kostka 
polynomials for all Kirillov--Reshetikhin crystals of type $A_{n-1}^{(1)}$.
This fermionic formula can be naturally interpreted in terms of a new set
of unrestricted rigged configurations for type $A_{n-1}^{(1)}$.
Rigged configurations are combinatorial objects originating from the Bethe
Ansatz, that label solutions of the Bethe equations. The simplest version of
rigged configurations appeared in Bethe's original paper~\cite{Bethe:1931}
and was later generalized by Kerov, Kirillov and Reshetikhin~\cite{KKR:1986,KR:1988} 
to models with $\mathrm{GL}(n)$ symmetry.
Since the solutions of the Bethe equations label highest weight vectors, one expects 
a bijection between rigged configurations and semi-standard Young tableaux in 
the $\mathrm{GL}(n)$ case. Such a bijection was given in~\cite{KR:1988,KSS:2002}. 
Here we extend this bijection to a bijection $\Phi$ between the new set of unrestricted
rigged configurations and unrestricted paths. It should be noted that $\Phi$ is 
defined algorithmically. In~\cite{Sch:2005b} the bijection was established in a different 
manner by constructing a crystal structure on the set of rigged configurations.
Here we show that the crystal structures are compatible under the algorithmically
defined $\Phi$ and use this to prove that $\Phi$ preserves the statistics.

The bijection $\Phi$ has been implemented as a C++ program 
and has been incorporated into the combinatorics package of MuPAD-Combinat by
Francois Descouens~\cite{MuPAD:2005}.  The program is given in 
chapter~\ref{chap:program}.

This chapter is structured as follows. In Section~\ref{sec:paths} we review crystals
of type $A_{n-1}^{(1)}$, highest weight paths, unrestricted paths and the definition 
of generalized Kostka polynomials and  unrestricted Kostka 
polynomials as generating functions of highest weight paths and unrestricted paths
respectively with energy statistics.
In Section~\ref{sec:RC} we give our new definition of unrestricted rigged configurations
(see Definition~\ref{def:uRC}) and derive from this a fermionic expression for the 
generating function of unrestricted rigged configurations graded by cocharge (see
Section~\ref{subsec:fermionic}). The statistic preserving bijection between unrestricted 
paths and unrestricted rigged configurations is established in Section~\ref{sec:main} 
(see Definition~\ref{def:bij} and Theorem~\ref{thm:bij}). As a corollary this
yields the equality of the unrestricted Kostka polynomials and the fermionic
formula of Section~\ref{sec:RC} (see Corollary~\ref{cor:X=M}). The result that
the crystal structures on paths and rigged configurations are compatible under
$\Phi$ is stated in Theorem~\ref{thm:commute}. Most of the technical proofs are relegated 
to the last  three sections. An extended abstract of this chapter can be found 
in~\cite{DS:2005}.

\section{Unrestricted paths and Kostka polynomials} \label{sec:paths}

\subsection{Crystals $B^{r,s}$ of type $A_{n-1}^{(1)}$}
Kashiwara~\cite{Kash:1990}  introduced the notion of crystals and crystal graphs 
as a combinatorial means to study representations of quantum algebras associated 
with any symmetrizable Kac--Moody algebra. In this paper we only consider the 
Kirillov--Reshetikhin crystal $B^{r,s}$ of type $A_{n-1}^{(1)}$ and hence restrict
to this case here.

As a set, the crystal $B^{r,s}$ consists of all column-strict Young tableaux of shape 
$(s^r)$ over the alphabet $\{1,2,\ldots,n\}$. As a crystal associated to the
underlying algebra of finite type $A_{n-1}$, $B^{r,s}$ is isomorphic to the highest
weight crystal with highest weight $(s^r)$. We will define the classical crystal
operators explicitly here. The affine crystal operators $e_0$ and $f_0$ are given
explicitly in~\cite{Sh:2002}. Since we do not use these operators we will omit the details.

Let $I=\{1,2,\ldots,n-1\}$ be the index set for the vertices of the Dynkin
diagram of type $A_{n-1}$, $P$ the weight lattice, $\{\La_i\in P \mid i\in I\}$
the fundamental roots, $\{\alpha_i\in P \mid i\in I\}$ the 
simple roots, and $\{h_i\in \Hom_{\Z}(P,\Z) \mid i\in I\}$ the simple coroots.  
As a type $A_{n-1}$ crystal, $B=B^{r,s}$ is equipped with maps 
$e_i,f_i: B\longrightarrow B\cup\{0\}$ and  $\wt: B\longrightarrow P$ 
for all $ i\in I$ satisfying 
\begin{equation*}
\begin{split}
 &f_i(b)=b' \Leftrightarrow e_i(b')=b \text{ if } b,b'\in B\\
 &\wt(f_i(b))=\wt(b)-\alpha_i \text{ if } f_i(b)\in B\\
 &\langle h_i,\wt(b) \rangle=\varphi_i(b)-\varepsilon_i(b),
\end{split}
\end{equation*}
where $\langle \cdot,\cdot \rangle$ is the natural pairing. The maps $f_i,e_i$ are known as
the Kashiwara operators. Here for $b\in B$,
\begin{align*}
 \varepsilon_i(b)&=\max\{k\ge 0 \mid e_i^k(b)\ne 0\}\\
 \varphi_i(b)&=\max\{k\ge 0 \mid f_i^k(b)\ne 0\}.
\end{align*}
Note that for type $A_{n-1}$, $P=\Z^{n}$ and $ \alpha_i=\epsilon_i-\epsilon_{i+1}$
where $\{\epsilon_i \mid i\in I\}$ is the standard basis in $P$. Here $\wt(b)=
(w_1,\ldots,w_n)$ is the weight of $b$ where $w_i$ counts the number of letters 
$i$ in $b$. 

Following~\cite{KN:1994} let us give the action of $e_i$ and $f_i$ 
for $i\in I$.  Let $b \in B^{r,s}$ be a tableau of shape $(s^r)$. The row word of $b$ is 
defined by $\word(b)=w_r\cdots w_2w_1$ where $w_k$ is the word obtained by reading the 
$k$-th row of $b$ from left to right. 
To find $f_i(b)$ and $e_i(b)$ we only consider the subword consisting of the letters $i$ 
and $i+1$ in the word of $b$. First view each $i+1$ in the subword as an opening bracket and
each $i$ as a closing bracket. Then we ignore each adjacent pair of matched 
brackets  successively. At the end of this process we are left with a subword of the
form $i^p(i+1)^q$. If $p>0$ (resp. $q>0$) then $f_i(b)$ (resp. $e_i(b)$) 
is obtained from $b$ by replacing the unmatched subword $i^p(i+1)^q$ 
by $i^{p-1}(i+1)^{q+1}$ (resp. $i^{p+1}(i+1)^{q-1}$). If $p=0$ (resp. $q=0$) then $f_i(b)$ 
(resp. $e_i(b)$) is undefined and we write $f_i(b)=0$ (resp. $e_i(b)=0$). 
 
 A crystal $B$ can be viewed as a directed edge-colored graph. The vertices of the
 graph are the elements of $B$ and there is a directed edge between vertices $b$
 and $b'$ labeled $i\in I\cup \{0\}$, if and only if $f_i(b)=b'$. This directed graph is known as
 the {\bf crystal graph}.   
 \begin{example}
 The crystal graph for $B=B^{1,1}$ is given in Figure~\ref{fig:1}.

\begin{figure}[h]
\center{
\scalebox{.8}{
\begin{picture}(250,62)(-10,-12)
\BText(0,0){1} \LongArrow(10,0)(40,0) \BText(50,0){2}
\LongArrow(60,0)(90,0) \BText(100,0){3} \LongArrow(110,0)(140,0)
\Text(160,0)[]{$\cdots$} \LongArrow(175,0)(205,0)
\BText(220,0){n} \LongArrowArc(110,-181)(216,62,118)
\PText(25,2)(0)[b]{1} \PText(75,2)(0)[b]{2} \PText(125,2)(0)[b]{3}
\PText(190,2)(0)[b]{n-1} \PText(110,38)(0)[b]{0}
\end{picture}}
\caption{Crystal $B^{1,1}$.\label{fig:1}}
}
\end{figure}
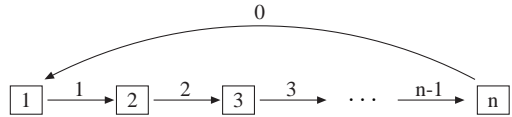
\end{example}
 
Given two crystals $B$ and $B'$, we can also define a new crystal by taking the tensor 
product $B\otimes B'$. As a set $B\otimes B'$ is just the Cartesian product of the
sets $B$ and $B'$. The weight function $\wt$ for $b\otimes b' \in B\otimes B'$
is $\wt(b\otimes b')=\wt(b)+\wt(b')$ and the Kashiwara operators $e_i$, $f_i$ are defined 
as follows
\begin{equation*}
\begin{split}
e_i(b\otimes b')&=
\begin{cases} e_ib\otimes b' & \text{if $\varepsilon_i(b)>\varphi_i(b')$,}\\
              b\otimes e_i b' & \text{otherwise,}\end{cases}\\
f_i(b\otimes b')&=
\begin{cases} f_ib\otimes b' & \text{if $\varepsilon_i(b)\ge \varphi_i(b')$,}\\
              b\otimes f_i b' & \text{otherwise.}\end{cases}
\end{split}
\end{equation*}
This action of $f_i$ and $e_i$ on the tensor product is compatible with the 
previously defined action on $\word(b\otimes b')=\word(b)\word(b')$.

\begin{example}
Let $i=2$ and  
\begin{equation*}
b\;=\; \young(12,23) \otimes \young(23,34,45).
\end{equation*}
Then $\word(b)=2312453423$, the relevant subword is $23-2--3-23$, and the unmatched 
subword is $2--------3$. Hence
\begin{equation*}
f_2(b)\;=\;\young(12,33) \otimes \young(23,34,45)\quad \text{and} \quad
e_2(b)\;=\;\young(12,23) \otimes \young(22,34,45). 
\end{equation*}
\end{example} 

\subsection{Paths and unrestricted paths}
Let $B=B^{r_k,s_k}\otimes B^{r_{k-1},s_{k-1}}\otimes \cdots \otimes B^{r_1,s_1}$.

A \textbf{highest weight path} or simply \textbf{path} is an element $b \in B$ 
 such that $e_i(b)=0$ for all $1\le i \le n-1$. It is known that the
 weight vector of a highest weight element for type $A_{n-1}$ is
 a partition with at most $n$ nonzero parts. Let $\la=(\la_1,\la_2,\ldots,\la_n)$ 
 be a partition with at most $n$ nonzero parts, then the set of 
all highest weight paths of weight $\la$ and shape $B$ is defined as
 \begin{equation*}
\overline{\Path}(B,\la)=\{b\in B\mid \wt(b)=\la \quad \text{and $e_i(b)=0$ for all $1\le i\le n-1$}  \}.
\end{equation*}

An unrestricted path is an element in the tensor product of crystals
$B=B^{r_k,s_k}\otimes B^{r_{k-1},s_{k-1}}\otimes \cdots \otimes B^{r_1,s_1}$.
Let $\la=(\la_1,\la_2,\ldots,\la_n)$ be an $n$-tuple of nonnegative integers. The set 
of \textbf{unrestricted paths} is defined as
\begin{equation*}
\Path(B,\la)=\{b\in B\mid \wt(b)=\la\}.
\end{equation*}
Note that the weight of an unrestricted path need not be a partition.

\begin{example}
For $B=B^{1,1}\otimes B^{2,2}\otimes B^{3,1}$ of type $A_3$  
the path 
\begin{equation*}
b\;=\; \young(1) \otimes \young(13,24) \otimes \young(1,2,3)
\end{equation*}
is a highest weight path of weight $\la=(2,2,2,1)$.
The path
\begin{equation*}
b\;=\; \young(2) \otimes \young(12,24) \otimes \young(1,3,4)
\end{equation*}
is an unrestricted path of weight $\la=(2,3,1,2)$.
\end{example}  

There exists a crystal isomorphism $R:B^{r,s}\otimes B^{r',s'} \to B^{r',s'} \otimes B^{r,s}$,
called the \textbf{combinatorial $R$-matrix}. Combinatorially it is given as follows. 
Let $b\in B^{r,s}$ and $b'\in B^{r',s'}$.
The product $b\cdot b'$ of two tableaux is defined as the Schensted insertion of $b'$ into
$b$. Then $R(b\otimes b')=\tilde{b}'\otimes \tilde{b}$ is the unique pair of tableaux
such that $b\cdot b'=\tilde{b}'\cdot\tilde{b}$.

The \textbf{local energy function} $H:B^{r,s}\otimes B^{r',s'}\to \Z$ is defined as
follows. For $b\otimes b'\in B^{r,s}\otimes B^{r',s'}$, $H(b\otimes b')$ is the number
of boxes of the shape of $b\cdot b'$ outside the shape obtained by concatenating 
$(s^r)$ and $({s'}^{r'})$.

\begin{example}
For 
\begin{equation*}
b\otimes b'= \young(12,24) \otimes \young(1,3,4)
\end{equation*}
we have 
\begin{equation*}
b\cdot b' = \young(113,224,4) = \young(1,2,4) \cdot \young(13,24) = \tilde{b}'\cdot\tilde{b}.
\end{equation*}
so that
\begin{equation*}
R(b\otimes b')=\tilde{b}'\otimes\tilde{b}=\young(1,2,4) \otimes \young(13,24).
\end{equation*}
Since the concatentation of $\yng(2,2)$ and $\yng(1,1,1)$ is $\yng(3,3,1)$, the local
energy function $H(b\otimes b')=0$.
\end{example}

Now let $B=B^{r_k,s_k}\otimes \cdots\otimes B^{r_1,s_1}$ be a $k$-fold tensor
product of crystals. Let $b=b_k\otimes b_{k-1}\otimes \cdots \otimes b_1 \in B$.
The \textbf{tail energy function} $\Dt:B\to \Z$ is given by
\begin{equation*}
  \Dt(b) = \sum_{1\le i<j\le k} H_{j-1} R_{j-2} \dotsm R_{i+1} R_i(b),
\end{equation*}
where $H_i$ (resp. $R_i$) is the local energy function (resp. combinatorial $R$-matrix)
acting on the $i$-th and $(i+1)$-th tensor factors of $b\in B$ from the right.

\begin{definition}
The generalized Kostka polynomial $\tilde{K}_{\la,\mu}(q)$ with $\mu=(\mu_k,\cdots,\mu_1)$ 
where $\mu_j$ is a rectangular partition of height $r_j$ and width $s_j$, is the generating
function of highest weight paths with the energy function
$$\tilde{K}_{\la,\mu}(q)=\sum_{b\in \overline{P}(\la, B)}q^{\Dt(b)}.$$ 
\end{definition}
We should note that the generalized Kostka polynomial defined here 
corresponds to the cocharge Kostka polynomial where cocharge \cite{Macdonald,SW:1999} 
is a new statistic on the set of tableaux which is a shift of the charge statistic. 

$A_{n-1}^{(1)}$-unrestricted Kostka polynomials or supernomial coefficients were first 
introduced in~\cite{SW:1999} as generating functions of unrestricted paths graded by an 
energy function. 
\begin{definition}
The $q$-\textbf{supernomial coefficient} or the \textbf{unrestricted Kostka polynomial} is 
defined as
\begin{equation*}
X(B,\la)=\sum_{b\in \Path(B,\la)} q^{\Dt(b)}.
\end{equation*}
\end{definition}

\section{Unrestricted rigged configurations and fermionic formula}
\label{sec:RC}
Rigged configurations are combinatorial objects invented to label
the solutions of the Bethe equations, which give the eigenvalues of
the Hamiltonian of the underlying physical model \cite{Bethe:1931}.
Motivated by the fact that representation theoretically the eigenvectors
and eigenvalues can also be labelled by Young tableaux, Kirillov
and Reshetikhin~\cite{KR:1988} gave a bijection between tableaux
and rigged configurations. This result and generalizations thereof
were proven in~\cite{KSS:2002}.

In terms of crystal base theory, the bijection is between highest
weight paths and rigged configurations. The new result of this paper
is an extension of this bijection to a bijection between unrestricted paths
and a new set of rigged configurations. The new set of unrestricted rigged
configurations is defined in this section, whereas the bijection is given in 
section~\ref{sec:main}.
In~\cite{Sch:2005b}, a crystal structure on the new set of unrestricted rigged 
configurations is given, which provides a different description of the bijection.

\subsection{Unrestricted rigged configurations}
Let $B=B^{r_k,s_k}\otimes \cdots \otimes B^{r_1,s_1}$ and denote by
$L=(L_i^{(a)}\mid (a,i)\in \HH)$ the multiplicity array of $B$, where
$L_i^{(a)}$ is the multiplicity of $B^{a,i}$ in $B$. Here 
$\HH=I \times \Z_{>0}$ and $I=\{1,2,\ldots,n-1\}$ is the index set of the
Dynkin diagram $A_{n-1}$.
The sequence of partitions $\nu=\{\nu^{(a)}\mid a\in I \}$ is a
\textbf{$(L,\la)$-configuration} if
\begin{equation}\label{eq:size1}
\sum_{(a,i)\in\HH} i m_i^{(a)} \alpha_a = \sum_{(a,i)\in\HH} i
L_i^{(a)} \La_a- \la,
\end{equation}
where $m_i^{(a)}$ is the number of parts of length $i$ in partition
$\nu^{(a)}$. Note that we do not require $\la$ to be a dominant weight here.
The \textbf{(quasi-)vacancy number} of a configuration is defined as
\begin{equation*}
p_i^{(a)}=\sum_{j\ge 1} \min(i,j) L_j^{(a)}
 - \sum_{(b,j)\in \HH} (\alpha_a | \alpha_b) \min(i,j)m_j^{(b)}.
\end{equation*}
Here $(\cdot | \cdot )$ is the normalized invariant form on the weight lattice $P$
such that $(\alpha_i | \alpha_j)$ is the Cartan matrix. Let
$\Conf(L,\la)$ be the set of all $(L,\la)$-configurations.
We call $p_i^{(a)}$ quasi-vacancy number to indicate that they can actually
be negative in our setting. For the rest of the paper we will simply call
them vacancy numbers. 

When the dependence of $m_i^{(a)}$ and $p_i^{(a)}$ on the configuration
$\nu$ is crucial, we also write $m_i^{(a)}(\nu)$ and $p_i^{(a)}(\nu)$,
respectively.

In the usual setting a rigged configuration $(\nu,J)$ consists of a configuration
$\nu\in \Conf(L,\la)$ together with a double sequence of partitions
$J=\{J^{(a,i)}\mid (a,i)\in\HH \}$ such that the partition
$J^{(a,i)}$ is contained in a $m_i^{(a)}\times p_i^{(a)}$ rectangle.
In particular this requires that $p_i^{(a)}\ge 0$. For unrestricted paths
we need a bigger set, where the lower bound on the parts in $J^{(a,i)}$
can be less than zero.

To define the lower bounds we need the following notation. Let 
$\la'=(c_1,c_2,\ldots,c_{n-1})^t$ where $c_k=\la_{k+1}+\la_{k+2}+\cdots+\la_n$.
We also set $c_0=c_1$.
Let $\A(\la')$ be the set of tableaux of shape $\la'$ such that the entries 
in column $k$ are from the set $\{1,2,\ldots,c_{k-1}\}$ and are strictly decreasing along 
each column.

\begin{example}  \label{ex:lbtab_example} 
For $n=4$ and $\la=(0,1,1,1)$, the set $\A(\la')$
consists of the following tableaux
\begin{equation*}
\young(332,22,1) \quad \young(332,21,1) \quad \young(322,21,1) \quad 
\young(331,22,1) \quad \young(331,21,1) \quad \young(321,21,1).
\end{equation*}
\end{example}
Note that each $t\in \A(\la')$ is weakly decreasing along each row. This is
due to the fact that $t_{j,k}\ge c_k-j+1$ since column $k$ of height $c_k$ is 
strictly decreasing and $c_k-j+1 \ge t_{j,k+1}$ since the entries in column $k+1$ 
are from the set $\{1,2,\ldots,c_{k}\}$.

Given $t\in\A(\la')$, we define the \textbf{lower bound} as
\begin{equation*}
M_i^{(a)}(t)=-\sum_{j=1}^{c_a} \chi(i\ge t_{j,a})
+\sum_{j=1}^{c_{a+1}} \chi(i\ge t_{j,a+1}),
\end{equation*}
where $t_{j,a}$ denotes the entry in row $j$ and column $a$ of $t$,
and $\chi(S)=1$ if the the statement $S$ is true and $\chi(S)=0$ otherwise.

\begin{example}
For the tableau $t=\young(332,22,1) $ from example~\ref{ex:lbtab_example} 
some of the lower bounds are given by 
$$M_1^{(1)}(t)=-1, M_4^{(1)}(t)=-1, M_1^{(2)}(t)=0, M_3^{(2)}(t)=-1, M_2^{(3)}(t)=-1.$$
\end{example}

Let $M,p,m\in \Z$ such that $m\ge 0$.
A $(M,p,m)$-quasipartition $\mu$ is a tuple of integers $\mu=(\mu_1,\mu_2,\ldots,\mu_m)$
such that $M\le \mu_m\le \mu_{m-1}\le \cdots\le \mu_1\le p$. Each $\mu_i$ is called
a part of $\mu$. Note that for $M=0$ this would be a partition with at most $m$ parts each 
not exceeding $p$.

\begin{definition} \label{def:uRC}
An \textbf{unrestricted rigged configuration} $(\nu,J)$ associated to a multiplicity
array $L$ and weight $\la$ is a configuration
$\nu\in\Conf(L,\la)$ together with a sequence $J=\{J^{(a,i)}\mid (a,i)\in\HH\}$
where $J^{(a,i)}$ is a $(M_i^{(a)}(t),p_i^{(a)},m_i^{(a)})$-quasipartition
for some $t\in \A(\la')$. Denote the set of all unrestricted rigged configurations
corresponding to $(L,\la)$ by $\RC(L,\la)$.
\end{definition}

\begin{remark}\mbox{}
\begin{enumerate}
\item
Note that this definition is similar to the definition of level-restricted rigged 
configurations~\cite[Definition 5.5]{SS:2001}. Whereas for level-restricted
rigged configurations the vacancy number had to be modified according to 
tableaux in a certain set, here the lower bounds are modified.
\item 
For type $A_1$ we have $\la=(\la_1,\la_2)$ so that $\A=\{ t\}$
contains just the single tableau
\begin{equation*}
t=\begin{array}{|c|} \hline \la_2\\ \hline \la_2-1\\ \hline \vdots\\ \hline 1\\
\hline \end{array}.
\end{equation*}
In this case $M_i(t)=-\sum_{j=1}^{\la_2} \chi(i\ge t_{j,1})=-i$. This agrees
with the findings of~\cite{T:2004}.
\end{enumerate}
\end{remark}
The quasipartition $J^{(a,i)}$ is called \textbf{singular} if it has a
part of size $p_i^{(a)}$. 
It is often useful to view an (unrestricted) rigged configuration $(\nu,J)$ as a
sequence of partitions $\nu$ where the parts of size $i$ in
$\nu^{(a)}$ are labeled by the parts of $J^{(a,i)}$. The pair
$(i,x)$ where $i$ is a part of $\nu^{(a)}$ and $x$ is a part of
$J^{(a,i)}$ is called a \textbf{string} of the $a$-th rigged
partition $(\nu,J)^{(a)}$. The label $x$ is called a \textbf{rigging}.

\begin{example}
Let $n=4$, $\la=(2,2,1,1)$, $L_1^{(1)}=6$ and all other $L_i^{(a)}=0$. Then
\begin{equation*}
(\nu,J) \;=\; \yngrc(3,-2,1,0) \quad \yngrc(2,0) \quad \yngrc(1,-1)
\end{equation*}
is an unrestricted rigged configuration in $\RC(L,\la)$, where we have written
the parts of $J^{(a,i)}$ next to the parts of length $i$ in partition $\nu^{(a)}$.
To see that the riggings form quasipartitions, let us write
the vacancy numbers $p_i^{(a)}$ next to the parts of length $i$ in partition $\nu^{(a)}$:
\begin{equation*}
\yngrc(3,0,1,3) \quad \yngrc(2,0) \quad \yngrc(1,-1).
\end{equation*}
This shows that the labels are indeed all weakly below the vacancy numbers. For
\begin{equation*}
\young(441,33,2,1) \in \A(\la')
\end{equation*}
we get the lower bounds
\begin{equation*}
\yngrc(3,-2,1,-1) \quad \yngrc(2,0) \quad \yngrc(1,-1),
\end{equation*}
which are less or equal to the riggings in $(\nu,J)$.
\end{example}

Let $B=B^{r_k,s_k}\otimes \cdots \otimes B^{r_1,s_1}$ and $L$
the corresponding multiplicity array. Let $(\nu,J)\in\RC(L,\la)$.
Note that rewritting \eqref{eq:size1} we get
\begin{equation}\label{eq:size}
|\nu^{(a)}| = \sum_{j>a} \la_j -\sum_{j=1}^k s_j \max (r_j-a,0).
\end{equation}
Hence for large $i$, by definition of vacancy numbers we have
\begin{equation}\label{eq:p limit}
\begin{split}
p_i^{(a)}&=|\nu^{(a-1)}|-2|\nu^{(a)}|+|\nu^{(a+1)}|+\sum_j \min(i,j)L_j^{(a)}\\
&=\la_a-\la_{a+1}
\end{split}
\end{equation} 
and
\begin{equation}\label{eq:M limit}
\begin{split}
M_i^{(a)}(t)&=-\sum_{j=1}^{c_a} \chi(i\ge t_{j,a})
+\sum_{j=1}^{c_{a+1}} \chi(i\ge t_{j,a+1})\\
&=-c_a+c_{a+1}=-\la_{a+1}.
\end{split}
\end{equation}
For a given $t\in\A(\la')$ define
\begin{equation*}
\Dp_i^{(a)}(t)=p_i^{(a)}-M_i^{(a)}(t).
\end{equation*}
We write $\Dp_i^{(a)}$ for $\Dp_i^{(a)}(t)$ when there is no cause of confusion.
For large $i$, $\Dp_i^{(a)}(t)=\la_a$.

{}From the definition of $p_i^{(a)}$ one may easily verify that
\begin{equation}\label{eq:p ineq}
-p_{i-1}^{(a)}+2p_i^{(a)}-p_{i+1}^{(a)}\ge m_i^{(a-1)}-2m_i^{(a)}+m_i^{(a+1)}.
\end{equation}
Let $t_{\cdot,a}$ denote the $a$-th column of $t$. Then it follows from the
definition of $M_i^{(a)}(t)$ that
\begin{equation*}
M_i^{(a)}(t)=M_{i-1}^{(a)}(t)-\chi(i\in t_{\cdot,a})+\chi(i\in t_{\cdot,a+1}).
\end{equation*}
Hence \eqref{eq:p ineq} can be rewritten as 
\begin{multline}\label{eq:dp ineq}
-\Dp_{i-1}^{(a)}+2\Dp_i^{(a)}-\Dp_{i+1}^{(a)}
-\chi(i\in t_{\cdot,a})+\chi(i\in t_{\cdot,a+1})\\
+\chi(i+1\in t_{\cdot,a})-\chi(i+1\in t_{\cdot,a+1})
\ge m_i^{(a-1)}-2m_i^{(a)}+m_i^{(a+1)}.
\end{multline}

\begin{lemma}\label{lem:convex}
Suppose that for some $t\in \A(\la')$, $\Dp_i^{(a)}(t)\ge 0$ for all $a\in I$ 
and $i$ such that $m_i^{(a)}>0$. Then there exists a $t'\in\A(\la')$ such that
$\Dp_i^{(a)}(t')\ge 0$ for all $i$ and $a$.
\end{lemma}
\begin{proof}
By definition $\Dp_0^{(a)}(t)=0$ and $\Dp_i^{(a)}(t)=\la_a\ge 0$ for large $i$.
By~\eqref{eq:dp ineq}
\begin{multline}\label{eq:convex}
\Dp_i^{(a)}(t)\ge \frac{1}{2}\bigl\{\Dp_{i-1}^{(a)}(t)+\Dp_{i+1}^{(a)}(t)
+\chi(i\in t_{\cdot,a})-\chi(i\in t_{\cdot,a+1})\\
-\chi(i+1\in t_{\cdot,a})+\chi(i+1\in t_{\cdot,a+1})
+m_i^{(a-1)}+m_i^{(a+1)}\bigr\}
\end{multline}
when $m_i^{(a)}=0$. Hence $\Dp_i^{(a)}(t)<0$ is only possible if
$m_i^{(a-1)}=m_i^{(a+1)}=0$, column $a$ of $t$ contains $i+1$ but no $i$,
and column $a+1$ of $t$ contains $i$ but no $i+1$.
Let $k$ be minimal such that $\Dp_i^{(k)}(t)<0$. Note that $k>1$ since
the first column of $t$ contains all letters $1,2,\ldots, c_1$. 
Let $k'\le k$ be minimal such that $\Dp_i^{(a)}(t)=0$ for all $k'\le a<k$.
Then inductively for $a=k-1,k-2,\ldots,k'$ it follows from~\eqref{eq:convex}
that $m_i^{(a-1)}=0$ and column $a$ of $t$ contains $i+1$ but no $i$. Construct a 
new $t'$ from $t$ by replacing all letters $i+1$ in columns $k',k'+1,\ldots,k$ by $i$.
This accomplishes that $\Dp_j^{(a)}(t')\ge 0$ for all $j$ and $1\le a<k$,
$\Dp_i^{(k)}(t')\ge 0$, and $\Dp_j^{(a)}(t')\ge 0$ for all $a\ge k$ such that
$m_j^{(a)}>0$. Repeating the above construction, if necessary, eventually yields 
a new tableau $t''$ such that finally $\Dp_j^{(a)}(t'')\ge 0$ for all $j$ and $a$.
\end{proof}

\subsection{Fermionic formula} \label{subsec:fermionic}
The following statistics can be defined on the set of unrestricted rigged configurations.
For $(\nu,J)\in\RC(L,\la)$ let
\begin{equation*}
\cc(\nu,J)=\cc(\nu)+\sum_{(a,i)\in\HH}|J^{(a,i)}|,
\end{equation*}
where $|J^{(a,i)}|$ is the sum of all parts of the quasipartition $J^{(a,i)}$ and 
\begin{equation*}
\cc(\nu)=\frac{1}{2} \sum_{a,b\in I} \sum_{j,k\ge 1} (\alpha_a | \alpha_b) 
 \min(j,k) m_j^{(a)} m_k^{(b)}.
\end{equation*}

\begin{definition}
The RC polynomial is defined as
\begin{equation*}
M(L,\la)=\sum_{(\nu,J)\in\RC(L,\la)} q^{\cc(\nu,J)}.
\end{equation*}
\end{definition}
The RC polynomial is in fact $S_n$-symmetric in the weight $\la$. This is not obvious
from its definition as both~\eqref{eq:size1} and the lower bounds are not
symmetric with respect to $\la$.

Let $\SA(\la')$ be the set of all nonempty subsets of $\A(\la')$ and set
\begin{equation*}
M_i^{(a)}(S)=\max\{M_i^{(a)}(t) \mid t\in S\} \qquad \text{for $S\in\SA(\la')$.}
\end{equation*}
By inclusion-exclusion the set of all allowed riggings for a given $\nu\in\Conf(L,\la)$ is
\begin{equation*}
\bigcup_{S\in\SA(\la')} (-1)^{|S|+1} \{J\mid \text{$J^{(a,i)}$ is a 
 $(M_i^{(a)}(S),p_i^{(a)},m_i^{(a)})$-quasipartition}\}.
\end{equation*}
The $q$-binomial coefficient $\qbin{m+p}{m}$, defined as
\begin{equation*}
\qbin{m+p}{m}=\frac{(q)_{m+p}}{(q)_m(q)_p}
\end{equation*}
where $(q)_n=(1-q)(1-q^2)\cdots(1-q^n)$, is the generating function of partitions
with at most $m$ parts each not exceeding $p$. Hence the polynomial $M(L,\la)$
may be rewritten as 
\begin{multline*}
M(L,\la)=\sum_{S\in\SA(\la')} (-1)^{|S|+1} \sum_{\nu\in\Conf(L,\la)}
q^{\cc(\nu)+\sum_{(a,i)\in\HH} m_i^{(a)}M_i^{(a)}(S)}\\
\times \prod_{(a,i)\in\HH} \qbin{m_i^{(a)}+p_i^{(a)}-M_i^{(a)}(S)}{m_i^{(a)}}
\end{multline*}
called \textbf{fermionic formula}. This formula is different from the fermionic
formulas of~\cite{HKKOTY:1999,Kir:1998} which exist in the special case when
$L$ is the multiplicity array of $B=B^{1,s_k}\otimes \cdots \otimes B^{1,s_1}$
or $B=B^{r_k,1}\otimes \cdots \otimes B^{r_1,1}$.

\section{Bijection} \label{sec:main}

In this section we define the bijection $\Phi:\Path(B,\la)\to\RC(L,\la)$ from
paths to unrestricted rigged configurations algorithmically. The algorithm generalizes
the bijection of~\cite{KSS:2002} to the unrestricted case. The main result is
summarized in the following theorem.
\begin{theorem}\label{thm:bij} Let $B=B^{r_k,s_k}\otimes \cdots \otimes B^{r_1,s_1}$,
$L$ the corresponding multiplicity array and $\la=(\la_1,\ldots,\la_n)$ a sequence
of nonnegative integers.
There exists a bijection $\Phi:\Path(B,\la)\to\RC(L,\la)$ which preserves the statistics,
that is, $\Dt(b)=\cc(\Phi(b))$ for all $b\in\Path(B,\la)$.
\end{theorem}
A different proof of Theorem~\ref{thm:bij} is given in~\cite{Sch:2005b} by proving 
directly that the crystal structure on rigged configurations and paths coincide.
The results in~\cite{Sch:2005b} hold for all for all simply-laced types, not just type
$A_{n-1}^{(1)}$. Hence Theorem~\ref{thm:bij} holds whenever there is a corresponding 
bijection for the highest weight elements (for example for type $D_n^{(1)}$ for 
symmetric powers~\cite{SS:2005} and antisymmetric powers~\cite{Sch:2005b}). 
Using virtual crystals and the method of folding Dynkin diagrams, these results can 
be extended to other affine root systems. 

Here we use the crystal structure to prove that the statistics is preserved. 
It follows from Theorem~\ref{thm:commute} that the algorithmic definition for $\Phi$
of this dissertation and the definition of~\cite{Sch:2005b} agree.

An immediate corollary of Theorem~\ref{thm:bij} is the relation between the fermionic 
formula for the RC polynomial of section~\ref{sec:RC} and the unrestricted Kostka 
polynomials of section~\ref{sec:paths}.
\begin{corollary}\label{cor:X=M}
With the same assumptions as in Theorem~\ref{thm:bij}, $X(B,\la)=M(L,\la)$.
\end{corollary}

\subsection{Operations on crystals}
To define $\Phi$ we first need to introduce certain maps on paths and
rigged configurations. These maps correspond to the following operations on crystals:
\begin{enumerate}
\item If $B=B^{1,1}\otimes B'$, let $\lh(B)=B'$. This operation is called \textbf{left-hat}.
\item If $B=B^{r,s}\otimes B'$ with $s\ge 2$, let $\ls(B)=B^{r,1}\otimes B^{r,s-1}\otimes B'$.
This operation is called \textbf{left-split}.
\item If $B=B^{r,1}\otimes B'$ with $r\ge 2$, let $\lb(B)=B^{1,1}\otimes B^{r-1,1}\otimes B'$.
This operation is called \textbf{box-split}.
\end{enumerate}
In analogy we define $\lh(L)$ (resp. $\ls(L)$, $\lb(L)$) to be the multiplicity array of 
$\lh(B)$ (resp. $\ls(B)$, $\lb(B)$), if $L$ is the multiplicity array of $B$.
The corresponding maps on crystal elements are given by:
\begin{enumerate}
\item Let $b=c\otimes b'\in B^{1,1}\otimes B'$. Then $\lh(b)=b'$.
\item Let $b=c\otimes b'\in B^{r,s}\otimes B'$, where $c=c_1c_2\cdots c_s$ and $c_i$
denotes the $i$-th column of $c$. Then $\ls(b)=c_1\otimes c_2\cdots c_s\otimes b'$.
\item Let $b=\begin{array}{|c|} \hline b_1\\ \hline b_2\\ \hline \vdots\\ \hline b_r\\ \hline
\end{array}\otimes b'\in B^{r,1}\otimes B'$, where $b_1<\cdots<b_r$.
Then $\lb(b)=\begin{array}{|c|} \hline b_r\\ \hline \end{array} \otimes 
\begin{array}{|c|} \hline b_1\\ \hline \vdots \\ \hline b_{r-1}\\ \hline \end{array} 
\otimes b'$.
\end{enumerate}

In the next subsection we define the corresponding maps on rigged configurations,
and give the bijection in subsection~\ref{ss:bij}.

\subsection{Operations on rigged configurations}
Suppose $L_1^{(1)}>0$. The main algorithm on rigged configurations as defined
in~\cite{KR:1988,KSS:2002} for admissible rigged configurations can be extended
to our setting. For a tuple of nonnegative integers $\la=(\la_1,\ldots,\la_n)$,
let $\lm$ be the set of all nonnegative tuples $\mu=(\mu_1,\ldots,\mu_n)$ such that
$\la-\mu=\epsilon_r$ for some $1\le r\le n$ where $\epsilon_r$ is the canonical $r$-th unit
vector in $\Z^n$. Define $\delta:\RC(L,\la)\to \bigcup_{\mu\in\lm} \RC(\lh(L),\mu)$
by the following algorithm. Let $(\nu,J)\in\RC(L,\la)$. Set $\ell^{(0)}=1$ and repeat the
following process for $a=1,2,\ldots,n-1$ or until stopped. Find the smallest index 
$i\ge \ell^{(a-1)}$ such that $J^{(a,i)}$ is singular. If no such $i$ exists, set 
$\rk(\nu,J)=a$ and stop. Otherwise set $\ell^{(a)}=i$ and continue with $a+1$.
Set all undefined $\ell^{(a)}$ to $\infty$.

The new rigged configuration $(\tilde{\nu},\tilde{J})=\delta(\nu,J)$ is obtained by
removing a box from the selected strings and making the new strings singular
again. Explicitly
\begin{equation*}
 m_i^{(a)}(\tilde{\nu})=m_i^{(a)}(\nu)+\begin{cases}
 1 & \text{if $i=\ell^{(a)}-1$}\\
 -1 & \text{if $i=\ell^{(a)}$}\\
 0 & \text{otherwise.} \end{cases}
\end{equation*}
The partition $\tilde{J}^{(a,i)}$ is obtained from $J^{(a,i)}$ by removing
a part of size $p_i^{(a)}(\nu)$ for $i=\ell^{(a)}$,
adding a part of size $p_i^{(a)}(\tilde{\nu})$ for $i=\ell^{(a)}-1$, 
and leaving it unchanged otherwise. Then $\delta(\nu,J)\in \RC(\lh(L),\mu)$
where $\mu=\la-\epsilon_{\rk(\nu,J)}$.

\begin{prop} \label{prop:delta}
$\delta$ is well-defined.
\end{prop}

The proof is given in section~\ref{appn:delta}.

\begin{example}\label{ex:delta}
Let $L$ be the multiplicity array of $B=B^{1,1}\otimes B^{2,1}\otimes B^{2,3}$
and $\la=(2,2,2,1,1,1)$. Then
\begin{equation*}
(\nu,J)= \yngrc(2,-1,1,0) \quad \yngrc(3,0,1,-1,1,-1) \quad \yngrc(3,0)
\quad \yngrc(2,-1) \quad \yngrc(1,-1) \in \RC(L,\la).
\end{equation*}
Writing the vacancy numbers next to each part instead of the riggings we get
\begin{equation*}
\yngrc(2,-1,1,0) \quad \yngrc(3,0,1,-1,1,-1) \quad \yngrc(3,1)
\quad \yngrc(2,-1) \quad \yngrc(1,-1).
\end{equation*}
Hence $\ell^{(1)}=\ell^{(2)}=1$ and all other $\ell^{(a)}=\infty$, so that
\begin{equation*}
\delta(\nu,J)= \yngrc(2,-1) \quad \yngrc(3,0,1,-1) \quad \yngrc(3,0)
\quad \yngrc(2,-1) \quad \yngrc(1,-1).
\end{equation*}
Also $\rk(\nu,J)=3$ and $\cc(\nu,J)=2$.
\end{example}

The inverse algorithm of $\delta$ denoted by $\delta^{-1}$ is defined as follows.
Let $L_1^{(1)}=\overline{L}_1^{(1)}+1, L_i^{(k)}=\overline{L}_i^{(k)}$ for all $i,k\ne1$. 
Let $\lab$ be a weight and $\la=\lab+\epsilon_r$   for some $1\le r\le n$. 
Define $\delta^{-1}: \RC(\overline{L},\lab)\to\RC(L,\la) $ 
by the following algorithm. Let $(\nub,\Jb) \in\RC(\overline{L},\lab)$. 
Let $s^{(r)}=\infty$. For $k=r-1$ down to $1$, select the longest singular string in 
$(\nub,\Jb)^{(k)}$ of length $s^{(k)}$ (possibly of zero length) such that $s^{(k)} \le 
s^{(k+1)}$. With the convention $s^{(0)}=0$ we have $s^{(0)}\le s^{(1)}$ as well. 
$\delta^{-1}(\nub,\Jb)=(\nu,J)$ is obtained from $(\nub,\Jb)$ by adding a box 
to each of the selected strings, and resetting their labels to make them 
singular with respect to the new vacancy number for $\RC(L,\la)$, and 
leaving all other strings unchanged.

\begin{example}
Let $n=6$, $ \overline{L}_1^{(2)}=\overline{L}_3^{(2)}=1$ and $\lab=(1,1,1,2,2,1)$.
$$(\nu,J)= \yngrc(3,-1), \yngrc(3,0,3,0), \yngrc(3,-1,2,-1),\yngrc(2,-1,1,-1),
\yngrc(1,0)$$
is a rigged configuration in $\RC(L,\la)$. For $r=4$ 
 $$\delta^{-1}(\nu,J)=\yngrc(3,-1,1,-1), \yngrc(4,0,3,0), \yngrc(4,-2,2,-1),
 \yngrc(2,-1,1,-1),\yngrc(1,0).$$ 
\end{example}

\begin{prop}\label{prop:inv delta}
$\delta^{-1}$ is well defined.
\end{prop}
This proposition will also be proved in section~\ref{appn:delta}.

Let $s\ge2$. Suppose $B=B^{r,s}\otimes B'$ and $L$ the corresponding 
multiplicity array. Note that $\Conf(L,\la)\subset \Conf(\ls(L),\la)$. Under this
inclusion map, the vacancy number $p_i^{(a)}$ for $\nu$ increases by
$\delta_{a,r} \chi(i<s)$. Hence there is a well-defined injective map
$\rcls:\RC(L,\la)\rightarrow \RC(\ls(L),\la)$ given by the identity map 
$\rcls(\nu,J)=(\nu,J)$.

Suppose $r\ge2$ and $B=B^{r,1}\otimes B'$ with multiplicity array $L$.
Then there is an injection $\rclb:\RC(L,\la)\to \RC(\lb(L),\la)$ defined by adding 
singular strings of length $1$ to $(\nu,J)^{(a)}$ for $1\le a < r$. Note that the
vacancy numbers remain unchanged under $\rclb$.

\subsection{Bijection}\label{ss:bij}
The map $\Phi:\Path(B,\la)\to\RC(L,\la)$ is defined recursively by various commutative 
diagrams. Note that it is possible to go from $B=B^{r_k,s_k}\otimes 
B^{r_{k-1},s_{k-1}}\otimes \cdots \otimes B^{r_1,s_1}$ to the empty crystal 
via successive application of $\lh$, $\ls$ and
$\lb$.

\begin{definition} \label{def:bij}
Define that map $\Phi:\Path(B,\la)\rightarrow \RC(L,\la)$ such that 
the empty path maps to the empty rigged configuration and such that the following
conditions hold:
\begin{enumerate}
\item \label{bij:1} Suppose $B=B^{1,1} \otimes B'$. Then the following diagram commutes:
\begin{equation*}
\begin{CD}
\Path(B,\la) @>{\Phi}>> \RC(L,\la) \\
@V{\lh}VV @VV{\delta}V \\
\displaystyle{\bigcup_{\mu\in\lm} \Path(\lh(B),\mu)} @>>{\Phi}> 
\displaystyle{\bigcup_{\mu\in\lm}
\RC(\lh(L),\mu)}
\end{CD}
\end{equation*}
\item \label{bij:2} Suppose $B=B^{r,s} \otimes B'$ with $s\ge 2$. Then the 
following diagram commutes:
\begin{equation*}
\begin{CD}
\Path(B,\la) @>{\Phi}>> \RC(L,\la) \\
@V{\ls}VV @VV{\rcls}V \\
\Path(\ls(B),\la) @>>{\Phi}> \RC(\ls(L),\la)
\end{CD}
\end{equation*}
\item \label{bij:3} Suppose $B=B^{r,1} \otimes B'$ with $r\ge2$. Then the 
following diagram commutes:
\begin{equation*}
\begin{CD}
\Path(B,\la) @>{\Phi}>> \RC(L,\la) \\
@V{\lb}VV @VV{\rclb}V \\
\Path(\lb(B),\la) @>>{\Phi}> \RC(\lb(L),\la)
\end{CD}
\end{equation*}
\end{enumerate}
\end{definition}

\begin{prop}\label{prop:bij} The map $\Phi$ of Definition~\ref{def:bij}
is a well-defined bijection.
\end{prop}
The proof is given in section~\ref{appn:phi}.

\begin{example} \label{ex:b}
Let $B=B^{1,1}\otimes B^{2,1}\otimes B^{2,3}$ and $\la=(2,2,2,1,1,1)$. Then
\begin{equation*}
b=\young(3) \otimes \young(1,2) \otimes \young(123,456)\in\Path(B,\la)
\end{equation*}
and $\Phi(b)$ is the rigged configuration $(\nu,J)$ of Example~\ref{ex:delta}.
We have $\Dt(b)=\cc(\nu,J)=2$.
\end{example}

\begin{example}
Let $n=4$, $B=B^{2,2}\otimes B^{2,1}$ and $\la=(2,2,1,1)$. Then the multiplicity
array is $L_1^{(2)}=1,L_2^{(2)}=1$ and $L_i^{(a)}=0$ for all other $(a,i)$. There are 
7 possible unrestricted paths in $\Path(B,\la)$. For each path $b\in\Path(B,\la)$ the 
corresponding rigged configuration $(\nu,J)=\Phi(b)$ together with the tail energy 
and cocharge is summarized below.
\begin{equation*}
\begin{array}{llllll}
b\;=\; \young(11,22) \otimes \young(3,4)
&
\quad (\nu,J) \;=
& \yngrc(1,0) & \yngrc(1,-1,1,-1) & \yngrc(1,0)
&
\Dt(b) =0= \cc(\nu,J)\\[5mm]
b\;=\; \young(11,24) \otimes \young(2,3)
&
\quad (\nu,J) \;=
& \yngrc(1,-1) & \yngrc(1,0,1,0) & \yngrc(1,0)
&
\Dt(b) =1= \cc(\nu,J)\\[5mm]
b\;=\; \young(12,23) \otimes \young(1,4)
&
\quad (\nu,J) \;=
& \yngrc(1,0) & \yngrc(1,0,1,0) & \yngrc(1,-1)
&
\Dt(b) =1= \cc(\nu,J)\\[5mm]
b\;=\; \young(12,24) \otimes \young(1,3)
&
\quad (\nu,J) \;=
& \yngrc(1,0) & \yngrc(1,0,1,-1) & \yngrc(1,0)
&
\Dt(b)=1= \cc(\nu,J)\\[5mm]
b\;=\; \young(13,24) \otimes \young(1,2)
&
\quad (\nu,J) \;=
& \yngrc(1,0) & \yngrc(1,0,1,0) & \yngrc(1,0)
&
\Dt(b)=2= \cc(\nu,J)\\[5mm]
b\;=\; \young(11,23) \otimes \young(2,4)
&
\quad (\nu,J) \;=
& \yngrc(1,-1) & \yngrc(2,0) & \yngrc(1,-1)
&
\Dt(b)=0= \cc(\nu,J)\\[5mm]
b\;=\; \young(12,34) \otimes \young(1,2)
&
\quad (\nu,J) \;=
& \yngrc(1,-1) & \yngrc(2,1) & \yngrc(1,-1)
&
\Dt(b)=1= \cc(\nu,J)
\end{array}
\end{equation*}
The unrestricted Kostka polynomial in this case is $M(L,\la)=2+4q+q^2=X(B,\la)$.
\end{example}

\subsection{Crystal operators on unrestricted rigged configurations}
Let $B=B^{r_k,s_s}\otimes \cdots \otimes B^{r_1,s_1}$ and $L$ be the multiplicity 
array of $B$. Let $\Path(B)=\bigcup_{\la}\Path(B,\la)$ and $\RC(L)=\bigcup_{\la}\RC(L,\la)$. 
Note that the bijection $\Phi$ of Definition~\ref{def:bij} extends to 
a bijection from $\Path(B)$ to $\RC(L)$. Let $f_a$ and $e_a$ for 
$1\le a<n$ be the crystal operators acting on the paths in $\Path(B)$. 
In~\cite{Sch:2005b} analogous operators $\ft_a$ and $\et_a$ for $1\le a<n$ acting on rigged 
configurations in $\RC(L)$ were defined.   

\begin{definition}\cite[Definition 3.3]{Sch:2005b}
\begin{enumerate}
\item
Define $\et_a(\nu,J)$ by removing a box from a string of length $k$ in
$(\nu,J)^{(a)}$ leaving all colabels fixed and increasing the new
label by one. Here $k$ is the length of the string with the smallest
negative rigging of smallest length. If no such string exists,
$\et_a(\nu,J)$ is undefined.
\item
Define $\ft_a(\nu,J)$ by adding a box to a string of length $k$ in
$(\nu,J)^{(a)}$ leaving all colabels fixed and decreasing the new
label by one. Here $k$ is the length of the string with the smallest
nonpositive rigging of largest length. If no such string exists,
add a new string of length one and label -1.
If the result is not a valid unrestricted rigged configuration
$\ft_a(\nu,J)$ is undefined.
\end{enumerate}
\end{definition}

\begin{example}
Let $L$ be the multiplicity array of $B=B^{1,3}\otimes B^{3,2}\otimes B^{2,1}$ 
and let
\begin{equation*}
(\nu,J)= \yngrc(4,-3,1,-1) \quad \yngrc(3,0,1,1) \quad \yngrc(2,-1,1,-1) \in \RC(L).
\end{equation*} 
Then 
\begin{equation*}
\begin{split}
\ft_3(\nu,J)&= \yngrc(4,-3,1,-1) \quad \yngrc(3,1,1,1) \quad \yngrc(3,-2,1,-1)\\
\text{and} \qquad
\et_3(\nu,J)&= \yngrc(4,-3,1,-1) \quad \yngrc(3,-1,1,0) \quad \yngrc(2,1).
\end{split}
\end{equation*} 
\end{example}

Define $\widetilde{\varphi}_a(\nu,J)=\max\{k\ge 0 \mid \ft_a(\nu,J)\neq 0\}$ and
$\widetilde{\varepsilon}_a(\nu,J)=\max\{k\ge 0 \mid \et_a(\nu,J)\neq 0\}$.
The following Lemma is proven in~\cite{Sch:2005b}.
\begin{lemma} \cite[Lemma 3.6]{Sch:2005b}\label{lem:varphi}
Let $(\nu,J)\in \RC(L)$. For fixed $a\in \{1,2,\ldots,n-1\}$, let $p=p_i^{(a)}$ be the 
vacancy number for large $i$ and let $s\le 0$ be the smallest nonpositive label
in $(\nu,J)^{(a)}$; if no such label exists set $s=0$. 
Then $\widetilde{\varphi}_a(\nu,J)=p-s$.
\end{lemma}

\begin{theorem} \label{thm:commute}
Let $B=B^{r_k,s_k}\otimes \cdots \otimes B^{r_1,s_1}$ and $L$ the multiplicity array 
of $B$. Then  the following diagrams commute:
\begin{equation}\label{eq:commute}
\begin{CD}
\Path(B) @>{\Phi}>> \RC(L) \\
@V{f_a}VV @VV{\ft_a}V \\
\Path(B) @>>{\Phi}> \RC(L)
\end{CD}
\qquad \qquad
\begin{CD}
\Path(B) @>{\Phi}>> \RC(L) \\
@V{e_a}VV @VV{\et_a}V \\
\Path(B) @>>{\Phi}> \RC(L).
\end{CD}
\end{equation} 
\end{theorem}
The proof of  Theorem~\ref{thm:commute} is given in section~\ref{appn:crystal}.  
Note that Proposition~\ref{prop:bij} and Theorem~\ref{thm:commute} imply
that the operators $\ft_a,\et_a$ give a crystal structure on $\RC(L)$. 
In~\cite{Sch:2005b} it is shown directly that $\ft_a$ and $\et_a$ define a crystal 
structure on $\RC(L)$.

\subsection{Proof of Theorem~\ref{thm:bij}}
By Proposition~\ref{prop:bij} $\Phi$ is a bijection which proves the first part of
Theorem~\ref{thm:bij}. By Theorem~\ref{thm:commute} the operators $\ft_a$ and $\et_a$ give 
a crystal structure on $\RC(L)$ induced by the crystal structure on $\Path(B)$
under $\Phi$. The highest weight elements are given by the usual rigged configurations 
and highest weight paths, respectively, for which Theorem~\ref{thm:bij} is known to 
hold by~\cite{KSS:2002}. The energy function $\Dt$ is constant on classical
components. By~\cite[Theorem 3.9]{Sch:2005b} the statistics $\cc$ on rigged configurations
is also constant on classical components. Hence $\Phi$ preserves the statistic.

\section{Proof of Propositions~\ref{prop:delta} and~\ref{prop:inv delta}}\label{appn:delta}

In this section we prove Propositions~\ref{prop:delta} and~\ref{prop:inv delta},
namely that $\delta$ is a well-defined bijection.
The following remark will be useful.

\begin{remark}\label{remark:new tab}
Let $(\nu,J)$ be admissible with respect to $t\in\A(\la')$.
Suppose that $\Dp_{i-1}^{(k)}(t)+\Dp_{i+1}^{(k)}(t)\ge 1$ and $\Dp_i^{(k)}(t)=m_i^{(k)}(\nu)=0$. 
Then by \eqref{eq:dp ineq} there are five choices for the letters $i$ and $i+1$
in columns $k$ and $k+1$ of $t$:
\begin{enumerate}
\item \label{l:1} $i+1$ in column $k$;
\item \label{l:2} $i+1$ in column $k$ and $k+1$, $i$ in column $k+1$;
\item \label{l:3} $i$ in column $k+1$;
\item \label{l:4} $i$ in column $k$ and $k+1$, $i+1$ in column $k$;
\item \label{l:5} $i+1$ in column $k$, $i$ in column $k+1$.
\end{enumerate}
In cases \ref{l:1} and \ref{l:2} we have $m_i^{(k-1)}(\nu)=0$. Changing letter $i+1$
to $i$ in column $k$ to form a new tableau $t'$ has the effect
$M_i^{(k)}(t')=M_i^{(k)}(t)-1$, $M_i^{(k-1)}(t')=M_i^{(k-1)}(t)+1$ and all
other lower bounds remain unchanged.
In cases \ref{l:3} and \ref{l:4} we have $m_i^{(k+1)}(\nu)=0$. Changing letter $i$
to $i+1$ in column $k+1$ to form a new tableau $t'$ has the effect
$M_i^{(k)}(t')=M_i^{(k)}(t)-1$, $M_i^{(k+1)}(t')=M_i^{(k+1)}(t)+1$ and all
other lower bounds remain unchanged.
Finally in case \ref{l:5} either $m_i^{(k-1)}(\nu)=0$ or $m_i^{(k+1)}(\nu)=0$.
Changing $i+1$ to $i$ in column $k$ (resp. $i$ to $i+1$ in column $k+1$)
has the same effect as in case \ref{l:1} (resp. case \ref{l:3}).

This shows that under the replacement $t\mapsto t'$ we have $\Dp_i^{(k)}(t')>0$
and by Lemma~\ref{lem:convex} $(\nu,J)$ is admissible with respect to some
tableau $t''$.
\end{remark}

Let $\la $ be a weight such that  $\la_r>0$ for a given $1\le r\le n$.
Set $\lab=\la-\epsilon_r$. Recall that $c_k=\la_{k+1}+\la_{k+2}+\cdots+\la_n$ is the 
height of the $k$-th column of $t \in \A(\la')$. Let us define the map 
$\D_r:\A(\la')\to\A(\lab')$ with $\tb=\D_r(t)$ as follows. If $t_{1,r}<c_{r-1}$ then 
\begin{equation}\label{eq:t bar1}
\tb_{i,k} = \begin{cases}
 t_{i+1,k} & \text{for $1\le k\le r-1$ and $1\le i<c_k$,}\\
 t_{i,k} & \text{for $r\le k\le n$ and $1\le i \le c_k$.}
\end{cases}
\end{equation} 
If $t_{1,r}=c_{r-1}$ then there exists $1\le j\le c_r$ such that $t_{i,r} =t_{i-1,r}-1$ for 
$2\le i\le j$ and $t_{j+1,r} < t_{j,r}-1$ if $j<c_r$. In this case
\begin{equation}\label{eq:t bar2}
\tb_{i,k} = \begin{cases}
 t_{i+1,k} & \text{for $1\le k\le r-1$ and $1\le i<c_k$,}\\
 t_{i,r}-1 & \text{for $k=r$ and $1\le i\le j$,}\\
 t_{i,r} &\text{for $k=r$ and $j<i\le c_r$,}\\
 t_{i,k} & \text{for $r<k\le n$ and $1\le i \le c_k$.}
\end{cases}
\end{equation}

Note that by definition the entries of $\D_r(t)$ are strictly decreasing along columns.
Let $\cb_k=\lab_{k+1}+\cdots +\lab_n$. Then we have $\cb_k=c_k-1$ for $1\le k\le r-1$ 
and $\cb_k=c_k$ for $r\le k\le n$. Again by definition
$\tb_{j,1} \in \{1,2,\cdots , \cb_{1}\}$ for all $1\le j\le \cb_1$ and  
$\tb_{j,k} \in \{1,2,\cdots , \cb_{k-1}\}$ for all $2\le j\le \cb_k$ and $1\le k \le n$. 
Therefore,  $\D_r(t)\in\A(\lab')$.

\begin{example}
Let $t=\young(332,21,1)$ and $r=3$. Then $\D_r(t)=\young(211,1)$.
\end{example}

We will use the following lemma and remark in the proofs.
\begin{lemma} \label{lem:biggest part}
Let $B=B^{r_l,s_l}\otimes \cdots \otimes B^{r_1,s_1}$ with 
 $r_l=1=s_l$. Let  $(\overline{\nu},\overline{J})=\delta(\nu,J)$ and let
$\rk(\nu,J)=r$. For $1<k<r$ let $i=t_{1,k}$. Then one of the following 
conditions hold:
\begin{enumerate}
\item $ m_i^{(k)}(\nu)=0$ or
\item $m_i^{(k)}(\nu)=1$, in which case $\delta$ selects the part of length $i$ 
in $\nu^{(k)}$.
\end{enumerate}
\end{lemma} 
\begin{proof} 
Note that $i=t_{1,k} \ge c_k$. 
By \eqref{eq:size} we have $|\nu^{(k)}|\le c_{k}$, so
that either $m_i^{(k)}(\nu)=0$ or $i=c_k$ and
$\nu^{(k)}$ consists of just one part of size $i$. In this case
$m_i^{(k)}(\nu)=1$ and $\delta$ has to select this single part.
\end{proof} 

\begin{remark}\label{remark:rigsize} By \eqref{eq:size} we have
\begin{equation*}
\begin{split}
|\nu^{(r)}|&=|\nu^{(r-1)}|-\la_r+ \sum_{i\ge 1}s_i\chi(r_i\ge r)\\ 
|\nu^{(r+1)}| &=|\nu^{(r-1)}|-\la_r-\la_{r+1}+2\sum_{i\ge 1}s_i\chi(r_i\ge r)-
\sum_{i\ge 1}s_i\delta_{r_i,r}.
\end{split}
\end{equation*}
Note that for $a>0$ 
$$\sum_{i\ge 1}\min(a,i)L_i^{(r)}=\sum_{i\ge 1}s_i\chi(s_i \le a)\delta_{r_i,r}+\sum_
{i\ge 1}a\chi(s_i>a)\delta_{r_i,r}.$$
Then if $|\nu^{(r-1)}|=c_{r-1}-k$ for some $k\ge 0$ it follows that
\begin{equation*}
 -2|\nu^{(r)}|+|\nu^{(r+1)}|+\sum_{i\ge 1}\min(a,i)L_i^{(r)}=-2\la_{r+1}-c_{r+1}+
k-\sum_{i\ge 1}\max(s_i-a,0)\delta_{r_i,r}.
\end{equation*}
\end{remark}

\begin{proof}[Proof of Proposition~\ref{prop:delta}]
To prove that $\delta$ is well-defined it needs to be shown that
$(\nub,\Jb)=\delta(\nu,J)\in\RC(\overline{L},\lab)$. Here 
$\overline{L}$ is given by $\overline{L}_1^{(1)}=L_1^{(1)}-1$,
$\overline{L}_i^{(a)}=L_i^{(a)}$ for all other $i,a$, and 
$\lab=\la-\epsilon_r$ where $r=\rk(\nu,J)$. 

Let us first show that $\lab$ indeed has nonnegative entries. Assume the 
contrary that  $\lab_r<0$. This can happen only if $\la_r=0$ . Suppose 
$t \in \A(\la')$ is such that $M_j^{(k)}(t)\le p_j^{(k)}(\nu)$ for all $j,k$. 
By~\eqref{eq:p limit}, 
$p_i^{(r)}(\nu)=-\la_{r+1}$ for large $i$. Let $\ell$ be the size of the largest 
part in $\nu^{(r)}$, so that $m_j^{(r)}(\nu)=0$ for $j> \ell$. 
By definition of vacancy numbers, $p_i^{(r)}(\nu)\ge p_j^{(r)}(\nu)$ for $i\ge j \ge \ell$. 
Also we have $M_j^{(r)}(t)\ge -\la_{r+1}$ for all $j$. Hence, 
$-\la_{r+1}\le M_{j}^{(r)}(t)\le p_{j}^{(r)}(\nu)\le p_i^{(r)}(\nu)=-\la_{r+1}$ implies 
\begin{equation} \label{eq:peq}
M_i^{(r)}(t)=M_j^{(r)}(t)= p_{j}^{(r)}(\nu)= p_i^{(r)}(\nu) \quad \text{for all $\ell \le j \le i$.}
\end{equation}
This means that the string of length $\ell$ in $(\nu,J)^{(r)}$ is singular and 
$\Dp_j^{(r)}(t)=0$ for all $j\ge \ell$. We claim that $m_j^{(r-1)}(\nu)=0$ for $ j >\ell$. 
By \eqref{eq:dp ineq} we get 
\begin{equation*}
\begin{split}
S:=&-\chi(j\in t_{\cdot,r})+\chi(j\in t_{\cdot,r+1})
+\chi(j+1\in t_{\cdot,r})-\chi(j+1\in t_{\cdot,r+1})\\
\ge & m_j^{(r-1)}(\nu)+m_j^{(r+1)}(\nu)
\end{split}
\end{equation*}
for $j>\ell$. Clearly, $m_j^{(r-1)}(\nu)=0$ unless $1\le S \le 2$. If $S=2$ we have 
$j+1\in t_{\cdot,r}$ and $j\in t_{\cdot,r+1}$ which implies 
$M_j^{(r)}(t)=M_{j+1}^{(r)}(t)+1$, a contradiction to~\eqref{eq:peq}. 
Hence $S=2$ is not possible. Similarly, we can show that $S=1$ is not possible. 
This proves that $m_j^{(r-1)}(\nu)=0$ for $j>\ell$.
Hence $\ell^{(r-1)}\le \ell$ which contradicts the assumption that $r=\rk(\nu,J)$ 
since $(\nu,J)^{(r)}$ has a singular string of length $\ell$. Therefore $\la_r>0$.

Next we need to show that $(\nub,\Jb)$ is admissible, which means that the
parts of $\Jb$ lie between the corresponding lower bound for some $\tb\in \A(\lab')$
and the vacancy number. Let $t\in\A(\la')$ be such that $(\nu,J)$ is admissible
with respect to $t$.
By the same arguments as in the proof of Proposition 3.12 of~\cite{KSS:2002}
the only problematic case is when
\begin{equation}\label{eq:problem}
m_{\ell-1}^{(k)}(\nu)=0, \quad \Dp_{\ell-1}^{(k)}(t)=0, \quad \ell^{(k-1)}<\ell
\quad \text{and $\ell$ finite}
\end{equation}
where $\ell=\ell^{(k)}$. 

Assume that $\Dp_{\ell-2}^{(k)}(t)+\Dp_{\ell}^{(k)}(t)\ge 1$ and \eqref{eq:problem}
holds. By Remark~\ref{remark:new tab} with $i=\ell-1$, there exists a new tableau
$t'$ such that $\Dp_{\ell-1}^{(k)}(t')>0$ so that the problematic case is avoided.

Hence assume that $\Dp_{\ell-2}^{(k)}(t)+\Dp_{\ell}^{(k)}(t)=0$ and \eqref{eq:problem}
holds. Let $\ell'<\ell$ be maximal such that $m_{\ell'}^{(k)}(\nu)>0$. If no such
$\ell'$ exists, set $\ell'=0$. 

Suppose that there exists $\ell'<j<\ell$ such that $\Dp_{j-1}^{(k)}(t)>0$.
Let $i$ be the maximal such $j$. Then by Remark~\ref{remark:new tab}
we can find a new tableau $t'$ such that $\Dp_i^{(k)}(t')>0$ and 
$(\nu,J)$ is admissible with respect to $t'$.
Repeating the argument we can achieve $\Dp_{\ell-1}^{(k)}(t'')>0$ for
some new tableau $t''$, so that the problematic case does not occur.

Hence we are left to consider the case $\Dp_i^{(k)}(t)=0$ for all
$\ell'\le i\le \ell$. If $m_i^{(k-1)}(\nu)=0$ for all $\ell'<i<\ell$, then by the 
same arguments as in the proof of Proposition 3.12 of~\cite{KSS:2002} we arrive 
at a contradition since $\ell^{(k-1)}\le \ell'$, but the string of length $\ell'$ 
in $(\nu,J)^{(k)}$ is singular which implies that $\ell^{(k)}\le \ell'<\ell$. Hence 
there must exist $\ell'<i<\ell$ such that $m_i^{(k-1)}(\nu)>0$ and $\ell^{(k-1)}=i$.
By \eqref{eq:dp ineq} the same five cases as in Remark~\ref{remark:new tab}
occur as possibilities for the letters $i$ and $i+1$ in columns $k$ and $k+1$ of $t$. 
In cases~\ref{l:3}, \ref{l:4} and case~\ref{l:5} if $m_i^{(k-1)}(\nu)=2$,
we have $m_i^{(k+1)}(\nu)=0$. Replace $i$ in column $k+1$ by $i+1$ in $t$ to get
a new tableau $t'$. In all other cases $m_i^{(k-1)}(\nu)=1$; replace the letter $i+1$
in column $k$ by $i$ to obtain $t'$.
The replacement $t\mapsto t'$ yields $\Dp_i^{(k)}(t')>0$ in all cases.
The change of lower bound $M_i^{(k-1)}(t')=M_i^{(k-1)}(t)+1$ in cases \ref{l:1},
\ref{l:2} and \ref{l:5} when $m_i^{(k-1)}(\nu)\neq 2$ will not cause any
problems since $m_i^{(k-1)}(\nu)=1$ so that after the application
of $\delta$ there is no part of length $i$ in the $(k-1)$-th rigged partition. 
Then again repeated application of Remark~\ref{remark:new tab} achieves
$\Dp_{\ell-1}^{(k)}(t'')>0$ for some tableau $t''$, so that the problematic
case does not occur. 

Let $t''$ be the tableau we constructed so far. Note that in all constructions 
above,  either a letter $i+1$ in column $k$ is changed to $i$, or a letter $i$ 
in column $k+1$ is changed to $i+1$. In the latter case $i+1\le \ell\le 
|\nu^{(k)}|\le c_k$. Hence $t''$ satisfies the constraint that 
$t''_{i,k}\in\{1,2,\ldots,c_{k-1}\}$ for all $i,k$.

Now let  $\tb=\D_r(t'')$. We know $\tb \in \A(\lab')$.
We will show that  the parts of $\Jb$ lie between the corresponding lower 
bound with respect to $\tb \in \A(\lab')$ and the vacancy number. 

If $t''_{1,r}<c_{r-1}$ then by Lemma~\ref{lem:biggest part} $M_i^{(k)}(\tb)\le
M_i^{(k)}(t'')$  for all $k$ and $i$ such that $m_i^{(k)}(\nub)>0$.  Hence 
by Lemma ~\ref{lem:convex}  we have that  $(\nub,\Jb)$ is admissible with 
respect to  $\tb$. 

Let  $t''_{1,r}=c_{r-1}$. Then there exists $j$ as in the definition of $\D_r$. 
We claim that 
\begin{enumerate}
\item[(i)] $m_i^{(r-1)}(\nu)=0$ for $i>c_{r-1}-j$ and $m_{c_{r-1}-j}^{(r-1)}(\nu)\le 1$.
\item[(ii)] If $m_{c_{r-1}-j}^{(r-1)}(\nu)=1$, then $\ell^{(r-1)}=c_{r-1}-j$. 
\end{enumerate}
Note that $M_i^{(r-1)}(\tb)=M_i^{(r-1)}(t'')+1$ for $c_{r-1}-j\le i <c_{r-1}$ and 
$M_i^{(k)}(\tb)\le M_i^{(k)}(t'')$ for all other $k$ and $i$ such that 
$m_i^{(k)}(\nub)>0$. Hence if  the claim is true  using  Lemma~\ref{lem:biggest part}  
we  have $M_i^{(k)}(\tb)\le M_i^{(k)}(t'')$  for all $k$ and $i$ such that 
$m_i^{(k)}(\nub)>0$. Therefore by Lemma~\ref{lem:convex} we have that  $(\nub,\Jb)$ 
is admissible with respect to $\tb$.

It remains to prove the claim. Note that if $|\nu^{(r-1)}|<c_{r-1}-j$ then our claim 
is trivially true. Let $|\nu^{(r-1)}|=c_{r-1}-k$ for some $0\le k \le j$. If all parts of
$\nu^{(r-1)}$ are strictly less than $c_{r-1}-j$, again our claim is trivially true. 
Let the largest part in $\nu^{(r-1)}$ be $c_{r-1}-p\ge c_{r-1}-j$ for some $k\le p\le j$. 
Let $a$ be the largest part in $\nu^{(r)}$. 

First suppose $a> c_{r-1}-p$ and $a=c_r-q$ for some
$0\le q <c_r$. Then $a=c_{r-1}-(\la_r+q)$ which implies that
\begin{equation*}
M_a^{(r)}(t'')\ge-(c_r-\la_r-q)+(c_{r+1}-q)=  \la_r-\la_{r+1}.
\end{equation*}
This means $p_a^{(r)}(\nu)\le M_a^{(r)}(t'')$ since $p_b^{(r)}(\nu)\ge p_a^{(r)}(\nu)$ 
for all $b\ge a$ and $p_b^{(r)}=\la_r-\la_{r+1}$ for large $b$. If 
$p_a^{(r)}(\nu)< M_a^{(r)}(t'')$,  it contradicts that
$p_a^{(r)}(\nu)\ge M_a^{(r)}(t'')$. If $p_a^{(r)}(\nu)= M_a^{(r)}(t'')$,  it contradicts the 
fact that $r=\rk(\nu,J)$ since we get a singular part of length $a$ in $\nu^{(r)}$
which is larger than the largest part in $\nu^{(r-1)}$. 
Therefore $a> c_{r-1}-p$ is not possible. 

Hence $a\le c_{r-1}-p$. Using Remark~\ref{remark:rigsize} we get,
\begin{equation}\label{eq:max vacancy}
\begin{split} 
p_a^{(r)}(\nu)&=Q_a(\nu^{(r-1)})-2|\nu^{(r)}|+Q_a(\nu^{(r+1)}) +\sum_{i\ge 1}
\min(a,i)L_i^{(r)}\\
&\le  a+p-k-2|\nu^{(r)}|+|\nu^{(r+1)}| +\sum_{i\ge 1}\min(a,i)L_i^{(r)}\\
&=a+p-2\la_{r+1}-c_{r+1}-\sum_{i\ge 1}\max(s_i-a,0)\delta_{r_i,r}.
\end{split}
\end{equation}    
Since $ p_a^{(r)}(\nu)\ge M_a^{(r)}(t'')\ge -\la_{r+1}$ we get 
\begin{equation*}
c_{r}-(p-\sum_{i\ge 1}\max(s_i-a,0)\delta_{r_i,r})\le a\le c_r.
\end{equation*}
Hence $a=c_r-q$ for $0\le q\le p- \sum_{i\ge 1}\max(s_i-a,0)\delta_{r_i,r}$.
Then from \eqref{eq:max vacancy} with $a=c_r-q$ we get
\begin{equation}\label{eq:ub vacancy}
 p_a^{(r)}(\nu)\le p-q-\la_{r+1}- \sum_{i\ge 1}\max(s_i-a,0)\delta_{r_i,r}
\le \la_r-\la_{r+1},
\end{equation}
where we used that $0\le p-q \le \la_r$ which follows from $a=c_r-q\le c_{r-1}-p$.

If $a>c_{r-1}-j$, as in the case $a>c_{r-1}-p$ we have
$$ M_a^{(r)}(t'')\ge -(c_r-\la_r-q)+(c_{r+1}-q)= \la_r-\la_{r+1}\ge p_a^{(r)}(\nu).$$ 
Hence we get a contradiction unless $p_a^{(r)}(\nu)= M_a^{(r)}(t'')$. By 
\eqref{eq:ub vacancy} and the fact that $0\le p-q\le \la_r$ we know  
$p_a^{(r)}(\nu)=\la_r-\la_{r+1}$ happens only when $p-q=\la_r$ and 
$\sum_{i\ge 1}\max(s_i-a,0)\delta_{r_i,r}=0$. This means the largest
part in $\nu^{(r-1)}$ is of length $c_{r-1}-p=c_r-q=a$. Since we have a 
singular string of length $a$ in $\nu^{(r)}$  this contradicts the fact that $r=\rk(\nu,J)$. 

If $a\le c_{r-1}-j$ then $ M_a^{(r)}(t'')\ge -(c_r-j)+(c_{r+1}-q)= j-q-\la_{r+1}\ge 
p_a^{(r)}(\nu)$ 
because of \eqref{eq:ub vacancy} and the fact that $j\ge p$. Again we get a contradiction 
unless $p_a^{(r)}(\nu)= M_a^{(r)}(t'')$. But this happens only when 
$p_a^{(r)}(\nu)=j-q-\la_{r+1}$ which gives $p=j$ because $p_a^{(r)}(\nu)$ attains 
the right hand side of \eqref{eq:ub vacancy}.  This means the largest part in 
$\nu^{(r-1)}$ is $c_{r-1}-j$. Furthermore, for large $i$ we have
$p_i^{(r)}=\la_r-\la_{r+1}\ge j-q-\la_{r+1}+(c_{r-1}-j-a)=\la_r-\la_{r+1}$ which
shows that besides $c_{r-1}-j$ all parts in $\nu^{(r-1)}$ have to be less 
than or equal to $a$. But the part of length $a$ in $\nu^{(r)}$ is singular, so 
we have to have $c_{r-1}-j>a$ and $\ell^{(r-1)}=c_{r-1}-j$ else it will 
contradict the fact that $r=\rk(\nu,J)$. This proves our claim.   

Hence $(\nub,\Jb)$ is admissible with respect to $\tb \in \A(\lab')$ and therefore
$\delta$ is well-defined.           
\end{proof}

\begin{example}
Let $L$ be the multiplicity array of $B=(B^{1,1})^{\otimes 4}$ and $\la=
(0,1,0,1,2)$. Let 
\begin{equation*}
(\nu,J)= \yngrc(3,-1,1,2) \quad \yngrc(2,0,1,0) \quad \yngrc(2,-1,1,-1) \quad
\yngrc(2,-1)\in \RC(L,\la).
\end{equation*} 
Let $t=\young(4433,3222,211,1)$ be the corresponding lower bound 
tableau. Then 
\begin{equation*}
\delta(\nu,J)= \yngrc(3,-1) \quad \yngrc(2,0) \quad \yngrc(2,-1) \quad
\yngrc(1,-1).
\end{equation*} 
Note that in this example $\ell=\ell^{(4)}=2$ and it satisfies \eqref{eq:problem}
with $k=4$. Also $\Dp_{\ell-2}^{(4)}(t)+\Dp_{\ell}^{(4)}(t)=0$ with 
$\Dp_{i}^{(4)}(t)=0$ for all $0\le i\le \ell$. Since $m_1^{(3)}(\nu)=1$ and
$2\in t_{.,4}$ this is an example where we get the new tableau $t'$ by replacing
the $2\in t_{.,4}$ by $1$ and then the corresponding lower 
bound tableau for $\delta(\nu,J)$ is  $\D_5(t')=\young(3221,211,1)$.
\end{example}

\begin{proof}[Proof of Proposition~\ref{prop:inv delta}]
Similar to Proposition~\ref{prop:delta} we need to show that for $(\nub,\Jb)\in
\RC(\overline{L},\lab)$ we have $\delta^{-1}(\nub,\Jb)=(\nu,J)\in \RC(L,\la)$ 
where $\la=\lab+\epsilon_r$. Clearly $\la$
has nonnegative parts, so it suffices to show that $(\nu, J)$ is admissible which 
means that the parts of $J$ lie between the corresponding lower bound with 
respect to some $t \in \A (\la')$ and the vacancy number. Let 
$\tb \in \A(\lab')$ be a tableau such that $(\nub, \Jb)$ is admissible with 
respect to $\tb$.  By similar argument as in the proof of Propostion~\ref{prop:delta} 
the only problematic case occurs when
\begin{equation}\label{eq:inv problem}
m_{s+1}^{(k)}(\nub)=0, \quad \Dp_{s+1}^{(k)}(\tb)=0, \quad s<s^{(k+1)}
\quad \text{and $ s$ finite}
\end{equation}
where $s=s^{(k)}$.

Assume that $\Dp_{s}^{(k)}(\tb)+\Dp_{s+2}^{(k)}(\tb)\ge 1$ and \eqref{eq:inv problem}
holds. By Remark~\ref{remark:new tab} with $i=s+1$  there exists a new tableau
$\tb'$ such that $\Dp_{s+1}^{(k)}(\tb')>0$ so that the problematic case is avoided.

Hence assume that $\Dp_{s}^{(k)}(\tb)+\Dp_{s+2}^{(k)}(\tb)=0$ and \eqref{eq:inv problem}
holds. Let $s'>s$ be minimal such that $m_{s'}^{(k)}(\nub)>0$. If no such
$s'$ exists, set $s'=\infty$. 

Suppose that there exists $s'>j>s$ such that $\Dp_{j+1}^{(k)}(\tb)>0$.
Let $i$ be the minimal such $j$. Then by Remark~\ref{remark:new tab}
we can find a new tableau $\tb'$ such that $\Dp_i^{(k)}(\tb')>0$ and 
$(\nub,\Jb)$ is admissible with respect to $\tb'$.
Repeating the argument we can achieve $\Dp_{s+1}^{(k)}(\tb'')>0$ for
some new tableau $\tb''$, so that the problematic case does not occur.

Hence we are left to consider the case $\Dp_i^{(k)}(\tb)=0$ for all
$s'\ge i \ge s$.  First let us suppose $k<r-1$. If $m_i^{(k+1)}(\nub)=0$ for all $s'>i>s$, 
then by the similar arguments as in the proof of Proposition ~\ref{prop:delta} we arrive 
at a contradiction since $s^{(k+1)}\ge s'$, but the string of length $s'$ 
in $(\nub,\Jb)^{(k)}$ is singular which implies that $s^{(k+1)}>s^{(k)}\ge s'>s$. Hence 
there must exist $s'>i>s$ such that $m_i^{(k+1)}(\nub)>0$ and $s^{(k+1)}=i$.
By \eqref{eq:dp ineq} the same five cases as in Remark~\ref{remark:new tab}
occur as possibilities for the letters $i$ and $i+1$ in columns $k$ and $k+1$ of $\tb$. 
In cases~\ref{l:1}, \ref{l:2} and case~\ref{l:5} if $m_i^{(k+1)}(\nub)=2$,
we have $m_i^{(k-1)}(\nub)=0$. Replace $i+1$ in column $k$ by $i$ in $\tb$ to get
a new tableau $\tb'$. In all other cases $m_i^{(k-1)}(\nub)=1$; replace the letter $i$
in column $k+1$ by $i+1$ to obtain $\tb'$.
The replacement $\tb\mapsto \tb'$ yields $\Dp_i^{(k)}(\tb')>0$ in all cases.
The change of lower bound $M_i^{(k+1)}(\tb')=M_i^{(k+1)}(\tb)+1$ in cases \ref{l:3},
\ref{l:4} and \ref{l:5} when $m_i^{(k+1)}\neq 2$ will not cause any
problems since $m_i^{(k+1)}=1$ so that after the application
of $\delta^{-1}$ there is no part of length $i$ in the $(k+1)$-th rigged partition. 
Then again repeated application of Remark~\ref{remark:new tab} achieves
$\Dp_{s+1}^{(k)}(\tb'')>0$ for some tableau $\tb''$, so that the problematic
case does not occur. 

Now let us consider the case $k=r-1$. Note that $s'=\infty$ here. Else $s^{(r-1)}>s$, 
a contradiction. So, $\Dp_i^{(r-1)}(\tb)=0$ for $i>s$ which implies 
$m_i^{(r-1)}(\nub)=0$ for $i>s$, else $s^{(r-1)}>s$.  Then
by~\eqref{eq:dp ineq} with $i\ge s+1$ and $k=r-1$ we have 
\begin{equation} \label{eq:S ineq}
\begin{split}
-\chi(i\in \tb_{\cdot,r-1})+&\chi(i\in \tb_{\cdot,r})+\chi(i+1\in \tb_{\cdot,r-1})
-\chi(i+1\in \tb_{\cdot,r})\\
&\ge m_{i}^{(r-2)}(\nub)+m_{i}^{(r)}(\nub)\ge 0.
\end{split}
\end{equation}

If  $s+1\in \tb_{.,r}$  by \eqref{eq:S ineq} with $i=s+1$ there are seven choices for the 
letters $s+1$ and $s+2$ in columns $r-1$ and $r$ of $\tb$.
\begin{enumerate}
\item \label{s:1} $s+1$ in both columns $r-1$ and $r$;
\item \label{s:2} Both $s+1,s+2$ in  column $r$; 
\item \label{s:3} Both $s+1,s+2$ in columns $r-1,r$;
\item \label{s:4}  $s+1$ in  columns $r-1,r$ and $s+2$ in column $r-1$; 
\item \label{s:5} $s+1$ in column  $r$;
\item \label{s:6} $s+1$ in  column $r$ and $s+2$ in columns $r-1,r$;
\item \label{s:7} $s+1$ in  column $r$ and $s+2$ in column $r-1$.
\end{enumerate} 
First note that by \eqref{eq:S ineq} $m_{s+1}^{(r-2)}(\nub)=m_{s+1}^{(r)}(\nub)=0$ 
for cases \ref{s:1}, \ref{s:2} and \ref{s:3}. For case \ref{s:4} we have 
$m_{s+1}^{(r)}(\nub)=0$ again, else  $p_{s+1}^{(r-1)}(\tb)>p_{s}^{(r-1)}(\tb)=
M_{s}^{(r-1)}(\tb)=M_{s+1}^{(r-1)}(\tb)$,
contradiction to $\Dp_{s+1}^{(r-1)}(\tb)=0$. In cases~\ref{s:5} and~\ref{s:6}  
either $m_{s+1}^{(r)}(\nub)=0$ or $m_{s+1}^{(r-2)}(\nub)=0$ by \eqref{eq:S ineq}. 
When  $m_{s+1}^{(r-2)}(\nub)=0$ and  $m_{s+1}^{(r)}(\nub)>0$ in case~\ref{s:5} 
we have $m_i^{(r-2)}(\nub)=0$ for all $i>s+1$, else 
$p_{s+1}^{(r-1)}(\nub)\ge p_{s}^{(r-1)}(\nub)+2=
M_{s}^{(r-1)}(\tb)+2\ge M_{s+1}^{(r-1)}(\tb)-1+2>M_{s+1}^{(r-1)}(\tb)$, 
a contradiction. In case~\ref{s:7} by the same string of inequalities 
either $m_{s+1}^{(r)}(\nub)=0$ or $m_{s+1}^{(r-2)}(\nub)=0$. 

When $m_{s+1}^{(r)}(\nub)=0$ we construct a new tableau $\tb'$ from $\tb$ 
by replacing $s+1$ in column $r$ by the smallest number $i>s+1$ that  
does not appear in column $r$ of $\tb$. The effect of this change is 
$M_{s+1}^{(r)}(\tb') =M_{s+1}^{(r)}(\tb)+1$ and $M_{s+1}^{(r-1)}(\tb') =
M_{s+1}^{(r-1)}(\tb)-1$. Since $m_{s+1}^{(r)}(\nub)=0$ the first change
does not create any problem.  When $m_{s+1}^{(r)}(\nub)>0$  in cases~\ref{s:6} 
and~\ref{s:7} we change the $s+2$ in column $r-1$ to 
$s+1$. The effect of this replacement is $M_{s+1}^{(r-2)}(\tb') 
=M_{s+1}^{(r-2)}(\tb)+1$ and $M_{s+1}^{(r-1)}(\tb') =M_{s+1}^{(r-1)}(\tb)-1$. 
Since $m_{s+1}^{(r-2)}(\nub)=0$ there is no problem. When 
$m_{s+1}^{(r)}(\nub)>0$ in case~\ref{s:5} we replace the smallest 
$\tb_{j,r-1}>s+1$ by $s+1$. This has the effect that 
$M_{i}^{(r-2)}(\tb') =M_{i}^{(r-2)}(\tb)+1$ for $s+1\le i< \tb_{j,r-1}$. Since 
we have $m_i^{(r-2)}=0$ for all $i\ge s+1$ we do not have any problem. 
In all cases, replacing $\tb$ by $\tb'$ the problematic case~\eqref{eq:inv problem} 
is avoided and we have $\Dp_i^{(k)}(\tb')\ge 0$ for 
all other $i,k$ such that $m_i^{(k)}(\nub)>0$. 

Let us consider the case $s+1 \not \in \tb_{.,r}$. Note that
$M_{s}^{(r-1)}(\tb)\ge M_{s+1}^{(r-1)}(\tb)$. 
We have $m_i^{(r)}(\nub)=0=m_i^{(r-2)}(\nub)$ for all 
$i>s$, else $p_{s+1}^{(r-1)}(\nub)>p_{s}^{(r-1)}(\nub)=M_{s}^{(r-1)}(\tb)\ge 
M_{s+1}^{(r-1)}(\tb)$, contradiction to $\Dp_{s+1}^{(r-1)}(\tb)=0$. Using   
\eqref{eq:S ineq}   for $i=s+1, k=r-1$ we have four possible cases for the 
choice of the letters $s+1$ and $s+2$ in columns $r-1$ and $r$ of $\tb$. 
 \begin{enumerate}
 \item \label{i:1} $s+2$ in column $r-1$;
 \item \label{i:2} $s+2$ in columns $r-1$ and $r$;
 \item \label{i:3} $s+1$ and $s+2$ in column $r-1$;
 \item \label{i:4} no $s+1,s+2$ in both columns $r-1$ and $r$.
\end{enumerate}
We first argue that case~\ref{i:3} cannot occur.
Suppose case~\ref{i:3} holds. Then $M_{s+1}^{(r-1)}(\tb)= M_{s}^{(r-1)}(\tb)-1$ and
$M_{s+2}^{(r-1)}(\tb)= M_{s+1}^{(r-1)}(\tb)-1$. But we also have 
$\Dp_i^{(r-1)}(\tb)=0$ for $i>s$ and $m_{i}^{(r-1)}(\nub)= m_{i}^{(r-2)}(\nub)= 
m_{i}^{(r)}(\nub)$ for $i>s$. Note that $\Dp_i^{(r-1)}(\tb)=0$ implies that 
$p_{s+2}^{(r-1)}(\nub)= p_{s+1}^{(r-1)}(\nub)-1=p_s^{(r-1)}(\nub)-2$. On the other hand
$m_{i}^{(r-1)}(\nub)= m_{i}^{(r-2)}(\nub)= m_{i}^{(r)}(\nub)$ implies that 
$p_{s+2}^{(r-1)}(\nub)\ge p_s^{(r-1)}(\nub)$ and $p_{s+1}^{(r-1)}(\nub)\ge 
p_s^{(r-1)}(\nub)$ which yields a contradiction.

In cases~\ref{i:1} and~\ref{i:2} we replace the letter $s+2$ in column $r-1$ to 
$s+1$ to get a new tableau $\tb'$. The change from $\tb$ to $\tb'$ yields 
$\Dp_{s+1}^{(r-1)}(\tb')>0$ without any other change. 
In case~\ref{i:4} if there exists $\tb_{j,r-1}>s+2$ for some $j$ then we 
replace the smallest such $\tb_{j,r-1}$ by $s+1$ to construct $\tb'$. Then 
again we get $\Dp_{s+1}^{(r-1)}(\tb')>0$ without any other change since 
$m_i^{(r-2)}(\nub)=0$ for all $i>s$. On the other hand if $\tb_{1,r-1}\le s$ then  
$\cb_{r-1} \le s\le |\nub^{(r-1)}|\le \cb_{r-1}$ implies $\tb_{1,r-1}=s$. 
Note that $\tb_{1,r-2}\ge s$. Here  we will avoid the problematic 
case \eqref{eq:inv problem} by constructing a new 
tableau $t\in \A(\la')$. Let 
\begin{equation}\label{eq:sp t}
t_{i,k} = \begin{cases}
 \cb_1+1 & \text{for $k=1=i$}\\
 \cb_{k-1} +1& \text{for $2\le k\le r-2$ and $i=1$},\\
 s+1 & \text{for  $k=r-1$ and $i=1$},\\ 
 \tb_{i-1,k} & \text{for $1\le k\le r-1$ and $1< i\le \cb_k$},\\
 \tb_{i,k} & \text{for $r\le k\le n$ and $1\le i \le \cb_k$.}
\end{cases}
\end{equation} 
Note that $c_k=\cb_k+1$ for $1\le k \le r-1$
and $c_k=\cb_k$ for $r\le k \le n$. Clearly $t_{i,k}\in \{1,2,\ldots,c_{k-1}\}$ 
for all $i,k$. Column-strictness of $t$ follows since $\tb_{1,1}< \cb_1+1$ 
and $\tb_{1,k}<\cb_k+1 \le \cb_{k-1}+1$ for
$2\le k\le r-1$ and $s+1>\tb_{1,r}$ .  Hence $t \in \A(\la')$. 
Note that we have $M_{s+1}^{(r-1)}(t)=  M_{s+1}^{(r-1)}(\tb)-1<  
p_{s+1}^{(r-1)}(\nub)$, so the problematic case~\eqref{eq:inv problem} is
avoided. The fact that $(\nu,J)$ is admissible with respect to $t$ is shown
later. 

Let us now define $t\in \A(\la')$ in all other cases. 
Let $\tb''\in \A(\lab')$ be the tableau we constructed from $\tb$ so far 
except in the last case. Note that in all constructions above,  either a 
letter $i+1$ in column $k$ is changed to $i$, or 
a letter $i$ in column $k+1$ is changed to $i+1$. In the latter case 
$m_i^{(k+1)}=0$ means $i+1\le s^{(k+1)}\le |\nu^{(k+1)}|\le \cb_{k+1}\le 
\cb_k$. Hence $\tb''$ satisfies the constraint that 
$\tb''_{i,k}\in\{1,2,\ldots,\cb_{k-1}\}$ for all $i,k$. 

Let us define a new tableau $t$ from $\tb''$ in the following way:
\begin{equation}\label{eq:t}
t_{i,k} = \begin{cases}
 \cb_1+1 & \text{for $k=1=i$}\\
 \cb_{k-1} +1& \text{for $2\le k\le r-1$ and $i=1$},\\ 
 \tb''_{i-1,k} & \text{for $1\le k\le r-1$ and $1< i\le \cb_k$,}\\
 \tb''_{i,k} & \text{for $r\le k\le n$ and $1\le i \le \cb_k$.}
\end{cases}
\end{equation}    
Similarly as in \eqref{eq:sp t} we have  $t \in \A(\la')$. 

Next we show that $(\nu,J)$ is admissible with respect to $t$, that is,
the parts of $J$ lie between the corresponding lower bound with 
respect to $t \in \A(\la')$ and the vacancy number.  
Note that $s^{(k)}+1 \le |\nu^{(k)}| \le c_k \le c_{k-1}$. We distinguish the
three cases $s^{(k)}+1<c_k$, $s^{(k)}+1=c_k=c_{k-1}$ and $s^{(k)}+1=c_k<c_{k-1}$.

If $s^{(k)}+1<c_k$ for all $1\le k \le r-1$, then $M_i^{(k)}(t)=M_i^{(k)}(\tb'')$ for 
all $i,k$ such that $m_i^{(k)}(\nu)>0$.
If $s^{(k)}+1=c_{k-1}$ for some $1\le k \le r-2$, then $M_{s^{(k)}+1}^{(k)}(t)=
M_{s^{(k)}+1}^{(k)}(\tb'')$ since $c_{k-1}\ge c_{k}$. Also if
$s^{(r-1)}+1=c_{r-2}$, then $M_{s^{(r-1)}+1}^{(r-1)}(t)=
M_{s^{(r-1)}+1}^{(r-1)}(\tb'')-1$. In both cases $(\nu,J)$ is admissible 
since $M_i^{(k)}(t)\le M_i^{(k)}(\tb'')$ for all $i,k$ such that $m_i^{(k)}(\nu)>0$.

Now suppose $s^{(k)}+1=c_{k} < c_{k-1}$  for some $1\le k < r-1$. Then 
$M_{s^{(k)}+1}^{(k)}(t)=M_{s^{(k)}+1}^{(k)}(\tb'')+1$. Suppose $k$ is 
minimal satisfying this condition. Note that in this situation, $s^{(k)}= 
c_k-1=\cb_k$. This means $|\nub^{(k)}|=\cb_k$ which implies by 
definition of $|\nub^{(k)}|$ that $|\nub^{(a)}|=\cb_a$ for $a\ge k$.
Using this we get 
$$ \cb_k=s^{(k)}\le s^{(k+1)}\le\cdots \le s^{(a)}\le \cdots \le s^{(r-1)}\le |
\nub^{(r-1)}|=\cb_{r-1}\le \cb_k.$$ 
This implies  $\cb_a=s^{(a)}=s^{(a+1)}=\cb_{a+1}$ for all $k\le a \le r-2$. 
When  $s^{(a)}=s^{(a+1)}$ we have   $p_{s^{(a)}+1}^{(a)}(\nu)=
p_{s^{(a)}+1}^{(a)}(\nub)$. Hence we only need to worry  
when $\Dp_{s^{(k)}+1}^{(k)}(\tb'')= 0$.  Let $\ell$ be the largest  
part in $\nub^{(k-1)}$. If $\ell>s^{(k)}$ then by definition 
$p_{s^{(k)}+1}^{(k)}(\nub)>p_{s^{(k)}}^{(k)}(\nub)$. 
But we have $M_{s^{(k)}}^{(k)}(\tb'')\ge M_{s^{(k)}+1}^{(k)}(\tb'')$, hence  
$\Dp_{s^{(k)}+1}^{(k)}(\tb'')> 0$. Suppose $\ell\le s^{(k)}$, then  
$p_{s^{(k)}+1}^{(k)}(\nub)\ge p_{s^{(k)}}^{(k)}(\nub)$ since 
$m_i^{(k)}(\nub) =0$ for $i>s^{(k)}$. If $s^{(k)}+1\in \tb''_{.,k}$ 
then $M_{s^{(k)}}^{(k)}(\tb'')= M_{s^{(k)}+1}^{(k)}(\tb'')+1$ and we 
get $\Dp_{s^{(k)}+1}^{(k)}(\tb'')> 0$.  If  $s^{(k)}+1\not \in \tb''_{.,k}$
then there exists $\tb''_{j,k}>s^{(k)}+1$ for some $j$ and we replace 
the  smallest such $\tb''_{j,k}$  by $s^{(k)}+1$ to get a new tableau $t'$ 
from $t\in \A(\la')$ . This has the effect that  $M_{s^{(k)}+1}^{(k)}(t')=
M_{s^{(k)}+1}^{(k)}(t)-1=M_{s^{(k)}+1}^{(k)}(\tb'')$ so that 
$\Dp_{s^{(k)}+1}^{(k)}(t')\ge 0$. 

This proves  that  $(\nu,J)$ is admissible with respect to 
$t$ or $t' \in\A(\la') $. Hence $\delta^{-1}$ is well-defined.
\end{proof}

\begin{example}
Let $\overline{L}$ be the multiplicity array of $B=(B^{1,1})^{\otimes 4}$ and $\lab=
(0,1,1,1,1)$. Let 
\begin{equation*}
(\nub,\Jb)= \yngrc(3,-1,1,1) \quad \yngrc(2,-1,1,0) \quad \yngrc(1,-1,1,-1) \quad
\yngrc(1,0)\in \RC(\overline{L},\lab).
\end{equation*} 
Let $\tb=\young(4432,321,21,1)$ be the corresponding lower bound 
tableau. Then with $r=3$,
\begin{equation*}
\delta^{-1}(\nub,\Jb)= \yngrc(3,-1,1,1,1,1) \quad \yngrc(3,-1,1,1) \quad \yngrc(1,-1,1,-1)
 \quad \yngrc(1,0).
\end{equation*} 
Note that in this example we have $k=r-1=2$ and $s=s^{(2)}=2$ which satisfies 
\eqref{eq:inv problem}. Also $s+1=3\in \tb_{.,r}$, hence this is the situation when $k=r-1$ 
in \eqref{eq:inv problem} with $\Dp_i^{(r-1}(\tb)=0$ for all $i>s$ and since 
$s+1\in \tb_{.,r}$ this is case~\ref{s:7} discussed in the proof. So we get the 
corresponding 
lower bound tableau for $(\nu,J)$ by replacing $3\in \tb_{.,r}$ by $4$ and then doing the 
construction defined in \eqref{eq:t}. The lower bound tableau we get is 
$\young(5542,441,32,21,1)$.
\end{example}

\section{Proof of Proposition~\ref{prop:bij}}\label{appn:phi}
In this section a proof of Proposition~\ref{prop:bij} is given stating that
the map $\Phi$ of Definition~\ref{def:bij} is a well-defined bijection.

The proof proceeds by induction on $B$ using the fact that 
it is possible to go from $B=B^{r_k,s_k}\otimes B^{r_{k-1},s_{k-1}}
\otimes \cdots \otimes B^{r_1,s_1}$ to the 
empty crystal via successive application of $\lh$, $\ls$ and $\lb$. 
Suppose that $B$ is the empty crystal. Then both sets $\Path(B,\la)$ and 
$\RC(L,\la)$ are empty unless $\la$ is the empty partition, in which case 
$\Path(B,\la)$ consists of the empty partition and $\RC(L,\la)$ consists of 
the empty rigged configuration. In this case $\Phi$ is the unique
bijection mapping the empty partition to the empty rigged configuration.

Consider the commutative diagram~\eqref{bij:1} of Definition~\ref{def:bij}.
By induction
\begin{equation*}
 \Phi: \displaystyle{\bigcup_{\mu\in\lm} \Path(\lh(B),\mu)} \longrightarrow  
\displaystyle{\bigcup_{\mu\in\lm} \RC(\lh(L),\mu)}
\end{equation*}
is a bijection. By Propositions~\ref{prop:delta} and~\ref{prop:inv delta}
$\delta$ is a bijection, and by definition it is clear that $\lh$ is a bijection
as well. Hence $\Phi=\delta^{-1} \circ \Phi \circ \lh$ is a well-defined bijection.

Suppose that $B=B^{r,1} \otimes B'$ with $r\ge2$. By induction $\Phi$ is a bijection
for $\lb(B)=B^{1,1}\otimes B^{r-1,1} \otimes B'$. Hence to prove 
that~\eqref {bij:3} uniquely determines $\Phi$ for $B$ it suffices to show that 
$\Phi$ restricts to a bijection between the image of 
$\lb:  \Path(B,\la)\longrightarrow \Path(\lb(B),\la)$ and the
image of $\rclb : \RC(L,\la) \longrightarrow \RC(\lb(L),\la)$.
Let $b=\begin{array}{|c|} \hline b_r\\ \hline \end{array} \otimes 
\begin{array}{|c|} \hline b_1\\ \hline \vdots \\ \hline b_{r-1}\\ \hline \end{array} 
\otimes b' \in \Path(\lb(B),\la)$ with $b_{r-1}<b_r$. Let $(\nu,J)=\Phi(b)$ which
is in $\RC(\lb(L),\la)$. We will show that $(\nu,J)^{(a)}$ has a 
singular string of length one for $1\le a\le r-1$.  

By induction we know  for $(\nub,\Jb)=\Phi(\overline{b})$ where $\overline{b}
=\begin{array}{|c|} \hline b_{r-1}\\ \hline \end{array} \otimes 
\begin{array}{|c|} \hline b_1\\ \hline \vdots \\ \hline b_{r-2}\\ \hline \end{array} 
\otimes b' \in \lb(B^{r-1,1}\otimes B')$ with $b_{r-2}<b_{r-1}$, $(\nub,\Jb)^{(a)}$ 
has a singular string of length one for $1\le a\le r-2$.  Let $\overline{b}'= 
\begin{array}{|c|} \hline b_1\\ \hline \vdots \\ \hline b_{r-1}\\ \hline \end{array} 
\otimes b'$ and $(\nub',\Jb')=\Phi(\overline{b}')$. This "unsplitting" on the
rigged configuration  side removes the singular string of length one from 
$(\nub,\Jb)^{(a)}$  for $1\le a \le r-2$ yielding $(\nub',\Jb')$.

Let  $\s^{(a)}$ be the length of the selected strings by $\delta^{-1}$ 
associated with $b_{r-1}$. Note that $\s^{(a)}=0$  for $1\le a \le r-2$. 
Now let $s^{(a)}$ be the selected strings by $\delta^{-1}$ 
associated with $b_r$. Since $b_{r-1}<b_r$ we have by 
construction that $s^{(a+1)}\le \s^{(a)}$. In
particular $s^{(r-1)}\le \s^{(r-2)}=0$ and therefore,
$s^{(r-1)}=0$. This implies that $s^{(a)}=0$ for $1\le a\le r-1$. 
Hence $(\nu,J)^{(a)}$ has a singular 
string of length one for $1\le a \le r-1$. 

Conversely, let $(\nu,J)\in \rclb(\RC(L,\la))$, that is, $(\nu,J)^{(a)}$ has 
singular string of length one for $1\le a \le r-1$. Let $b=\Phi^{-1}(\nu,J)=
\begin{array}{|c|} \hline b_r\\ \hline \end{array} \otimes 
\begin{array}{|c|} \hline b_1\\ \hline \vdots \\ \hline b_{r-1}\\ \hline \end{array} 
\otimes b' \in \Path(\lb(B),\la)$. We want to show that $b_{r-1}<b_r$. 
Let $(\nub,\Jb)=\delta(\nu,J)$ and $\ell^{(a)}$ be the length of the 
selected string in $(\nu,J)^{(a)}$ by $\delta$. Then $\ell^{(a)}=1$
for $1\le a \le r-1$ and the change of vacancy numbers from $(\nu,J)$ 
to $(\nub,\Jb)$ is given by
\begin{equation}\label{eq:l change in vac}
p_i^{(a)}(\nub)= p_i^{(a)}(\nu)-\chi(\ell^{(a-1)}\le i < \ell^{(a)})
 +\chi( \ell^{(a)}\le i < \ell^{(a+1)}).
\end{equation}
This implies that $(\nub,\Jb)^{(r-1)}$ has no singular string of 
length less than $\ell^{(r)}$ since $\ell^{(r-1)}=1$. Let 
$(\nub',\Jb')=\rclb(\nub,\Jb)$. Denote by $\ellb^{(a)}$ the length of the singular 
string selected by $\delta$ in $(\nub',\Jb')^{(a)}$.  
Then by induction $\ellb^{(a)}=1$ for $1\le a \le r-2$ and by 
\eqref{eq:l change in vac} we get $\ellb^{(a)}\ge \ell^{(a+1)}$ for 
$a\ge r-1$. Therefore $\ellb^{(a)}\ge \ell^{(a+1)}$ for all $1\le a\le n$. 
Hence $b_{r-1}<b_r$. This proves that $\Phi$ in~\eqref{bij:3} is
uniquely determined.

Let us now consider the case $B=B^{r,s}\otimes B'$ where $s\ge 2$. 
Any map $\Phi$ satisfying \eqref{bij:2} is injective by definition and 
unique by induction. To prove the existence and surjectivity it 
suffices to prove that  bijection $\Phi$ maps the image of $\ls: \Path(B,\la)
\longrightarrow \Path(\ls(B),\la)$ to the image of $\rcls: \RC(L,\la) 
\longrightarrow \RC(\ls(L),\la)$. 
Let $b=c_1\otimes c \otimes b' \in \ls(\Path(B,\la))$ 
where $c=c_2c_3\cdots c_{s}$ and $c_i$
denotes the $(i-1)$-th column of $c\in B^{r,s-1}$. Let $c_1= \begin{array}{|c|} 
\hline a_1\\ \hline \vdots \\ \hline a_{r}\\ \hline \end{array}\in B^{r,1}$
and $c_2=\begin{array}{|c|} \hline b_1\\ \hline \vdots \\ \hline b_{r}\\ \hline 
\end{array}$, so that we have $a_i\le b_i$ 
for $1\le i \le r$.  Let $(\nu,J)=\Phi(b)$.
We want to show that $(\nu,J) \in \rcls(\RC(L,\la))$. To do that by definition 
of $\rcls$ it is enough to show that $(\nu,J)^{(r)}$  has no singular string of 
length less than $s$. 

Let us introduce some further notation.  Let $\overline{b}=c
\otimes b'$ and $(\nub_0,\Jb_0)=\Phi(c_3\cdots c_{s}\otimes b')$. Define  
$(\nub_i,\Jb_i)=(\rclb^{-1}\circ \delta^{-1})^{i-1}\circ \delta^{-1} (\nub_0,\Jb_0)$ 
for $1\le i \le r$ and let 
$\s_i^{(a)}$ be the length of the singular strings associated to $b_i$. Similarly 
define $(\nu_i,J_i)= (\rclb^{-1}\circ \delta^{-1})^{i-1}\circ \delta^{-1}(\nu_0,J_0)$ 
for $1\le i \le r$ and let  $s_i^{(a)}$ be the length of the singular strings 
associated to $a_i$ where $(\nu_0,J_0)=\Phi(\overline{b})$. The change of 
vacancy number from $(\nub_0,\Jb_0)$ to $(\nub_i,\Jb_i)$ is given by
\begin{equation}\label{eq:1 change in vac}
p_k^{(a)}(\nub_i)= p_k^{(a)}(\nub_0)+\sum_{m=1}^{i}\chi(\s_{m}^{(a-1)}< k \le 
\s_{m}^{(a)})-\sum_{m=1}^{i}\chi( \s_{m}^{(a)}<k \le \s_{m}^{(a+1)}),
\end{equation}
and the change of vacancy number from $(\nub_0,\Jb_0)$ to $(\nu_i,J_i)$ is given by
\begin{equation}\label{eq:2 change in vac}
\begin{split}
p_k^{(a)}(\nu_i)&= p_k^{(a)}(\nub_0)+ \sum_{m=1}^{r}\chi(
 \s_{m}^{(a-1)}< k \le \s_{m}^{(a)})-\sum_{m=1}^{r}\chi( \s_{m}^{(a)}<k \le 
 \s_{m}^{(a+1)})\\ 
 &-\delta_{a,r}\chi (k<s-1)+\sum_{m=1}^{i}\chi(s_{m}^{(a-1)}< k \le 
s_{m}^{(a)})-\sum_{m=1}^{i}\chi( s_{m}^{(a)}<k \le s_{m}^{(a+1)}).
\end{split}
\end{equation}
 
Using this we will show that $s_i^{(a)}>\s_i^{(a)}$ for all $a\ge i$ and $1\le i\le r$ 
by induction on $i$.
Note that by \eqref{eq:1 change in vac}  in $(\nu_0,J_0)^{(a)}$ the strings of
length $\s_i^{(a)}+1$ remain singular for all $i,a$. Since $a_1\le b_1$ we 
have $s_1^{(a)}>\s_1^{(a)}$ for all $a$, this starts the induction. Let  $s_i^{(a)}
>\s_i^{(a)}$ for all $a$ and for $1\le i \le k$. Then by induction hypothesis and 
\eqref{eq:2 change in vac} in $(\nu_k,J_k)^{(a)}$ the strings of length 
$\s_i^{(a)}+1$ remain singular for all $a$ and $k+1\le i \le r$,
which implies that $s_{k+1}^{(a)}\ge \s_{k+1}^{(a)}+1$. Hence $s_{k+1}^{(a)}
> \s_{k+1}^{(a)}$  which proves our claim by induction.
In particular $s_r^{(r)}>\s_r^{(r)}$. By induction $(\nub_r,\Jb_r)^{(r)}$ has no 
singular string of length strictly less than $s-1$, so $\s_r^{(r)}\ge s-1$ which
implies $s_r^{(r)}\ge s$. But note that by construction of the algorithm 
$s_r^{(a)}=0$ for $1\le a\le r-1$ and the change of vacancy numbers
from $(\nu_{r-1},J_{r-1})$ to $(\nu_{r},J_{r})=(\nu,J)$ is given by,
\begin{equation*}
p_k^{(a)}(\nu)=p_k^{(a)}(\nu_{r-1})+\chi(s_{r}^{(a-1)}< k \le s_{r}^{(a)})
                              -\chi( s_{r}^{(a)}<k \le s_{r}^{(a+1)}).
\end{equation*}                              
This implies that $(\nu,J)^{(r)}$ has no singular string less than $s_r^{(r)}$ which
means $(\nu,J)^{(r)}$ has no singular string less than $s$ and we are done.

Conversely let $(\nu,J)\in \rcls(\RC(L,\la))$ and $b=\Phi^{-1}(\nu,J)=c_1\otimes 
c \otimes b'$, same notation as before. We will show that $a_i\le b_i$ for 
$1\le i \le r$. Set $(\nu_i,J_i)=(\delta \circ \lb)^{r-i}(\nu,J)$ for $1\le i\le r$ and 
set  $(\nu_0,J_0)=\delta \circ (\delta \circ \lb)^{r-1}(\nu,J)$.  Let us denote 
the length of the string selected by $\delta$ in $(\nu_i,J_i)^{(a)}$ by $\ell_i^{(a)}$. 
Similarly set $(\nub,\Jb)=\rcls (\nu_0,J_0)$ and  $(\nub_i,\Jb_i)=
(\delta \circ \lb)^{r-i}(\nub,\Jb)$ for $1\le i\le r$ and  $(\nub_0,\Jb_0)=
\delta \circ (\delta \circ \lb)^{r-1}(\nub,\Jb)$. Denote the length of the 
string selected by $\delta$ in $(\nub_i,\Jb_i)^{(a)}$ by $\ellb_i^{(a)}$. 
We claim that $\ell_i^{(a)}> \ellb_i^{(a)}$ for all $1\le i \le r$ and all 
$i\le a \le n$. We will show this by reverse induction on $i$.

First note that the change in vacancy number from $(\nu,J)$ to $(\nu_i,J_i)$ is 
given by
\begin{equation}\label{eq:lrs1 change in vac}
p_k^{(a)}(\nu_i)= p_k^{(a)}(\nu)-\sum_{m=i+1}^{r}\chi(\ell_{m}^{(a-1)}\le k < 
\ell_{m}^{(a)})+\sum_{m=i+1}^{r}\chi( \ell_{m}^{(a)}\le k < \ell_{m}^{(a+1)}).
\end{equation}
The change in vacancy number from $(\nu,J)$ to $(\nub_i,\Jb_i)$ is given by
\begin{equation}\label{eq:lrs2 change in vac}
\begin{split}
p_k^{(a)}&(\nub_i)= p_k^{(a)}(\nu)-\sum_{m=1}^{r}\chi(\ell_{m}^{(a-1)}\le k < 
\ell_{m}^{(a)})+\sum_{m=1}^{r}\chi( \ell_{m}^{(a)}\le k < \ell_{m}^{(a+1)})\\
&+\delta_{a,r}\chi (k<s-1)-\sum_{m=i+1}^{r}\chi(\ellb_{m}^{(a-1)}\le k < 
\ellb_{m}^{(a)})+\sum_{m=i+1}^{r}\chi( \ellb_{m}^{(a)}\le k < \ellb_{m}^{(a+1)}).
\end{split}
\end{equation}
\eqref{eq:lrs1 change in vac} implies that $\ell_i^{(a)}<\ell_{i-1}^{(a)}$ and   
the string of length $\ell_j^{(a)}-1$ remains singular in $(\nu_i,J_i)^{(a)}$ 
for $i+1\le j\le r$. 
Recall that $(\nu,J)^{(r)}$ has no singular string of length less than $s$. So,
$\ell_r^{(r)}\ge s$.  By construction of the algorithm $\ellb_r^{(a)}=1$ for 
$1\le a\le r-1$. By induction $(\nub,\Jb)^{(r)}$ has no singular string of 
length less than $s-1$ and hence by \eqref{eq:lrs2 change in vac} $s-1\le 
\ellb_r^{(r)}< \ell_r^{(r)}$ since the string of length $\ell_r^{(r)}-1\ge s-1$ is 
singular. Now by using \eqref{eq:lrs1 change in vac} the 
algorithm of $\delta$ 
acting on $(\nub,\Jb)$ gives that  $\ellb_r^{(a)}< \ell_r^{(a)}$ for $a\ge r$. This
starts the induction. Suppose  $\ell_i^{(a)}> \ellb_i^{(a)}$ for all $k\le i \le r$ 
and all $i< a \le n$. 
Induction hypothesis along with \eqref{eq:lrs2 change in vac} implies that 
in $(\nub_{k-1},\Jb_{k-1})^{(a)}$ we have
$\ellb_i^{(a)}<\ellb_{i-1}^{(a)}$ for $i\ge k+1$ and   the string of length 
$\ell_j^{(a)}-1$ remains singular  for $1\le j\le k-1$. Therefore 
$\ellb_{k-1}^{(a)}=1$ for $1\le a \le k-2$ and in $(\nub_{k-1},\Jb_{k-1})^{(k-1)}$,
the smallest singular string  we know is of length $\ell_{k-1}^{(k-1)}-1$. Hence 
$\ellb_{k-1}^{(k-1)}\le \ell_{k-1}^{(k-1)}-1< \ell_{k-1}^{(k-1)}$. Then by using 
\eqref{eq:lrs2 change in vac} the algorithm of $\delta$ acting on 
$(\nub_k,\Jb_k)$ gives that  $\ellb_{k-1}^{(a)}< \ell_{k-1}^{(a)}$ for $a>k-1$. 
This proves our claim.      

But  $\ell_i^{(a)}> \ellb_i^{(a)}$ for all 
$1\le i \le r$ and all $i\le a \le n$ implies $a_i\le b_i$. So we are done.

\section{Proof of Theorem~\ref{thm:commute}}\label{appn:crystal}
 
In this section we prove that the crystal operators on paths and rigged
configurations commute with the bijection $\Phi$.

The following Lemma is a result of~\cite[Lemma 3.11]{KSS:2002} about the convexity of
the vacancy numbers.
\begin{lemma}\label{lem:convexity}({\bf Convexity})
Let $(\nu,J)\in \RC(L)$.
\begin{enumerate}
\item  For all $i,k \ge 1$ we have
 $-p_{k-1}^{(i)}(\nu)+2p_{k}^{(i)}(\nu)-p_{k+1}^{(i)}(\nu)\ge m_{k}^{(i-1)}(\nu)-
 2m_k^{(i)}(\nu)+m_{k}^{(i+1)}(\nu)$.
\item Let $m_k^{(i)}(\nu)=0$ for $a<k<b$. Then
 $p_{k}^{(i)}(\nu)\ge \min(p_{a}^{(i)}(\nu), p_{b}^{(i)}(\nu))$.
\item Let $m_k^{(i)}(\nu)=0$ for $a<k<b$. If $p_a^{(i)}(\nu)=p_{a+1}^{(i)}(\nu)$ and 
$p_{a+1}^{(i)}(\nu)\le p_{b}^{(i)}(\nu)$
then $p_{a+1}^{(i)}(\nu)=p_{k}^{(i)}(\nu)$ for all $a\le k\le b$.
\item Let $m_k^{(i)}(\nu)=0$ for $a<k<b$. If $p_b^{(i)}(\nu)=p_{b-1}^{(i)}(\nu)$ and 
$p_{b-1}^{(i)}(\nu)\le p_{a}^{(i)}(\nu)$
then $p_{b-1}^{(i)}(\nu)=p_{k}^{(i)}(\nu)$ for all $a\le k\le b$.
\end{enumerate}
\end{lemma}
\begin{proof} The proof of (1) is given in \cite[Appendix]{KS:1998} 
(see also~\eqref{eq:p ineq}), (2) follows from repeated use of (1), and the proof of (3)
and (4) follow from (1) and (2).
\end{proof}  

\begin{lemma}\label{lem:delta commute}
Let $B=B^{1,1}\otimes B'$ and let $L$ and $L'$ be the multiplicity 
arrays of $B$ and $B'$. For $1\le i<n$ the following diagrams 
commute if $\ft_i$ is always defined:
\begin{equation} \label{eq:delta commute}
\begin{CD}
\RC(L) @>{\delta}>> \RC(L') \\
@V{\ft_i}VV @VV{\ft_i}V \\
\RC(L) @>>{\delta}> \RC(L')
\end{CD}
\qquad \quad
\begin{CD}
\RC(L) @>{\delta}>> \RC(L') \\
@V{\et_i}VV @VV{\et_i}V \\
\RC(L) @>>{\delta}> \RC(L')
\end{CD}
\end{equation} 
\end{lemma}   

\begin{proof} We prove~\eqref{eq:delta commute} for $\ft_i$ here; the proof for $\et_i$ is 
similar. Let us introduce some notation. Let $(\nu,J)\in \RC(L)$ and let $\ell^{(a)}$ be 
the length of the singular string selected by $\delta$ in $(\nu,J)^{(a)}$ for $1\le a<n$. 
Let $(\nub,\Jb)=\delta(\nu,J)$ and $(\nut,\Jt)=\ft_i(\nu,J)$.  Let $\ellt^{(a)}$ be 
the length of the singular string selected by $\delta$ in $(\nut,\Jt)^{(a)}$ for 
$1\le a<n$ and $\ell$ (respectively $\ellb$) be the length of the string selected by 
$\ft_i$ in $(\nu,J)^{(i)}$ (respectively in $(\nub,\Jb)^{(i)}$). A string of length $k$ 
and label $x_k$ in $(\nu,J)^{(a)}$ is denoted by $(k,x_k)$. 

Using the definition of $\ft_i$ it is easy to see that the diagram \eqref{eq:delta commute}
commutes trivially  except when $\ell^{(i-1)}-1\le \ell\le \ell^{(i)}$. We list the 
nontrivial cases as follows:
\begin{enumerate} 
\item[(a)] $\ell^{(i-1)}<\infty, \ell^{(i)}=\infty$, $\ell+1\ge \ell^{(i-1)}$. 
\item[(b)] $\ell^{(i)}<\infty, \ell^{(i-1)}\le \ell+1\le \ell^{(i)}$. 
\item[(c)] $\ell^{(i)}<\infty$ and $\ell^{(i)}=\ell$.
\end{enumerate}
Note that since $\ft_i$ fixes all the colabels, the singular strings (except the new
string of length $\ell+1$) remain singular under the action of $\ft_i$.  Let
$(\ell,x_{\ell})$ be the string selected by $\ft_i$ in $(\nu,J)^{(i)}$.
The new string of length $\ell+1$ can be singular in $(\nut,\Jt)^{(i)}$ only if 
$p_{\ell+1}^{(i)}(\nu)=x_{\ell}+1$. Also note that by the definition of $\ft_i$ if 
$m_k^{(i)}(\nu)>0$ and $(k,x_k)$ is a string in $(\nu,J)^{(i)}$
then
\begin{equation}\label{eq:ft condition}
\begin{split}
x_{\ell} &<x_k\le p_{k}^{(i)}(\nu), \quad \text{if $k>\ell$},\\
x_{\ell} &\le x_k\le p_{k}^{(i)}(\nu),  \quad \text{if $k<\ell$}.
\end{split}
\end{equation} 
Let us now consider all the nontrivial cases.

\noindent
\textbf{Case (a):} If the new string of length $\ell+1$ in $(\nut,\Jt)^{(i)}$ is 
nonsingular, then \eqref{eq:delta commute}  commutes trivially. Let us 
consider the case when the new string of length $\ell+1$ in 
$(\nut,\Jt)^{(i)}$ is singular. We have  $p_{\ell+1}^{(i)}(\nu)=x_{\ell}+1$ 
and since $\ell^{(i-1)}<\infty, \ell^{(i)}=\infty$ we have $ p_{j}^{(i)}(\nub)
=p_{j}^{(i)}(\nu)-1$  for $j\ge \ell^{(i-1)}$.  
In particular $p_{\ell+1}^{(i)}(\nub)=p_{\ell+1}^{(i)}(\nu)-1=x_{\ell}$.
The labels in $(\nub,\Jb)^{(i)}$ are the same as in $(\nu,J)^{(i)}$. Hence $\ellb=\ell$,  
but the result is not a valid rigged configuration since  
$p_{\ell+1}^{(i)}(\nub)-2<x_{\ell}-1$.
So, $\ft_i(\nub,\Jb)$ is undefined, which contradicts the assumptions of 
Lemma~\ref{lem:delta commute}.

\noindent 
\textbf{Case (b):} If the new string of length $\ell+1$ in 
$(\nut,\Jt)^{(i)}$ is singular, we show that the following conditions hold:
\begin{enumerate}
 \item[(i)] $p_{\ell^{(i)}-1}^{(i)}(\nub)\le x_{\ell}$;
 \item[(ii)] $m_j^{(i+1)}(\nu)=0$ for $\ell<j<\ell^{(i)}$.     
\end{enumerate} 
The above conditions imply that diagram \eqref{eq:delta commute} with $\ft_i$ 
commutes for the following reason. Condition (i) implies that $\ft_i$ acts on the new 
string of length $\ell^{(i)}-1$ in $(\nub,\Jb)^{(i)}$. Condition (ii) implies that  
if $\ell^{(i+1)}<\infty$ then $\ellt^{(i+1)}=\ell^{(i+1)}$. Hence $\ellt^{(a)}=\ell^{(a)}$ 
for $ a \ne i$ and $\ellt^{(i)}=\ell+1$. This gives $\ft_i\circ \delta(\nu,J)=
\delta \circ \ft_i(\nu,J)$. 
 
If  the new string of length $\ell+1$ in $(\nut,\Jt)^{(i)}$ is nonsingular 
then the diagram \eqref{eq:delta commute} with $\ft_i$ commutes if $\ft_i$ acts on 
the same string of length $\ell$ in $(\nub,\Jb)^{(i)}$ as it did on $(\nu,J)^{(i)}$. 
In this case if $(\ell^{(i-1)}-1, p_{\ell^{(i)}-1}^{(i)}(\nub))$ is the new string created 
by $\delta$ we need to show that $x_{\ell}<p_{ \ell^{(i)}-1}^{(i)}(\nub)$.
 
Let us now consider the proof of conditions (i) and (ii) in the case when the new string 
of length $\ell+1$ in $(\nut,\Jt)^{(i)}$ is singular. Note that $p_{\ell+1}^{(i)}(\nu)=
x_{\ell}+1\le x_j$ for $j>\ell$ and $m_j^ {(i)}(\nu)>0$ by~\eqref{eq:ft condition}.
In particular if $m_{\ell+1}^ {(i)}(\nu)>0$ and $(\ell+1,x_{\ell+1})$ is 
a string in $(\nu,J)^{(i)}$  then $p_{\ell+1}^{(i)}(\nu)\le x_{\ell+1}\le 
p_{\ell+1}^{(i)}(\nu)$. This implies $p_{\ell+1}^{(i)}(\nu)=x_{\ell+1}$, 
hence $(\ell+1, x_{\ell+1})$ is a singular string which is a contradiction 
if $\ell^{(i-1)}\le \ell+1< \ell^{(i)}$. If $\ell+1=\ell^{(i)}$, it is easy to
see that~\eqref{eq:delta commute} commutes. Hence we may assume that $\ell+1< \ell^{(i)}$,
so that $m_{\ell+1}^{(i)}(\nu)=0$. 
Let $k>\ell$ be smallest so that $m_{k}^ {(i)}(\nu)>0$. 
Then by Lemma ~\ref{lem:convexity}  (2) we have
\begin{equation} \label{eq:convex1}
p_{\ell+1}^{(i)}(\nu)\ge \min(p_{\ell}^{(i)}(\nu), p_{k}^{(i)}(\nu)).
\end{equation}
If $p_{\ell}^{(i)}(\nu)>p_{k}^{(i)}(\nu)$ then by \eqref{eq:convex1} we 
get $p_{\ell+1}^{(i)}(\nu)\ge p_{k}^{(i)}(\nu)$. But 
\begin{equation}\label{eq:convex2}
\begin{split}
p_{\ell+1}^{(i)}(\nu) &\le x_k< p_{k}^{(i)}(\nu) \quad \text{if $\ell<k<\ell^{(i)}$},\\
p_{\ell+1}^{(i)}(\nu) &\le x_k= p_{k}^{(i)}(\nu) \quad \text{if $k=\ell^{(i)}$}.
\end{split}
\end{equation}
Hence $k=\ell^{(i)}$ which implies $p_{\ell+1}^{(i)}(\nu)=p_{\ell^{(i)}}^{(i)}(\nu)
<p_{\ell}^{(i)}(\nu)$ and $m_{j}^ {(i)}(\nu)=0$ for $\ell<j<\ell^{(i)}$. But now 
using Lemma ~\ref{lem:convexity} (1) we get
the following contradiction:
\begin{equation*}
0>-p_{\ell}^{(i)}(\nu)+2p_{\ell+1}^{(i)}(\nu)-p_{\ell+2}^{(i)}(\nu)\ge 
m_{\ell+1}^{(i-1)}(\nu)+ m_{\ell+1}^{(i+1)}(\nu)\ge 0.
\end{equation*}
Hence $p_{\ell}^{(i)}(\nu)\le p_{k}^{(i)}(\nu)$ and by \eqref{eq:convex1} 
we $p_{\ell+1}^{(i)}(\nu)\ge p_{\ell}^{(i)}(\nu)$. Recall that we have
\begin{align*}
p_{\ell+1}^{(i)}(\nu)&=x_{\ell}+1\le p_{\ell}^{(i)}(\nu),
 &&\text{if $\ell^{(i-1)}<\ell<\ell^{(i)}$ or 
         $(\ell, x_{\ell})$ is nonsingular},\\
p_{\ell+1}^{(i)}(\nu)&=x_{\ell}+1=p_{\ell}^{(i)}(\nu)+1,
 &&\text{if $\ell=\ell^{(i-1)}-1$ and $(\ell, x_{\ell})$ is singular}.
\end{align*}
This gives us two possible situations:
\begin{enumerate}
\item $p_{\ell+1}^{(i)}(\nu)=p_{\ell}^{(i)}(\nu)$ if $\ell^{(i-1)}<\ell<\ell^{(i)}$ or 
                            $(\ell, x_{\ell})$ is nonsingular,
\item  $p_{\ell+1}^{(i)}(\nu)= p_{\ell}^{(i)}(\nu)+1$  if $\ell=\ell^{(i-1)}-1$ and 
                          $(\ell, x_{\ell})$ is singular.
\end{enumerate}

In situation (1) using Lemma~\ref{lem:convexity} (3) we get $p_{\ell+1}^{(i)}(\nu)=
p_{j}^{(i)}(\nu)$ for $\ell+1< j \le k$. Using \eqref{eq:convex2} this implies 
$k=\ell^{(i)}$ and by convexity we get condition (ii). Also this
gives $p_{\ell^{(i)}-1}^{(i)}(\nu)=p_{\ell+1}^{(i)}(\nu)=x_{\ell}+1$ and hence 
$p_{\ell^{(i)}-1}^{(i)}(\nub)=x_{\ell}$, which proves condition (i).                        

In situation (2) we have 
\begin{align*}
p_{\ell+1}^{(i)}(\nub)&=p_{\ell+1}^{(i)}(\nu)-1
 &&\text{since $\ell^{(i-1)}=\ell+1< \ell^{(i)}$}\\
 &=p_{\ell}^{(i)}(\nu) 
 &&\text{since $p_{\ell+1}^{(i)}(\nu)=p_{\ell}^{(i)}(\nu)+1$}\\
 &=p_{\ell}^{(i)}(\nub) &&\text{since $\ell<\ell^{(i-1)}$}.
\end{align*}
Also note that $p_{k}^{(i)}(\nub)=p_{k}^{(i)}(\nu)-1$ if $k<\ell^{(i)}$ and 
$p_{k}^{(i)}(\nub)=p_{k}^{(i)}(\nu)+1$ if $k=\ell^{(i)}$. Now using \eqref{eq:convex2}
and \eqref{eq:ft condition}
we get $p_{\ell+1}^{(i)}(\nub)<p_{k}^{(i)}(\nub)$. Since $m_{\ell}^{(i)}(\nub)=
m_{\ell}^{(i)}(\nu)>0$ using Lemma~\ref{lem:convexity} (3)
we get  $p_{\ell+1}^{(i)}(\nub)=p_{\ell}^{(i)}(\nub)=p_{j}^{(i)}(\nub)$
for $\ell+2\le j \le k$.  This contradicts that   $p_{\ell+1}^{(i)}(\nub)<p_{k}^{(i)}(\nub)$,
hence situation (2) cannot occur.                                           

Now let us consider the case when  the new string of length $\ell+1$ in $(\nut,\Jt)^{(i)}$ is 
nonsingular. If $\ell+1=\ell^{(i)}$ the commutation of~\eqref{eq:delta commute} is again
fairly easy to see. Hence assume that $\ell+1<\ell^{(i)}$. Then we have 
$p_{\ell^{(i)}-1}^{(i)}(\nub)=p_{\ell^{(i)}-1}^{(i)}(\nu)-1$. If 
$m_{\ell^{(i)}-1}^{(i)}(\nu)>0$  and $(\ell^{(i)}-1,x_{\ell^{(i)}-1})$
is a string in $(\nu,J)^{(i)}$  then $x_{\ell^{(i)}-1}<p_{\ell^{(i)}-1}^{(i)}(\nu)$ since 
$\ell^{(i-1)}\le \ell+1\le \ell^{(i)}-1<\ell^{(i)}$. Hence by \eqref{eq:ft condition} we 
have $x_{\ell}<x_{\ell^{(i)}-1}<p_{\ell^{(i)}-1}^{(i)}(\nu)$ which implies $x_{\ell}<
p_{\ell^{(i)}-1}^{(i)}(\nub)$ and we are done.

If $m_{\ell^{(i)}-1}^{(i)}(\nu)=0$ let $\ell \le j<\ell^{(i)}-1$ be smallest such that 
$m_j^{(i)}(\nu)>0$.
By Lemma ~\ref{lem:convexity} (2) we get 
\begin{equation}\label{eq:min}
p_{\ell^{(i)}-1}^{(i)}(\nu)\ge \min( p_{j}^{(i)}(\nu), p_{\ell^{(i)}}^{(i)}(\nu)).
\end{equation}
Note that if $\ell<j<\ell^{(i)}$ then the string $(j,x_j)$ in $(\nu,J)^{(i)}$ is 
nonsingular and therefore  $p_{j}^{(i)}(\nu)>x_j>x_{\ell}$ by \eqref{eq:ft condition}.
Also if $(\ell^{(i)},x_{\ell^{(i)}})$ is the singular string $p_{\ell^{(i)}}^{(i)}(\nu)=
x_{\ell^{(i)}}>x_{\ell}$ by  \eqref{eq:ft condition}.
So $\min( p_{j}^{(i)}(\nu), p_{\ell^{(i)}}^{(i)}(\nu))\ge x_{\ell}+1$. 
Hence by \eqref{eq:min} $p_{\ell^{(i)}-1}^{(i)}(\nu)\ge x_{\ell}+1$.
Suppose $p_{\ell^{(i)}-1}^{(i)}(\nu)= x_{\ell}+1$. Since $p_{j}^{(i)}(\nu)>x_{\ell}+1$
we get by \eqref{eq:min}  $x_{\ell}+1=p_{\ell^{(i)}-1}^{(i)}(\nu)\ge 
p_{\ell^{(i)}}^{(i)}(\nu) \ge x_{\ell}+1$
which implies $p_{\ell^{(i)}-1}^{(i)}(\nu)= p_{\ell^{(i)}}^{(i)}(\nu)$. Since 
$p_{\ell^{(i)}-1}^{(i)}(\nu)=x_{\ell}+1\le  p_{a}^{(i)}(\nu)$ for all $j<a<\ell^{(i)}$ by
Lemma ~\ref{lem:convexity} (4) we get $p_{j}^{(i)}(\nu)=x_{\ell}+1$ which is a
contradiction. Hence $p_{\ell^{(i)}-1}^{(i)}(\nu)> x_{\ell}+1$ and we get $x_{\ell}<
p_{\ell^{(i)}-1}^{(i)}(\nub)$ as desired.

Let us consider the case $j=\ell$. If the string $(\ell, x_{\ell})$ is nonsingular
by similar argument as in the previous case we have that $p_{\ell^{(i)}-1}^{(i)}(\nu)\ge
x_{\ell}+1$. Suppose   $p_{\ell^{(i)}-1}^{(i)}(\nu)=
x_{\ell}+1$. By \eqref{eq:min} if $p_{\ell^{(i)}-1}^{(i)}(\nu)\ge p_{\ell^{(i)}}^{(i)}(\nu)\ge 
x_{\ell}+1$ we get as before that $p_{\ell^{(i)}-1}^{(i)}(\nu)= p_{\ell^{(i)}}^{(i)}(\nu)$. 
Using Lemma~\ref{lem:convexity} (4)  we can show as before that  $p_{\ell+1}^{(i)}(\nu)
=x_{\ell}+1$ which is a contradiction since the string of length $\ell+1$ is not singular
in $(\nut,\Jt)^{(i)}$. By \eqref{eq:min} if $p_{\ell^{(i)}-1}^{(i)}(\nu)\ge 
p_{\ell}^{(i)}(\nu)\ge x_{\ell}+1$ we get  $p_{\ell^{(i)}-1}^{(i)}(\nu)= p_{\ell}^{(i)}(\nu)
=x_{\ell}+1$. This implies that $p_{\ell^{(i)}-1}^{(i)}(\nu)\le p_{a}^{(i)}(\nu)$ for all 
$a>\ell$. If we use this in Lemma~\ref{lem:convexity} (1) for $k=\ell^{(i)}-1$ we get
$p_{\ell^{(i)}-1}^{(i)}(\nu)=p_{\ell^{(i)}}^{(i)}(\nu)$ and then using 
Lemma~\ref{lem:convexity} (4) we get $p_{\ell+1}^{(i)}(\nu)=x_{\ell}+1$ which is a 
contradiction as before.

Hence the only case left to be considered is when $j=\ell=\ell^{(i-1)}-1$ and the string 
$(\ell, x_{\ell})$ is singular in $(\nu,J)^{(i)}$. Here $ \min( p_{\ell}^{(i)}(\nu), 
p_{\ell^{(i)}}^{(i)}(\nu))= p_{\ell}^{(i)}(\nu)$ and therefore by \eqref{eq:min}
$p_{\ell^{(i)}-1}^{(i)}(\nu)\ge x_{\ell}$. Suppose $p_{\ell^{(i)}-1}^{(i)}(\nu)=x_{\ell}$.
Since $p_{\ell^{(i)}}^{(i)}(\nu)\ge x_{\ell}+1$ we have $p_{\ell^{(i)}-1}^{(i)}(\nu)<
p_{\ell^{(i)}}^{(i)}(\nu)$. Also, $p_{\ell^{(i)}-1}^{(i)}(\nu)\ge \min(p_{\ell}^{(i)}(\nu),
p_{\ell^{(i)}}^{(i)}(\nu))=p_{\ell}^{(i)}(\nu)=x_\ell=p_{\ell^{(i)}-1}^{(i)}(\nu)$. Using 
this in Lemma~\ref{lem:convexity} (1) for $k=\ell^{(i)}-1$ we get the following contradiction:
\begin{equation}\label{eq:contradiction1} 
 0>-p_{\ell^{(i)}-2}^{(i)}(\nu)+2p_{\ell^{(i)}-1}^{(i)}(\nu)-p_{\ell^{(i)}}^{(i)}(\nu)\ge 
m_{\ell^{(i)}-1}^{(i-1)}(\nu)+ m_{\ell^{(i)}-1}^{(i+1)}(\nu)\ge 0.
\end{equation}
Hence $p_{\ell^{(i)}-1}^{(i)}(\nu)> x_{\ell}$. Suppose $p_{\ell^{(i)}-1}^{(i)}(\nu)=
 x_{\ell}+1$. Here $p_{\ell^{(i)}-1}^{(i)}(\nu)\le p_{\ell^{(i)}}^{(i)}(\nu)$. If 
 $p_{\ell^{(i)}-1}^{(i)}(\nu)= p_{\ell^{(i)}}^{(i)}(\nu)$ as before we can show that
$p_{\ell+1}^{(i)}(\nu)=x_{\ell}+1$, which is a contradiction. Suppose
$p_{\ell^{(i)}-1}^{(i)}(\nu)< p_{\ell^{(i)}}^{(i)}(\nu)$ then $p_{\ell^{(i)}-2}^{(i)}(\nu)
\ge \min( p_{\ell}^{(i)}(\nu),p_{\ell^{(i)}}^{(i)}(\nu))=p_{\ell}^{(i)}(\nu)=x_{\ell}=
p_{\ell^{(i)}-1}^{(i)}(\nu)-1$. If $p_{\ell^{(i)}-2}^{(i)}(\nu)>p_{\ell^{(i)}-1}^{(i)}(\nu)$
we again get  the contradiction  \eqref{eq:contradiction1}. If 
$p_{\ell^{(i)}-2}^{(i)}(\nu)=p_{\ell^{(i)}-1}^{(i)}(\nu)$ using 
Lemma~\ref{lem:convexity} (1) for $k=\ell^{(i)}-1$ we get  $p_{\ell^{(i)}}^{(i)}(\nu)
=p_{\ell^{(i)}-1}^{(i)}(\nu)$ which is a contradiction to our assumption. Hence 
$p_{\ell^{(i)}-1}^{(i)}(\nu)> x_{\ell}+1$ giving $x_{\ell}<p_{\ell^{(i)}-1}^{(i)}(\nub)$.

\noindent
\textbf{Case (c):}  Note that since $\ft_i$ acts on the string   $(\ell, x_{\ell})$ in 
$(\nu,J)^{(i)}$ we have  
\begin{equation}\label{eq:convex3}
p_{\ell+1}^{(i)}(\nu)\ge x_{\ell}+1=p_{\ell}^{(i)}(\nu)+1.
\end{equation} 
If $\ft_i$ and $\delta$ select the same string of length $\ell$ in $(\nu,J)^{(i)}$ 
then $m_{\ell}^{(i)}(\nu)= 1$.
But if $\ft_i$ and $\delta$ select different strings of length $\ell$ in $(\nu,J)^{(i)}$ 
then  $m_{\ell}^{(i)}(\nu)>1$. We will consider each of these two cases separately.

If  $m_{\ell}^{(i)}(\nu)>1$  let $(\ell, x_{\ell})$ be the string selected by $\ft_i$ 
and $(\ell, p_{\ell}^{(i)}(\nu))$ be the string selected by $\delta$ in $(\nu,J)^{(i)}$. 
Note that $x_{\ell}\le p_{\ell}^{(i)}(\nu)$. To prove that the 
diagram \eqref{eq:delta commute} with $\ft_i$ commutes it is enough to show that 
$\ft_i$ acts on the same string $(\ell, x_{\ell})$ in $(\nub,\Jb)^{(i)}$ as it did 
in $(\nu,J)^{(i)}$. Hence it suffices to show that the new label in $(\nub,\Jb)^{(i)}$ 
satisfies $p_{\ell-1}^{(i)}(\nub)\ge x_{\ell}$. Note that
\begin{align*}
p_{\ell-1}^{(i)}(\nub)&=p_{\ell-1}^{(i)}(\nu)-1 &&\text{if $\ell>\ell^{(i-1)}$,}\\
p_{\ell-1}^{(i)}(\nub)&=p_{\ell-1}^{(i)}(\nu) &&\text{if $\ell=\ell^{(i-1)}$}.
\end{align*}

If $m_{\ell-1}^{(i)}(\nu)>0$  let $(\ell-1,x_{\ell-1})$ be a string in $(\nu,J)^{(i)}$. 
Then
\begin{align*}
x_{\ell}&\le x_{\ell-1}< p_{\ell-1}^{(i)}(\nu) &&\text{if $\ell>\ell^{(i-1)}$,}\\
x_{\ell}&\le x_{\ell-1}\le p_{\ell-1}^{(i)}(\nu) &&\text{if $\ell=\ell^{(i-1)}$,}
\end{align*}
which implies $p_{\ell-1}^{(i)}(\nub)\ge x_{\ell}$. 

If $m_{\ell-1}^{(i)}(\nu)=0$ let $j<\ell-1$ be largest such that $m_j^{(i)}(\nu)>0$ 
and $(j,x_j)$ be a string in $(\nu,J)^{(i)}$. Then by Lemma~\ref{lem:convexity} (2) we 
have $p_{\ell-1}^{(i)}(\nu)\ge \min(p_{j}^{(i)}(\nu),
p_{\ell}^{(i)}(\nu))$. 

If $p_{j}^{(i)}(\nu)\le p_{\ell}^{(i)}(\nu)$ then using 
\eqref{eq:ft condition} we have
\begin{align*}
p_{\ell-1}^{(i)}(\nu)&\ge p_{j}^{(i)}(\nu)>x_j\ge x_{\ell}
 &&\text{if $\ell^{(i-1)}\le j<\ell-1$,}\\
p_{\ell-1}^{(i)}(\nu)&\ge p_{j}^{(i)}(\nu)\ge x_j\ge x_{\ell}
 &&\text{if $j<\ell^{(i-1)}$.}
\end{align*}
Hence $p_{\ell-1}^{(i)}(\nub)\ge x_{\ell}$ unless 
\begin{equation}\label{eq:special case}
p_{\ell-1}^{(i)}(\nu)=p_{j}^{(i)}(\nu)=x_j= x_{\ell} \le p_{\ell}^{(i)}(\nu) \qquad 
\text{with $j<\ell^{(i-1)} \le \ell-1$.}
\end{equation}
But if this happens by Lemma~\ref{lem:convexity} we get $p_{\ell^{(i-1)}}^{(i)}(\nu)=
 p_{\ell^{(i-1)}-1}^{(i)}(\nu)\le p_{\ell^{(i-1)}+1}^{(i)}(\nu)$. Note that here 
 $m_{\ell^{(i-1)}}^{(i)}(\nu)=0$ and $m_{\ell^{(i-1)}}^{(i-1)}(\nu)>0$. Using all these 
we get the following contradiction:
\begin{equation*}
0\ge -p_{\ell^{(i-1)}-1}^{(i)}(\nu)+2p_{\ell^{(i-1)}}^{(i)}(\nu)
 -p_{\ell^{(i-1)}+1}^{(i)}(\nu)
 \ge m_{\ell^{(i-1)}}^{(i-1)}(\nu)+m_{\ell^{(i-1)}}^{(i-1)}(\nu)\ge 1.
\end{equation*}
This shows that \eqref{eq:special case} can not happen.
 
If $p_{j}^{(i)}(\nu)> p_{\ell}^{(i)}(\nu)$ then $p_{\ell-1}^{(i)}(\nu)\ge 
\min(p_{j}^{(i)}(\nu), p_{\ell}^{(i)}(\nu))=p_{\ell}^{(i)}(\nu)\ge x_{\ell}$. 
Again $p_{\ell-1}^{(i)}(\nub)\ge x_{\ell}$ unless 
\begin{equation}\label{eq:special case2}
p_{\ell-1}^{(i)}(\nu)=p_{\ell}^{(i)}(\nu)= x_{\ell} \qquad 
\text{with $\ell^{(i-1)} \le \ell-1$.}
\end{equation}
But this implies by  Lemma~\ref{lem:convexity}  that $p_{j}^{(i)}(\nu)=p_{\ell}^{(i)}(\nu)
= x_{\ell}$ which is a contradiction to our assumption. Hence \eqref{eq:special case2}
does not occur. This completes the proof when $m_{\ell}^{(i)}(\nu)>1$. 
  
If  $m_{\ell}^{(i)}(\nu)= 1$  we claim that 
\begin{enumerate}
\item[(i)]  $p_{\ell+1}^{(i)}(\nu)=x_{\ell}+1=p_{\ell}^{(i)}(\nu)+1$,
\item[(ii)]  $p_{\ell-1}^{(i)}(\nub)=x_{\ell}$,
\item[(iii)] If $\ell^{(i+1)}<\infty$ then  $\ell+1\le \ell^{(i+1)}$.  
\end{enumerate}
It is easy to see that diagram \eqref{eq:delta commute} with $\ft_i$ commutes if 
our claim is true. Condition (i) implies that the new string 
$(\ell+1,x_{\ell}-1)$ in $(\nut,\Jt)^{(i)}$ is singular and $\ellt^{(i)}=\ell+1$. 
Condition (iii) implies that $\ellt^{(i+1)}=\ell^{(i+1)}$. On the other hand condition
(ii) implies $\ellb=\ell-1$, the new string created by $\delta$ in $(\nub,\Jb)^{(i)}$. 

Let us prove our claims now. Using Lemma~\ref{lem:convexity} (1) we have
\begin{equation*}
(p_{\ell}^{(i)}(\nu)-p_{\ell-1}^{(i)}(\nu))+(p_{\ell}^{(i)}(\nu)-p_{\ell+1}^{(i)}(\nu))\ge 
m_{\ell}^{(i-1)}(\nu)-2+m_{\ell}^{(i+1)}(\nu).
\end{equation*}
which can be rewritten as
\begin{equation}\label{eq:convex4}
(p_{\ell}^{(i)}(\nu)+1-p_{\ell-1}^{(i)}(\nu))+(p_{\ell}^{(i)}(\nu)+1
 -p_{\ell+1}^{(i)}(\nu))\ge m_{\ell}^{(i-1)}(\nu)+m_{\ell}^{(i+1)}(\nu)\ge 0.
\end{equation}
Suppose  $\ell^{(i-1)}<\ell=\ell^{(i)}$. If $m_{\ell-1}^{(i)}(\nu)>0$ then the string 
$(\ell-1,x_{\ell-1})$ is nonsingular and hence by \eqref{eq:ft condition} 
$p_{\ell}^{(i)}(\nu)=x_{\ell}\le x_{\ell-1}< p_{\ell-1}^{(i)}(\nu)$.  If 
$m_{\ell-1}^{(i)}(\nu)=0$  let $j<\ell-1$ be largest such that
$m_j^{(i)}(\nu)>0$. Note that $p_{j}^{(i)}(\nu)\ge x_{\ell}=p_{\ell}^{(i)}(\nu)$, so by 
Lemma~\ref{lem:convexity} (2) we have
$p_{\ell-1}^{(i)}(\nu)\ge \min(p_{j}^{(i)}(\nu),p_{\ell}^{(i)}(\nu))=p_{\ell}^{(i)}(\nu)$.
Hence $p_{\ell-1}^{(i)}(\nu)>p_{\ell}^{(i)}(\nu)$ unless
\begin{equation}\label{eq:special case3}
p_{\ell-1}^{(i)}(\nu)=p_{\ell}^{(i)}(\nu)=p_{j}^{(i)}(\nu)=x_{\ell} 
 \text{with $j<\ell^{(i-1)}<\ell$.}
\end{equation}
But if this happens by Lemma~\ref{lem:convexity} we get $p_{\ell^{(i-1)}}^{(i)}(\nu)
=p_{\ell^{(i-1)}-1}^{(i)}(\nu)=
p_{\ell^{(i-1)}+1}^{(i)}(\nu)$ which gives us the following contradiction since 
$m_{\ell^{(i-1)}}^{(i-1)}(\nu)>0$:
\begin{equation*}
0\ge -p_{\ell^{(i-1)}-1}^{(i)}(\nu)+2p_{\ell^{(i-1)}}^{(i)}(\nu)
-p_{\ell^{(i-1)}+1}^{(i)}(\nu)\ge m_{\ell^{(i-1)}}^{(i-1)}(\nu)
+m_{\ell^{(i-1)}}^{(i-1)}(\nu)\ge 1.
\end{equation*}
Hence  \eqref{eq:special case3} cannot happen and we have $p_{\ell-1}^{(i)}(\nu)>
p_{\ell}^{(i)}(\nu)$. Now using this and \eqref{eq:convex3} in~\eqref{eq:convex4}
we get
\begin{equation*}
0\ge (p_{\ell}^{(i)}(\nu)+1-p_{\ell-1}^{(i)}(\nu))+(p_{\ell}^{(i)}(\nu)+1
 -p_{\ell+1}^{(i)}(\nu))\ge m_{\ell}^{(i-1)}(\nu)+m_{\ell}^{(i+1)}(\nu)\ge 0,
\end{equation*}
which implies $p_{\ell}^{(i)}(\nu)=p_{\ell-1}^{(i)}(\nu)-1$, $p_{\ell+1}^{(i)}(\nu)
=p_{\ell}^{(i)}(\nu)+1$, 
$m_{\ell}^{(i-1)}(\nu)=0$ and $m_{\ell}^{(i+1)}(\nu)=0$. This proves (i) and (iii). Also 
$p_{\ell}^{(i)}(\nu)=p_{\ell-1}^{(i)}(\nu)-1$ implies $p_{\ell-1}^{(i)}(\nub)
=p_{\ell-1}^{(i)}(\nu)-1=p_{\ell}^{(i)}(\nu)=x_{\ell}$.  This proves (ii).

Suppose $\ell^{(i-1)}=\ell=\ell^{(i)}$. This means $m_{\ell}^{(i-1)}(\nu)\ge 1$ and 
as before if $m_{\ell-1}^{(i)}(\nu)>0$ we have $p_{\ell}^{(i)}(\nu)=
x_{\ell}\le x_{\ell-1} \le p_{\ell-1}^{(i)}(\nu)$.  If $m_{\ell-1}^{(i)}(\nu)=0$ 
again as in the previous case we have 
$p_{\ell-1}^{(i)}(\nu)\ge \min(p_{j}^{(i)}(\nu),p_{\ell}^{(i)}(\nu))=p_{\ell}^{(i)}(\nu)$. 
Using this and~\eqref{eq:convex3} in~\eqref{eq:convex4} we get 
$p_{\ell}^{(i)}(\nu)=p_{\ell-1}^{(i)}(\nu)$, $p_{\ell+1}^{(i)}(\nu)=
p_{\ell}^{(i)}(\nu)+1$, $m_{\ell}^{(i-1)}(\nu)=1$ and $m_{\ell}^{(i+1)}(\nu)=0$. 
Note that since $\ell^{(i-1)}=\ell$, $p_{\ell-1}^{(i)}(\nub)=p_{\ell-1}^{(i)}(\nu)
=p_{\ell}^{(i)}(\nu)=x_{\ell}$. So we proved (i), (ii) and (iii).
\end{proof}

\begin{lemma}\label{lem:j commute}
Let $B=B^{r,1}\otimes B'$, $r\ge 2$ and let $L$ be the multiplicity
array of $B$. For $1\le i<n$ the following diagrams commute:
\begin{equation}\label{eq:j commute}
\begin{CD}
\RC(L) @>{\rclb}>> \RC(\lb(L)) \\
@V{\ft_i}VV @VV{\ft_i}V \\
\RC(L) @>>{\rclb}> \RC(\lb(L))
\end{CD}
\qquad \quad
\begin{CD}
\RC(L) @>{\rclb}>> \RC(\lb(L)) \\
@V{\et_i}VV @VV{\et_i}V \\
\RC(L) @>>{\rclb}> \RC(\lb(L))
\end{CD}
\end{equation} 
\end{lemma}

\begin{proof}
Note that if $i>r-1$ then the proof of  \eqref{eq:j commute} is trivial. Suppose 
$1\le i\le r-1$. The proof for $\et_i$ is very similar to the proof for $\ft_i$, 
so here we only prove~\eqref{eq:j commute} for $\ft_i$.
Let $(\nu,J)\in \RC(L)$. Let $(\ell, x_{\ell})$ be the string 
selected by $\ft_i$ in $(\nu,J)^{(i)}$. Let $(\nub,\Jb)=\rclb(\nu,J)$. By 
definition of $\rclb$ we get $(\nub,\Jb)^{(k)}$ by adding a singular string 
of length one to $(\nu,J)^{(k)}$ for $1\le k \le r-1$. Hence to show that 
the diagram~\eqref{eq:j commute} commutes it suffices to show that the label for 
the new singular string of length one in $(\nub,\Jb)^{(i)}$ satisfies
$p_1^{(i)}(\nub)\ge x_{\ell}$. 
Note that $p_1^{(i)}(\nub)= p_1^{(i)}(\nu)$ for all $1\le i\le r-1$.
					    		    	 
If $m_1^{(i)}(\nu)>0$ then $ x_1^{(i)}\ge x_{\ell}$ by 
\eqref{eq:ft condition}.  So,  $p_1^{(i)}(\nub)=
p_1^{(i)}(\nu)\ge x_{1}^{(i)}\ge x_{\ell}$.  
If $m_1^{(i)}(\nu)=0$ let $j$ be smallest such that $m_j^{(i)}(\nu)>0$ and
$(j,x_j)$ be a string in $(\nu,J)^{(i)}$. By Lemma~\ref{lem:convexity} (2) we get
$p_{1}^{(i)}(\nu)\ge \min(p_{0}^{(i)}(\nu), p_{j}^{(i)}(\nu))$. Recall that
$p_{0}^{(i)}(\nu)=0$ and $x_{\ell}\le 0$ by the definition 
of $\ft_i$. So, if $p_{j}^{(i)}(\nu)\ge 0$ then
$p_{j}^{(i)}(\nub)=p_{1}^{(i)}(\nu)\ge 0\ge x_{\ell}$.  If $p_{j}^{(i)}(\nu)< 0$ then  
$p_{1}^{(i)}(\nu)\ge p_{j}^{(i)}(\nu)$. But $p_{j}^{(i)}(\nu)\ge x_{j}
\ge x_{\ell}$. Hence $p_1^{(i)}(\nub)= p_1^{(i)}(\nu)\ge p_{j}^{(i)}(\nu) \ge 
x_{\ell}$ and we are done.
\end{proof}   

\begin{lemma}\label{lem:i commute}
Let $B=B^{r,s}\otimes B'$, $r\ge 1, s\ge 2$ and let $L$ be the multiplicity
array of $B$. For $1\le i<n$ the following diagrams commute:
\begin{equation}\label{eq:i commute}
\begin{CD}
\RC(L) @>{\rcls}>> \RC(\ls(L)) \\
@V{\ft_i}VV @VV{\ft_i}V \\
\RC(L) @>>{\rcls}> \RC(\ls(L))
\end{CD}
\qquad \quad 
\begin{CD}
\RC(L) @>{\rcls}>> \RC(\ls(L)) \\
@V{\et_i}VV @VV{\et_i}V \\
\RC(L) @>>{\rcls}> \RC(\ls(L))
\end{CD}
\end{equation} 
\end{lemma}  
\begin{proof}
Let $(\nu,J)\in \RC(L)$. By definition $\rcls$ only changes the vacancy
numbers in $(\nu,J)^{(r)}$. Hence the proof of this lemma is trivial.  
\end{proof} 

Now we will prove Theorem~\ref{thm:commute}.
\begin{proof}[Proof of Theorem~\ref{thm:commute}]
To prove this theorem we will use a diagram of the form
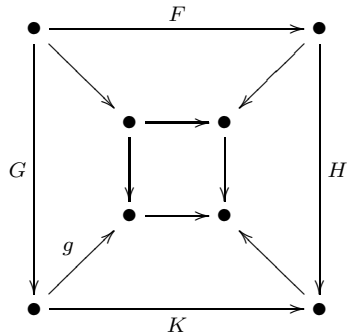
\begin{equation*}
\xymatrix{
 {\bullet} \ar[rrr]^{F} \ar[ddd]_{G} \ar[dr] & & &
        {\bullet} \ar[ddd]^{H} \ar[dl] \\
 & {\bullet} \ar[r] \ar[d] & {\bullet} \ar[d] & \\
 & {\bullet} \ar[r]   & {\bullet}  & \\
 {\bullet} \ar[rrr]_{K} \ar[ur]^{g} & & & {\bullet} \ar[ul]
}
\end{equation*}
We view this diagram as a cube with front face given by the large square. 
By~\cite[Lemma 5.3]{KSS:2002} if the squares given by all the faces of
the cube except the front commute and the map $g$ is injective then the
front face also commutes.

We will prove Theorem~\ref{thm:commute} by using induction on $B$
as we did in the proof of the bijection of Proposition~\ref{prop:bij}. 
First let $B=B^{1,1}\otimes B'$. We prove Theorem~\ref{thm:commute} for $\ft_i$ 
by using Lemma~\ref{lem:delta commute} and the following diagram when 
$f_i$ and $\ft_i$ are defined: 
\begin{equation*}
\xymatrix{
 {\Path(B)} \ar[rrr]^{\Phi} \ar[ddd]_{f_i} \ar[dr] ^{\lh}& & &
        {\RC(L)} \ar[ddd]^{\ft_i} \ar[dl]_{\delta} \\
 & {\Path(B')} \ar[r]^{\Phi} \ar[d]_{f_i} & {\RC(L')} \ar[d] ^{\ft_i}& \\
 & {\Path(B')} \ar[r]^{\Phi}  & {\RC(L')}  & \\
 {\Path(B)} \ar[rrr]_{\Phi} \ar[ur] ^{\lh}& & & {\RC(L)} \ar[ul]_{\delta}
}
\end{equation*}
Note the top and the bottom faces commute by Definition~\ref{def:bij} (1).
The right face commutes by Lemma~\ref{lem:delta commute}. The left face commutes
by definition of $f_i$ on the paths and we know $\lh$ is injective. By induction
hypothesis the back face commutes. Hence the front face must commute.

Let us now prove Theorem~\ref{thm:commute} when not all $f_i$ (resp. $\ft_i$)
in the above diagram are defined. Let $(\nu,J)\in\RC(L)$, $(\nub,\Jb)=\delta(\nu,J)$,
$b=\Phi^{-1}(\nu,J)$ and $b'=\Phi^{-1}(\nub,\Jb)$.
We need to show the following cases:
\begin{enumerate}
\item $f_i(b)$ is defined and $f_i(b')$ is undefined if and only if
      $\ft_i(\nu,J)$ is defined and $\ft_i(\nub,\Jb)$ is undefined.
      In addition $\Phi(f_i(b))=\ft_i(\nu,J)$.
\item $f_i(b)$ is undefined and $f_i(b')$ is defined if and only if
      $\ft_i(\nu,J)$ is undefined and $\ft_i(\nub,\Jb)$ is defined.
\item $f_i(b)$ and $f_i(b')$ are both undefined if and only if
      $\ft_i(\nu,J)$ and $\ft_i(\nub,\Jb)$ are both undefined.
\end{enumerate}

For Case (1) suppose that $\ft_i(\nu,J)=(\nut,\Jt)$ is defined, but $\ft_i(\nub,\Jb)$
is undefined. Then we are in the situation described in Case (a) of 
Lemma~\ref{lem:delta commute}.
That is  $\ell^{(i-1)}<\infty$, $\ell^{(i)}=\infty$, $\ell+1\ge \ell^{(i-1)}$ and 
the new string of length $\ell+1$ is singular in $(\nut,\Jt)^{(i)}$. In this situation
note that $m_{\ell+1}^{(i)}(\nub)=0$, else $p_{\ell+1}^{(i)}(\nub)\ge x_{\ell+1}>x_{\ell}$ 
by \eqref{eq:ft condition}, which is a contradiction to $p_{\ell+1}^{(i)}(\nub)=x_{\ell}$
as discussed in Case (a) of Lemma~\ref{lem:delta commute}. 
Suppose $j>\ell$ be smallest such that $m_{j}^{(i)}(\nub)>0$. Then
\begin{equation}\label{cond1}
p_j^{(i)}(\nub)\ge x_j>x_{\ell}=p_{\ell+1}^{(i)}(\nub).
\end{equation}
By Lemma~\ref{lem:convexity} (2), $p_{\ell+1}^{(i)}(\nub)\ge \min(p_{\ell}^{(i)}(\nub),
p_j^{(i)}(\nub))$. By \eqref{cond1} this implies $p_{\ell+1}^{(i)}(\nub)\ge 
p_{\ell}^{(i)}(\nub)$.
But $x_{\ell}=p_{\ell+1}^{(i)}(\nub)\ge p_{\ell}^{(i)}(\nub)\ge x_{\ell}$, hence we get
$p_{\ell+1}^{(i)}(\nub)=p_{\ell}^{(i)}(\nub)$. Again by Lemma~\ref{lem:convexity} (3) 
since $m_k^{(i)}(\nub)=0$ for $\ell< k < j$ we get $p_{\ell+1}^{(i)}(\nub)=p_{j}^{(i)}(\nub)$ 
which contradicts \eqref{cond1}.
Hence $m_j^{(i)}(\nub)=0$ for $j>\ell$. Also by Lemma~\ref{lem:convexity} (1) 
$p_{\ell+1}^{(i)}(\nub)=p_{\ell}^{(i)}(\nub)$ with $m_j^{(i)}(\nub)=0$ for $j>\ell$ 
implies that $m_j^{(i+1)}(\nub)=0$ for $j>\ell$. Since  $\nub^{(i+1)}$ and $\nut^{(i+1)}$  
have the same shape we get $m_j^{(i+1)}(\nut)=0$ for $j>\ell$. Hence $\ellt^{(a)}=\ell^{(a)}$ 
for $1\le a \le i-1$, $\ellt^{(i)}=\ell+1$ and $\ellt^{(i+1)}=\infty$.
Therefore we proved that if $\Phi^{-1}(\nub,\Jb)=b' \in B'$ then $\Phi^{-1}(\nu,J)=i 
\otimes b'$ and $\Phi^{-1}(\nut,\Jt)=i+1 \otimes b'$. But $\ft_i(\nub,\Jb)=0$ implies 
$f_i(\Phi^{-1}(\nub,\Jb))=0$ since by induction we have that $\Phi^{-1} \circ \ft_i
=f_i\circ \Phi^{-1}$ for $B'$. Hence $f_i(\Phi^{-1}(\nu,J))=\Phi^{-1}(\nut,\Jt)
=\Phi^{-1}(\ft_i(\nu,J))$, so that indeed $f_i(b)$ is defined, $f_i(b')$
and $\Phi(f_i(b))=\ft_i(\nu,J)$.

Now suppose that $f_i(b)$ is defined and $f_i(b')$ is undefined.
This implies that $b=i\otimes b'$. By induction $\ft_i(\nub,\Jb)$ is
undefined so that by Lemma~\ref{lem:varphi} we have $\overline{p}=\overline{s}$
where $\overline{p}=p_j^{(i)}(\nub)$ for large $j$ and $\overline{s}$ is the 
smallest label occurring in $(\nub,\Jb)^{(i)}$. Since $b$ is obtained from $b'$
by adding $i$ it follows that the vacancy numbers change as
$p:=p_j^{(i)}(\nu)=\overline{p}+1$ for large $j$ under $\delta^{-1}$
and the new smallest label occurring in $(\nu,J)^{(i)}$ is $s=\overline{s}$.
Hence $\widetilde{\varphi}_i(\nu,J)=p-s=1$, so that $\ft_i(\nu,J)$ is defined. 
It remains to prove that $\Phi(f_i(b))=\ft_i(\nu,J)$. Note that $f_i(b)=i+1\otimes b'$.
Let $\ell$ be the length of the largest part in $(\nub,\Jb)^{(i)}$.
Suppose that $\nub^{(i-1)}$ or $\nub^{(i+1)}$ has a part strictly bigger than $\ell$.
In this case $p_\ell^{(i)}(\nub)<\overline{p}=\overline{s}$ contradicting the fact
that $\overline{s}\le p_\ell^{(i)}(\nub)$ is the smallest label occurring in
$(\nub,\Jb)^{(i)}$. Hence both $\nub^{(i-1)}$ and $\nub^{(i+1)}$ have only parts
of length less or equal to $\ell$. Also by Lemma~\ref{lem:convex} we have
$p_\ell^{(i)}(\nub)=\overline{s}=s$ which shows that both $\delta^{-1}$ adding
$i+1$ and $\ft_i$ pick the string of length $\ell$ in $(\nub,\Jb)^{(i)}$.
Hence $\Phi(f_i(b))=\ft_i(\nu,J)$.

Let us now consider Case (2). Suppose that $\ft_i(\nu,J)$ is undefined and
$\ft_i(\nub,\Jb)$ is defined. Again by Lemma~\ref{lem:varphi} we have that
$p=s$ where $p=p_j^{(i)}(\nu)$ for large $j$ and $s$ is the smallest label
in $(\nu,J)^{(i)}$. If $\rk(\nu,J)<i+1$, then $s$ is still the smallest label in
$(\nub,\Jb)$ and by the change in vacancy numbers $\overline{p}\le p$.
Hence by Lemma~\ref{lem:varphi} $\widetilde{\varphi}_i(\nub,\Jb)=\overline{p}-s\le 0$
contradicting that $\ft_i(\nub,\Jb)$ is defined. Hence we must have
$\rk(\nu,J)\ge i+1$. In fact we want to show that $\rk(\nu,J)=i+1$.
Suppose $\rk(\nu,J)>i+1$. Then by the change in vacancy numbers by $\delta$
we have $\overline{p}=p=s$, so that $\widetilde{\varphi}_i(\nub,\Jb)=s-\overline{s}$.
So to achieve $\widetilde{\varphi}_i(\nub,\Jb)>0$ we need $\overline{s}<s$.
This can only happen if $p_{\ell^{(i)}-1}^{(i)}(\nu)=s$ and $\ell^{(i-1)}<\ell^{(i)}$.
If $m_{\ell^{(i)}-1}^{(i)}(\nu)>0$, then the string of length $\ell^{(i)}-1$ is
singular. Since $\ell^{(i-1)}<\ell^{(i)}$ this contradicts the fact that $\delta$
picks the string of length $\ell^{(i)}$ in $(\nu,J)^{(i)}$.
If $m_{\ell^{(i)}-1}^{(i)}(\nu)=0$, by convexity Lemma~\ref{lem:convexity},
we get a similar contradiction. Hence we have that $b=i+1\otimes b$.
Note that the above arguments also shows that $\widetilde{\varphi}_i(\nub,\Jb)=1$
since $\overline{s}\ge s$ and $\overline{p}=p-1$ if $\rk(\nu,J)=i+1$.
Hence $f_i(b)$ is undefined since $\varphi_i(b')=\widetilde{\varphi}_i(\nub,\Jb)
=1$.

Consider Case (2) where $f_i(b)$ is undefined and $f_i(b')$ is defined.
This implies that $b=i+1\otimes b'$. By induction $\widetilde{\varphi}_i(\nub,\Jb)
=\varphi_i(b')=1$ so that by Lemma~\ref{lem:varphi} we have $\overline{p}=\overline{s}+1$.
Hence $\widetilde{\varphi}_i(\nu,J)=p-s=\overline{p}-1-s=\overline{s}-s$ by the change 
of vacancy numbers. Therefore $\widetilde{\varphi}_i(\nu,J)=0$ if $\overline{s}=s$. It remains
to show that $p_{\ell+1}^{(i)}(\nu)\ge \overline{s}$ where $\ell:=s^{(i)}$ is the length 
of the string in $(\nub,\Jb)^{(i)}$ selected by $\delta^{-1}$.
Hence the only problem occurs if $p_{\ell+1}^{(i)}(\nub)=\overline{s}$ and
$s^{(i-1)}<\ell$. If $m_{\ell+1}^{(i)}(\nub)>0$, this means that
there is a singular string of length $\ell+1>s^{(i)}$ in $(\nub,\Jb)^{(i)}$
contradicting the maximality of $s^{(i)}$. If $m_{\ell+1}^{(i)}(\nub)=0$
one can again use convexity to arrive at similar contradiction.

By exclusion Case (3) follows from all the previous cases where at least 
one $f_i$ or $\ft_i$ is defined.

Now let $B=B^{r,1}\otimes B'$ where $r\ge 2$. Consider the following diagram:  
\begin{equation*}
\xymatrix{
 {\Path(B)} \ar[rrr]^{\Phi} \ar[ddd]_{f_i} \ar[dr] ^{\lb}& & &
        {\RC(L)} \ar[ddd]^{\ft_i} \ar[dl]_{\rclb} \\
 & {\Path(\lb(B))} \ar[r]^{\Phi} \ar[d]_{f_i} & {\RC(\lb(L))} \ar[d] ^{\ft_i}& \\
 & {\Path(\lb(B))} \ar[r]^{\Phi}  & {\RC(\lb(L))}  & \\
 {\Path(B)} \ar[rrr]_{\Phi} \ar[ur] ^{\lb}& & & {\RC(L)} \ar[ul]_{\rclb}
}
\end{equation*}
Again  the top and the bottom faces commute because of Definition~\ref{def:bij} (3).
The right face commutes by Lemma~\ref{lem:j commute}. The left face commutes
by definition of $f_i$ on the paths and we know $\lb$ is injective. By induction
hypothesis the back face commutes too. Hence the front face commutes.

Finally let $B=B^{r,s}\otimes B'$ where $s\ge 2$. Consider the following diagram:
\begin{equation*}
\xymatrix{
 {\Path(B)} \ar[rrr]^{\Phi} \ar[ddd]_{f_i} \ar[dr] ^{\ls}& & &
        {\RC(L)} \ar[ddd]^{\ft_i} \ar[dl]_{\rcls} \\
 & {\Path(\ls(B))} \ar[r]^{\Phi} \ar[d]_{f_i} & {\RC(\ls(L))} \ar[d] ^{\ft_i}& \\
 & {\Path(\ls(B))} \ar[r]^{\Phi}  & {\RC(\ls(L))}  & \\
 {\Path(B)} \ar[rrr]_{\Phi} \ar[ur] ^{\ls}& & & {\RC(L)} \ar[ul]_{\rcls}
}
\end{equation*}
As in the previous cases by Definition~\ref{def:bij} (2), Lemma~\ref{lem:i commute}
and induction hypothesis all the faces commute except the front. Since the map $\ls$
is injective the  front face of the above diagram commutes. This completes the proof
of Theorem~\ref{thm:commute}.  
\end{proof}

    %
    %

    \newchap{Fermionic formulas for the characters of $N=1$ and $N=2$ superconformal algebras}
    \label{chap:bailey}
            \section{Introduction}
\noindent
Bailey's lemma is a powerful method to prove $q$-series identities of the 
Rogers--Ramanujan-type~\cite{Bailey:1949}. One of the key features of Bailey's lemma
is its iterative structure which was first observed by Andrews~\cite{A:1984}
(see also~\cite{Paule:1985}). This iterative structure called the Bailey chain 
makes it possible to start with one seed identity and derive an infinite family 
of identities from it. The Bailey chain has been generalized to the Bailey 
lattice~\cite{AAB:1987} which yields a whole tree of identities from a single seed.

The relevance of the Andrews--Bailey construction to physics was first revealed
in the papers by Foda and Quano~\cite{FQ:1995,FQ:1996} in which they derived
identities for the Virasoro characters using Bailey's lemma.
By the application of Bailey's lemma to polynomial versions of the character
identity of one conformal field theory, one obtains character identities of
another conformal field theory. This relation between the two conformal
field theories is called the Bailey flow. In~\cite{BMS:1995} it was demonstrated 
that there is a Bailey flow from the minimal models $M(p-1,p)$ to
$N=1$ and $N=2$ superconformal models. More precisely, it was shown that there
is a Bailey flow from $M(p-1,p)$ to $M(p,p+1)$, and from $M(p-1,p)$ to the $N=1$ 
superconformal model $SM(p,p+2)$ and the unitary $N=2$ superconformal model with 
central charge $c=3(1-\frac{2}{p})$. In the conclusions of~\cite{BMS:1995} it was
conjecture that this construction can also be carried out for the nonunitary
minimal models $M(p,p')$ where $p$ and $p'$ are relatively prime.
In this chapter of the dissertation we consider the nonunitary case. We show that starting 
with character identities for the nonunitary minimal model $M(p,p')$ 
of~\cite{BMS:1997,W:2002}, characters of the $N=1$ superconformal models $SM(p',2p+p')$, 
$SM(p',3p'-2p)$ and of the $N=2$ superconformal model with central element 
$c=3(1-\frac{2p}{p'})$ can be obtained via the Bailey flow. 
We also give a new Ramond-sector character formula for a representation of the 
$N=2$ superconformal model with central element $c=3(1-\frac{2p}{p'})$ and calculate
the corresponding fermionic formula.

The chapter is organized as follows. In section~\ref{sec:bailey} we provide the 
necessary background about Bailey pairs. In section~\ref{sec:baileypairs} we derive 
new  Bailey pairs using the Bose-Fermi identity for the minimal model $M(p,p')$.
In section~\ref{sec:fermionic} we state the fermionic formulas of the $M(p,p')$ 
models following \cite{BMS:1995,BMS:1997, BMSW:1997}.  
In section~\ref{sec:N1} necessary background for $N=1$ superconformal algebras is stated 
and the characters of the $N=1$ supersymmetric models 
$SM(2p+p',p')$ and $SM(3p'-2p,p')$ are derived using the Bailey flow. 
Explicit fermionic expressions for these characters are given. 
In section~\ref{sec:N2} the background regarding $N=2$ superconformal models is stated
and a new character for the Ramond-sector is derived. Then it is demonstrated how to 
obtain the characters of the $N=2$ superconformal model with central element 
$c=3(1-\frac{2p}{p'})$ via the Bailey flow along with the explicit fermionic expressions
for these characters. In section~\ref{sec:conclusion} we conclude with some remarks.

\section{Bailey's lemma} \label{sec:bailey}

Following~\cite{BMS:1995}, we are going to use an extended definition of Bailey pair
called the bilateral Bailey pair.
A pair $(\alpha_n,\beta_n)$ of sequences $\{\alpha_n\}_{n\in\Z}$ and
$\{\beta_n\}_{n\in\Z}$ is said to be a \textbf{bilateral Bailey pair} with respect 
to $a$ if 
\begin{equation}\label{eq:defbail}
\beta_n=\sum_{j=-\infty}^n \frac{\alpha_j}{(q)_{n-j}(aq)_{n+j}}.
\end{equation}

\begin{theorem}[\textbf{Bilateral Bailey's lemma}~\cite{A:1984,Bailey:1949,BMS:1995}]
If $(\alpha_n,\beta_n)$ is a bilateral Bailey pair then 
\begin{equation}\label{eq:bilaterlemma}
\begin{split}
\sum_{n=-\infty}^{\infty}(\rho_1)_n(\rho_2)_n & (aq/\rho_1\rho_2)^n\beta_n\\
= & \frac{(aq/\rho_1)_{\infty}(aq/\rho_2)_{\infty}}{(aq)_{\infty}(aq/\rho_1 
\rho_2)_{\infty}}\sum_{n=-\infty}^{\infty} \frac{(\rho_1)_n(\rho_2)_n
(aq/\rho_1\rho_2)^n\alpha_n}{(aq/\rho_1)_n(aq/\rho_2)_n}.
\end{split}
\end{equation}
\end{theorem}

This lemma has been used with various Bailey pairs and different specializations 
of the parameters $\rho_1$ and $\rho_2$ to prove many $q$-series identities 
(see for example~\cite{AAB:1987,BMS:1995,FQ:1996,Sl:1951}).
In this chapter the bilateral Bailey's lemma is used to derive character identities
for $N=1,2$ superconformal algebras from nonunitary minimal models 
$M(p,p')$.

It was observed by looking at the  famous list of 130 identities of Slater~\cite{Sl:1951}  that the  
specialization \eqref{eq:infinity_spec} leads to characters of minimal model $M(p,p')$
and the second specialization \eqref{eq:half_finite_spec} leads to the characters of 
$N=1$ superconformal model. Hence  by putting suitable Bailey pairs and then 
using some appropriate specialization of the parameters one can derive Bose-Fermi 
type character identities for CFTs. Therefore, the main question
is: how do we find suitable Bailey pairs? The sources for the list of Bailey pairs used 
by Rogers, Bailey and Slater were some well known hypergeometric  series 
identities. In physics Foda and Quano \cite{FQ:1996} observed the remarkable 
fact that the finitized Bose-Fermi identities of the configuration sum of the 
ABF model are of the form \eqref{eq:baileypair}. Hence one can read off Bailey pairs
from this.  This fact has been used in \cite{AB:1997,BMS:1995,BMSW:1997,FQ:1996,
Wa:1997} to derive Bose-Fermi identities for some subset of the minimal models, 
superconformal models and higher level coset models. Berkovich, McCoy and Schilling 
explored this in ~\cite{BMS:1995} for the Minimal model $M(p-1,p)$. They 
calculated the Bailey pairs using the Bose-Fermi identity for the unitary 
minimal model $M(p-1,p)$ and used the specialization \eqref{eq:half_finite_spec} 
and the additional specialization:
\begin{equation}\label{eq:finite_spec}
\rho_1=\text{finite} ,\quad \rho_2=\text{finite}. 
\end{equation} 
to obtain the characters of $N=1,2$ superconformal models, hence demonstrated 
a Bailey flow between these CFTs. In this chapter we will extend their method to
the nonunitary minimal models $M(p,p')$ where $p,p'>0$ are any two coprime integers.

A useful way to obtain new Bailey pairs from old ones is the construction of
dual Bailey pairs. If $(\alpha_n,\beta_n)$ is a bilateral Bailey pair 
with respect to $a$, the \textbf{dual Bailey pair} $(A_n,B_n)$ is defined as  
\begin{equation}\label{eq:dualpair}
\begin{split}
A_n(a,q) &= a^nq^{n^2}\alpha_n(a^{-1},q^{-1}),\\
B_n(a,q) &= a^{-n}q^{-n^2-n}\beta_n(a^{-1},q^{-1}).
\end{split}
\end{equation} 
Then $(A_n,B_n)$ satisfies \eqref{eq:defbail} with respect to $a$.

\section{Bailey pairs from the minimal models $M(p,p')$}\label{sec:baileypairs}
In this section we derive new Bailey pairs using the Bose-Fermi identity for
the minimal model $M(p,p')$.

As shown by Foda and Quano~\cite{FQ:1996}, the Bose-Fermi character 
identities~\cite{BM:1994,BMS:1997,FW:2001,W:2002} for the minimal models $M(p,p')$
are of the form 
\begin{equation}\label{eq:bosefermi1}
B_{r(b),s}^{(p,p')}(L,b;q)=q^{-\mathcal{N}_{r(b),s}}F_{r(b),s}^{(p,p')}(L,b;q),
\end{equation}
where the bosonic side is given by
\begin{equation}\label{eq:brs1}
\begin{split}
B_{r(b),s}^{(p,p')}(L,b;q)
=\sum_{j=-\infty}^{\infty} & \Bigg( q^{j(jpp'+r(b)p'-sp)}\left[\begin{array}{c}
L\\ \frac{1}{2}(L+s-b)-jp' \end{array}\right]_q \\ 
&-q^{(jp-r)(jp'-s)}\left[\begin{array}{c}L\\ \frac{1}{2}(L-s-b)+jp' \end{array}\right]_q 
\Bigg).
\end{split}
\end{equation}
The function fermionic formula  $F_{r(b),s}^{(p,p')}(L,b;q)$  will be discussed in 
the next section. The normalization constant $\mathcal{N}_{r(b),s}$ is explicitly 
calculated in \cite{BMS:1997}. Since the explicit expression is not used any 
where in our calculations we will exclude the details.

We will now construct new Bailey pairs using \eqref{eq:bosefermi1}. For 
simplicity we are going to write $r$ for $r(b)$.  Let us set $L=2n+b-s+2x$ to 
rewrite the q-binomial coefficients in \eqref{eq:brs1},  
\begin{eqnarray*} 
\left[\begin{array}{c}2n+b-s+2x\\ n+x-p'j \end{array}\right]_q &=&
\frac{\left(q^{b-s+2x+1}\right)_{2n}}{\left(q\right)_{n-(p'j-x)}\left(q^{b-s+2x+1}\right)_
{n+(p'j-x)}},\\ 
 \left[\begin{array}{c}2n+b-s+2x\\ n+x-s+p'j \end{array}\right]_q &=&
\frac{\left(q^{b-s+2x+1}\right)_{2n}}{\left(q\right)_{n-(p'j-b-x)}\left(q^{b-s+2x+1}\right)_
{n+(p'j-b-x)}}.
\end{eqnarray*}
\noindent
Following~\cite{BMS:1995,FQ:1996}  we  rewrite (\ref{eq:bosefermi1}) as
\begin{equation*}
\begin{split}
q^{-\mathcal{N}_{r(b),s}}F_{r(b),s}^{(p,p')}(L,b;q)=\sum_{j=-\infty}^{\infty}
& \Bigg(q^{j(jpp'+rp'-sp)} \frac{\left(q^{b-s+2x+1}
\right)_{2n}}{\left(q\right)_{n-(p'j-x)}\left(q^{b-s+2x+1}\right)_{n+(p'j-x)}}\\
&-q^{(jp-r)(jp'-s)}\frac{\left(q^{b-s+2x+1}\right)_{2n}}{\left(q\right)
_{n-(p'j-b-x)}\left(q^{b-s+2x+1}\right)_{n+(p'j-b-x)}}\Bigg).
\end{split}
\end{equation*}
where $L=2n+b-s+2x$.
This is in the form of  \eqref{eq:defbail}, hence we can read off the bilateral 
Bailey pair relative to $a=q^{b-s+2x}$ 
\begin{equation}\label{eq:bail1-1}
\begin{split}
\alpha_n &= \begin{cases}
q^{j(jpp'+rp'-sp)} & \text{if $n=jp'-x$}\\
-q^{(jp-r)(jp'-s)} & \text{if $n=jp'-b-x$}\\
0 & \text{otherwise}
\end{cases}\\
\beta_n &= \frac{q^{-\mathcal{N}_{r,s}}}{\left(aq\right)_{2n}} 
F_{r,s}^{(p,p')}(2n+b-s+2x,b;q).
\end{split}
\end{equation}
where $x=\frac{L-2n-b+s}{2}$.

Using the definition \eqref{eq:dualpair} we calculate the  Bailey pair dual 
to \eqref{eq:bail1-1} relative to $a=q^{b-s+2x}$ and denote it by 
$(\hat{\alpha}_n,\hat{\beta}_n)$ where 
\begin{equation}\label{eq:dbail1-1}
\begin{split}
\hat{\alpha}_n &= \begin{cases}
q^{j^2p'(p'-p)-jp'(r-b)-js(p'-p)-x(b+x-s)} & \text{if $n=jp'-x$}\\
-q^{(jp'-s)(j(p'-p)+r-b)-x(b+x-s)} & \text{if $n=jp'-b-x$}\\
0 & \text{otherwise}
\end{cases}\\
\hat{\beta}_n &= \frac{q^{\mathcal{N}_{r,s}}}{\left(aq\right)_{2n}} 
a^n q^{n^2}F_{r,s}^{(p,p')}(2n+b-s+2x,b;q^{-1}).
\end{split}
\end{equation}
Inserting \eqref{eq:bail1-1} and \eqref{eq:dbail1-1} into the bilateral Bailey 
lemma yields
\begin{equation}\label{eq:ra1-1}
\begin{split}
&\sum_{n=0}^{\infty}(\rho_1)_n(\rho_2)_n(aq/\rho_1\rho_2)^n\frac{q^{-\mathcal{N}_
{r(b),s}}}{\left(aq\right)_{2n}} F_{r,s}^{(p,p')}(2n+b-s+2x,b;q)\\
&=\frac{(aq/\rho_1)_{\infty}(aq/\rho_2)_{\infty}}{(aq)_{\infty}
(aq/\rho_1 \rho_2)_{\infty}}
\sum_{j=-\infty}^{\infty} \Bigg(\frac{(\rho_1)_{jp'-x}(\rho_2)_{jp'-x}}{(aq/\rho_1)_
{jp'-x}(aq/\rho_2)_{jp'-x}}(aq/\rho_1\rho_2)^{jp'-x}\\ & \times q^{j(jpp'+rp'-sp)}
-\frac{(\rho_1)_{jp'-b-x}(\rho_2)_{jp'-b-x}}{(aq/\rho_1)_
{jp'-b-x}(aq/\rho_2)_{jp'-b-x}}\\ & \times (aq/\rho_1\rho_2)^{jp'-b-x} q^{(jp-r)(jp'-s)}
\Bigg)
\end{split}
\end{equation}
and
\begin{equation}\label{eq:ra1-2}
\begin{split}
&\sum_{n=0}^{\infty}(\rho_1)_n(\rho_2)_n(aq/\rho_1\rho_2)^n\frac{q^{\mathcal{N}_
{r(b),s}}}{\left(aq\right)_{2n}}a^n q^{n^2} F_{r,s}^{(p,p')}(2n+b-s+2x,b;q^{-1})\\
&=\frac{(aq/\rho_1)_{\infty}(aq/\rho_2)_{\infty}}{(aq)_{\infty}
(aq/\rho_1 \rho_2)_{\infty}}
\sum_{j=-\infty}^{\infty} \Bigg(\frac{(\rho_1)_{jp'-x}(\rho_2)_{jp'-x}}{(aq/\rho_1)_
{jp'-x}(aq/\rho_2)_{jp'-x}}(aq/\rho_1\rho_2)^{jp'-x}\\ & \times 
q^{j^2p'(p'-p)-jp'(r-b)-js(p'-p)-x(b+x-s)}
-\frac{(\rho_1)_{jp'-b-x}(\rho_2)_{jp'-b-x}}{(aq/\rho_1)_
{jp'-b-x}(aq/\rho_2)_{jp'-b-x}}\\ & \times (aq/\rho_1\rho_2)^{jp'-b-x} 
q^{(jp'-s)(j(p'-p)+r-b)-x(b+x-s)}
\Bigg).
\end{split}
\end{equation}
As in~\cite{BMS:1995}, we are going to consider different specializations of the 
parameters $\rho_1$ and $\rho_2$ in \eqref{eq:ra1-1} and \eqref{eq:ra1-2} to get
character identities for $N=1,2$ superconformal algebras.

\section{Fermionic formulas for $M(p,p')$} \label{sec:fermionic}
So far we have only considered the bosonic side of \eqref{eq:bosefermi1} explicitly. 
For the fermionic side we will consider two cases  $p<p'<2p$ and 
$p'>2p$ separately  with $p$ and $p'$ relatively prime 
and  $r,s$ being pure Takahashi length.

\subsection{Fermionic formula for $M(p,p')$ with $p<p'<2p$} 
We need to introduce a lot of notations to give the explicit fermionic formula and
we follow~\cite[Section 4]{BMSW:1997} here. The fermionic formula depends on 
the continued fraction decomposition
\begin{equation*}
\frac{p'}{p'-p}=1+\nu_0 +\cfrac{1}{\nu_1+ \cfrac{1}{\nu_2 +\cdots 
\cfrac{1}{\nu_{n_0}+2}}}.
\end{equation*}
Define $t_i=\sum_{j=0}^{i-1}\nu_j$ for $1\le i \le n_0+1$ and the fractional
level incidence matrix $\I_B$ and corresponding Cartan matrix $B$ as
\begin{equation*}
\begin{split}
(\I_B)_{j,k}&=\begin{cases}  
     \delta_{j,k+1}+\delta_{j,k-1} & \text{for $1\le j<t_{n_0+1},j\ne t_i$}\\
     \delta_{j,k+1}+\delta_{j,k}-\delta_{j,k-1} & \text{for $j= t_i, 1\le i 
                                               \le n_0-\delta_{\nu_{n_0},0}$}\\
     \delta_{j,k+1}+\delta_{\nu_{n_0},0}\delta_{j,k} & \text{for $j= t_{n_0+1}$}
\end{cases}\\
B&=2I_{t_{n_0+1}}-\I_B,
\end{split}
\end{equation*}
where $I_{n}$ is the identity matrix of dimension $n$. Recursively define
\begin{equation*}
\begin{split}
y_{m+1}&=y_{m-1}+(\nu_m+\delta_{m,0}+2\delta_{m,n_0}) y_m, \qquad y_{-1}=0,
\qquad y_0=1,\\
\yb_{m+1}&=\yb_{m-1}+(\nu_m+\delta_{m,0}+2\delta_{m,n_0}) \yb_m, \qquad \yb_{-1}=-1,
\qquad \yb_0=1.
\end{split}
\end{equation*}
Then the Takahashi length and truncated Takahashi length are given by
\begin{equation*}
\begin{array}{ll}
\ell_{j+1}=y_{m-1}+(j-t_m)y_m\\
\ellb_{j+1}=\yb_{m-1}+(j-t_m)\yb_m
\end{array}
\qquad \text{for $t_m<j\le t_{m+1}+\delta_{m,n_0}$ with $0\le m\le n_0$.}
\end{equation*}

Let us define the $t_{n_0+1}$-dimensional vectors $Q^{(j)}$ $(j=1,2,\cdots,t_{n_0+1} +1)$
which we will need to specify the parity of the summation variables in the fermionic formula. 
For $1\le i\le t_{n_0+1}$
and $0\le m \le n_0$ such that $t_m<j\le t_{m+1}+\delta_{m,n_0}$ the components of 
$Q^{(j)}$ are defined recursively as
\begin{equation}
Q_i^{(j)}=\begin{cases} 0 &\quad \text{for $j\le i\le t_{n_0+1}$},\\
                                 j-i &\quad  \text{for $t_m\le i< j$},\\
                                Q_{i+1}^{(j)}+Q_{t_{m'}+1}^{(j)}& \quad\text{for $t_{m'-1}\le i <t_{m'},
                                1\le m'\le m$}. 
              \end{cases}                  
\end{equation}
When $\nu_{n_0}=0$, so that $t_{n_0+1}=t_{n_0}$, we need to set the initial condition
$Q_{t_{n_0}+1}^{(t_{n_0}+1)}=0$.
Also define
$$Q_{\uu}=\sum_{j=1}^{t_{n_0+1}+1} u_jQ^{(j)},$$
and for $t_i<j\le t_{i+1}$,
 \begin{equation}
( A_{\uu,\vv})_j=\begin{cases}u_j&\quad \text{for $i$ odd},\\
                   v_j& \quad \text{for $i$ even}.
                   \end{cases}
 \end{equation} 

For $b=\ell_{\beta+1}$, $r(b)=\ellb_{\beta+1}$ with $t_{\xi}<\beta\le t_{\xi+1}
+\delta_{\xi,n_0}$ and $s=\ell_{\sigma+1}$ with $t_{\zeta}<\sigma\le t_{\zeta+1}
+\delta_{\zeta,n_0}$ the fermionic formula is given by
\begin{equation}\label{eq:fermi1-1}
F_{r,s}^{(p,p')}(L,b;q)=q^{k_{b,s}}\sum_{\m\equiv Q_{\uu+\vv}\pmod{2}}
q^{\frac{1}{4}\m^tB\m -\frac{1}{2}
A_{\uu,\vv}\m}\prod_{j=1}^{t_{n_0+1}}\left[ \begin{array}{c} n_j+m_j\\m_j \end{array} 
\right]_{q}^{'}
\end{equation}
where $k_{b,s}$ is a normalization constant and $\n,\m\in\Z^{t_{n_0+1}}$ such that
\begin{equation}\label{I-B}  
\n +\m=\frac{1}{2}\Big(\I_B\m+\uu+\vv+L \e_1\Big)
\end{equation}
with $\e_i$ the standard $i$-th basis element of $\Z^{t_{n_0+1}}$,
$\uu=\e_{\beta}-\sum_{k=\xi+1}^{n_0}\e_{t_k}$ and $\vv=\e_\sigma-\sum_{k=\zeta+1}^{n_0}
\e_{t_k}$. The notation $\m \equiv Q_{\uu+\vv}\pmod{2}$ stands for $m_j$ even when
 $(Q_{\uu+\vv})_j$ is even and $m_j$ is odd when $(Q_{\uu+\vv})_j$ is odd. 

The $q$-binomial is also defined for negative entries
\begin{equation*}
\left[ \begin{array}{c} n+m\\m \end{array}\right]_q^{'}=
\frac{(q^{n+1})_m}{(q)_m}.
\end{equation*}
Note that
\begin{equation}\label{eq:inverbinom}
\left[ \begin{array}{c} n+m\\m \end{array} 
\right]_{q^{-1}}^{'} = q^{-nm}\left[ \begin{array}{c} n+m\\m \end{array} 
\right]_q^{'}.
\end{equation}
In fact using \eqref{eq:inverbinom} we get the following dual form of the fermionic 
formula that will be useful later on
\begin{multline}\label{eq:fermi1-2}
F_{r,s}^{(p,p')}(L,b;q^{-1})=\\q^{-k_{b,s}}\sum_{\m\equiv Q_{\uu+\vv}}   
q^{\frac{1}{4}\m^tB\m -\frac{1}{2}Lm_1+\frac{1}{2}A_{\uu,\vv}\m-\frac{1}{2}\m^t(\uu+\vv)}
\prod_{j=1}^{t_{n_0+1}}\left[ 
\begin{array}{c} n_j+m_j\\m_j \end{array} \right]_q^{'}.
\end{multline}  
\subsection{Fermionic formula for $M(p,p')$ with $p'>2p$:} We use the fermionic formula 
$F_{r(b),s}^{(p,p')}(L,b;q)$ derived in \cite[section 10] {BMS:1997} with notations defined
in \cite[section 2, section 3] {BMS:1997}  .
In this case the fermionic formula depends on the continued fraction decomposition of
$$\frac{p'}{p}= \nu_0+1 +\cfrac{1}{\nu_1+ \cfrac{1}{\nu_2 +\cdots \cfrac{1}{\nu_{n_0}+2}}}$$
where $\nu_{n_0}\ge 0$ and $\nu_i\ge 1$.
Let $t_i$ for $1\le i \le n_0+1$ be the same as in the previous case along with $t_0=-1$.
Recursively define
\begin{equation*}
\begin{split}
y_{m+1}&=y_{m-1}+(\nu_m+\delta_{m,0}+2\delta_{m,n_0}) y_m,  \quad y_{-1}=0, \quad 
 y_0=1, \quad 0\le m \le n_0\\
z_{m+1}&=z_{m-1}+(\nu_{m+1}+2\delta_{m+1,n_0}) z_m,  \quad z_{-1}=0, \quad
 z_0=1, \quad 0\le m \le n_0-1.
\end{split}
\end{equation*}
Then the Takahashi length and truncated Takahashi length are given by
\begin{equation*}
\ell_{j+1}^{(m)}=\begin{cases} j+1 &\text{for} \quad m=0 \text{ and  $0\le j\le t_1$}\\
                                                      y_{m-1}+(j-t_m)y_m & \text{for $1\le m 
                                                      \le n_0$ and 
                                                        $1+t_m\le j\le t_{1+m}+\delta_{n,m}$}.
                               \end{cases}       
\end{equation*}
\begin{equation*}                                               
\tilde{\ell}_{j+1}^{(m)}=\begin{cases} z_{m-2}+(j-t_m)z_{m-1} &
                 \quad \text{for $1+t_m<j\le t_{m+1}+\delta_{m,n_0}$ with $1\le m\le n_0$}\\
                0&\quad \text{for $m=0$}.
                \end{cases}
\end{equation*}
Let us define the corresponding Cartan matrix $B$  in this situation. The nonzero
elements of the matrix $B$ are given by the $\nu_0\times \nu_0$ matrix
\begin{equation}
C_T^{-1}=\left( \begin{array}{ccccc}
                              1& 1&1 & \cdots & 1\\ 
			  1& 2& 2&\cdots & 2\\
			   1& 2& 3&\cdots & 3\\ 
			  \cdot &\cdot &\cdot &\cdots &\cdot \\
			  \cdot &\cdot &\cdot &\cdots &\cdot \\
			  1& 2 & 3&\cdots &\nu_0
			   \end{array}
         \right)\\
\end{equation}
as 
\begin{equation}
\begin{split}
&B_{i,j}=2\left(C_T^{-1}\right)_{i,j} \quad \text{ for $1\le i,j \le \nu_0$}\\
&B_{\nu_0+1,j}=B_{j,\nu_0+1}=j \quad \text{for $1\le j \le \nu_0$}\\
&B_{j,j}=\frac{\nu_0}{2}\delta_{j,\nu_0+1}+(1-\frac{1}{2}\sum_{i=2}^{t_{n_0}}\delta_{j,t_i})
\quad \text{for $\nu_0+1\le j\le t_{n_0+1}-1$}\\
&B_{t_{n_0+1},t_{n_0+1}}=1\\
&B_{j,j+1}=-\frac{1}{2}+\sum_{i=2}^{n_0}\delta_{j,t_i} \quad \text{for $j>\nu_0$}\\
&B_{j+1,j}=-\frac{1}{2} \quad \text{for $j>\nu_0$}
\end{split}
\end{equation}
Define the $t_{n_0+1}$-dimensional vector $\eb_k$  and $1+t_{n_0+1}$-dimensional
vector $\e_k$ by
\begin{equation}
(\eb_k)_j=\begin{cases} \delta_{j,k} & \quad \text{for} \quad 1\le k\le t_{n_0+1}\\
                          0 &\quad \text{for}\quad k=0,1+t_{n_0+1}
                          \end{cases}
 \end{equation}
 \begin{equation}
(\e_k)_j=\begin{cases} \delta_{j,k} & \quad \text{for} \quad 1\le k\le t_{n_0+1}+1\\
                          0 &\quad \text{for}\quad k=0.
                          \end{cases}
 \end{equation}
 
The $t_{n_0+1}$-dimensional vectors $\uub,\uub_+,\uub_-,V$ and  $\overline{E}_{a,b}$
are defined by
\begin{equation}
\begin{split}
\uub &=\uub_++\uub_-\\
\uub_+&=\sum_{i=1}^{\nu_0}\overline{u}_i\eb_i\\
\uub_- &=\sum_{i=\nu_0+1}^{t_{n_0+1}}\overline{u}_i\eb_i\\
V &=\sum_{i=1}^{\nu_0}i\eb_i\\
 \overline{E}_{a,b}&=\sum_{i=a}^{b}\eb_i
\end{split}
\end{equation}

We will write $A=A^{(b)}+A^{(s)}$ where the $t_{{n_0}+1}$-dimensional vectors
$A^{(b)}$ and $A^{(s)}$ are defined as
\begin{equation}
A_k^{(b)}=\begin{cases} -\frac{1}{2}u_k &\quad\text{for $k$ in an even, nonzero zone},\\
                                    0& \quad \text{otherwise}.
\end{cases}
\end{equation}
where $\uu$ is a $t_{{n_0}+1}+1$-dimensional vector. 
\begin{equation}
A_k^{(s)}=
\begin{cases} -\frac{1}{2}\uub(s)_k-\frac{1}{2}V^t \cdot \uub(s)
                        \delta_{k,\nu_0+1} &\quad \text{for $k$ in an odd zone},\\
                         -\frac{1}{2}\eb_k^tB\uub(s)_+& \quad \text{for $k$ in 
                         an even zone}.
\end{cases}
\end{equation}
where $t_{n_0+1}$-dimensional vector $\uub(s)$ is given by: for $1\le k\le t_{n_0+1}$
\begin{equation}
(\uub(s))_k=\begin{cases} \delta_{k,j_s}-\sum_{i=\mu_s+1}^{n_0}\delta_{k,t_i} &
                              \quad \text{for $1+t_{\mu_s}\le j_s\le t_{\mu_s+1}$ and $\mu_s\le n-1$},\\
                              \delta_{k,j_s} & \quad\text{for $t_{\mu}<j_s,\mu_s=n$}
                              \end{cases}
\end{equation}  

For $\n,\m \in  \bf{Z}^{t_{n_0+1}}$
 define
\begin{equation*}
\begin{split}
\tm&=(n_1,n_2,\cdots,n_{\nu_0},m_{\nu_0+1},m_{\nu_0+2},\cdots,m_{t_{n_0+1}})\\
\tn&=(m_1,m_2,\cdots,m_{\nu_0},n_{\nu_0+1},n_{\nu_0+2},\cdots,n_{t_{n_0+1}})
\end{split}
\end{equation*}
such that
\begin{equation}
\tn+\tm=(I_{t_{n_0}+1}-B)\tm+L\overline{E}_{1,\nu_0}+\frac{L}{2}\eb_{\nu_0+1}+
\frac{B}{2}\uub_{+}+\frac{1}{2}\eb_{\nu_0+1}(\uub_{+}^t.V)+\frac{1}{2}\uub_{-}
\end{equation}

For $b=\ell_{j_{\mu}+1}^{(\mu)}, r(b)=\delta_{\mu,0}+ \tilde{\ell}_{j_{\mu}+1}^{(\mu)}$ with $1+t_{\mu}
<j_{\mu}\le t_{\mu+1}+\delta_{\mu,n_0}$ and $s=\ell_{j_{\beta}+1}^{(\beta)}$ with $1+t_{\beta}
<j_{\beta}\le t_{\beta+1}+\delta_{\beta,n_0}$ the fermoinic formula can be written  as:

\begin{equation}\label{eq:fermip'>2p}
F_{r(b),s}^{(p,p')}(L,b;q)=q^{C_{b,s}}\sum_{\tm\equiv \overline{Q}_{\uu}\pmod{2}}
q^{\frac{1}{2}\tm^tB\tm +
A^t\tm}\prod_{j=1}^{t_{n_0+1}}\left[ \begin{array}{c} \tilde{n}_j+\tilde{m}_j\\\tilde{m}_j \end{array} 
\right]_{q}^{'}
\end{equation}
\noindent
where  the normalization constant
$C_{b,s}= C(j_{\mu})$ is defined in \cite[(8.33)]{BMS:1997} as
\begin{equation}
\begin{split}
C(j_{0})&= 0 \quad \text{for $1\le j_0\le t_1$}\\
C(j_{\mu})&= \frac{1}{2}(-1)^{\mu}+(j_{\mu}-t_{\mu}) \{-\frac{\nu_0+1}
                      {4}+\frac{3}{4}\theta(\mu \quad  \text{odd})+c(t_{\mu})\} +c(t_{\mu-1})\\
                     &\text{ for} \quad t_{\mu}+1\le j_{\mu} \le t_{\mu+1}+2\delta_{\mu,n_0} \quad
                     \text{and} \quad 1\le \mu \le n_0  
\end{split} 
\end{equation}
where $\theta(S)=1$ if $S$ is true and $\theta(S)=0$ if $S$ is false. 

Also in our case  $\uu=\e_{j_{\mu}}-\sum_{ i=\mu+1}^{n_0}\e_i$.  To explain the sum $\tm\equiv 
\overline{Q}_{\uu}\pmod{2}$ let us introduce the following notation.
For $1\le k\le t_{n_0+1}, t_{\mu_0}+\delta_{\mu_0,0}+1\le j\le t_{\mu_0+1}+
\delta_{\mu_0,n_0}$, for some $0\le \mu_0\le n_0$,
\begin{equation*}
w_k^{(j)}=\begin{cases} 0 &\quad \text{for $k\ge j$},\\
                                 j-k &\quad  \text{for $t_{\mu_0}\le k< j$},\\
                                w_{k+1}^{(j)}+w_{t_{\mu+1}+1}^{(j)}& 
                                \quad\text{for $t_{\mu}\le k <t_{\mu+1},
                                1\le \mu < \mu_0 $}. 
              \end{cases}                  
\end{equation*}
then define
\begin{equation*}
{\bf w}(u_{1+t_{n_0+1}},\uub)=\sum_{k=1}^{t_{n_0+1}}\eb_k\left(u_{1+t_{n_0+1}}
w_k^{(1+t_{n_0+1})}+
\sum_{j=1}^{t_{n_0+1}}w_k^{(j)}\overline{u}_j\right). 
\end{equation*}
With the notations above $\tm\equiv \overline{Q}_{\uu}\pmod{2}$
means $\tm_-\in 2\Z^{t_{n_0+1}}+{\bf w}_-(u_{1+ t_{n_0+1}},\uub)$, and $m_+ \in \Z^{\nu_0}$. 
This restriction makes sure that the entries of all $q$-binomials in 
\eqref{eq:fermip'>2p} are integers as long as $\uu\in \Z^{1+t_{n_0+1}}$.

Again using (\ref{eq:inverbinom}) we get the following useful fermionic formula:

\begin{equation}\label{eq:invfermip'>2p}
\begin{split}
F_{r(b),s}^{(p,p')}(L;q^{-1})&=q^{-C_{b,s}}\sum_{\tm \equiv \overline{Q}_{\uu}}   
q^{\frac{1}{2}\tm^t B\tm-L\tm^t \overline{E}_{1,\nu_0}-\frac{L}{2}\tm^t \e_{\nu_0+1}}\\
&\times q^{-A^t\tm-\frac{1}{2}\tm^tB\uu_+
-\frac{1}{2}\tm^t\e_{\nu_0+1}(\uu_{+}^t.V)-\frac{1}{2}\tm^t\uu_{-}}
\times \prod_{j=1}^{t_{n_0+1}}\left[ 
\begin{array}{c} \tilde{n}_j+\tilde{m}_j\\\tilde{m}_j \end{array} \right]_q^{'}
\end{split}
\end{equation}  
\noindent             
We will use this in the later sections. 
\begin{remark}
We like to mention that the fermionic formula described above  
for $p<p'<2p$ is not valid for $r=b=1$ and the fermionic formula for $p'>2p$
is not valid for $r=0,b=1$. The valid fermionic formulas in these special cases are
given in \cite{W:2002}. Let us state the formula in the case $r=b=1$ and $p<p'<2p$ which
will be used later.
\end{remark}

\subsection{Special fermionic formula for $r=b=1$ and $p<p'<2p$:} Here we state the 
fermionic formula when $r=b=1$ and $p<p'<2p$ as given in Welsh's paper \cite{W:2002}. 
Let us introduce the notations following \cite{W:2002}. 

In this case the normalization constant in \eqref{eq:bosefermi1} is $\mathcal{N}_{1,1}=0$. 

The formula depends on
the continued fraction decomposition given by 
$$\frac{p'}{p}= c_0 +\cfrac{1}{c_1+ \cfrac{1}{c_2 +\cdots \cfrac{1}{c_{n_0}}}}$$
where $c_{n_0}\ge 2$ and $c_i\ge 1$ for $0\le i< n_0$.
Define $t_i=-1+\sum_{j=0}^{i-1}c_j$ for $0\le i \le n_0+1$ and let $t=-2+\sum_{i=0}^{n_0}c_i$.

Let us define the $t \times t$ matrices $\CC=(C_{ji})_{0\le i,j< t} $ and $\CC^*=(C_{ji})_
{\substack{1\le j\le t\\0\le i <t}} $ where
for $0\le i,j\le t$,
\begin{equation}
\begin{split}
C_{j,j-1}&=-1, \quad C_{j,j}=1, \quad C_{j,j+1}=1, \quad \text{if $j=t_k$ for $k=1,2,\cdots, n_0$};\\
C_{j,j-1}&=-1, \quad C_{j,j}=2, \quad C_{j,j+1}=-1, \quad 0\le j <t, \quad \text{otherwise};\\
C_{ji}&=0 \quad \text{ for $|i-j|>1$}.
\end{split}
\end{equation} 
The $(t-t_1-1)\times (t-t_1-1)$ matrix $\Cb$ is defined as $\Cb=(C_{ji})_{t_1+1\le i,j< t}$.
Note that when $p<p'<2p$ we have $t_1=0$, hence $\Cb$ is a $t-1 \times t-1$ matrix.
In this section all the vectors are column vectors, but we write them as row vectors.

Let $\m=(m_1,\cdots, m_{t-1})$ be a $t-1$-dimensional vector $\hat{\m}=(L,m_1,\cdots, 
m_{t-1})$ and $\n=(n_1,\cdots, n_{t})$ be two $t$-dimensional
vectors satisfying
\begin{equation}\label{mn} 
\n=-\frac{1}{2}\CC^*\hat{\m}+ \frac{1}{2}\uu
\end{equation}
where  $\uu=(u_1,\cdots, u_t)=2\uu'$ and
$\uu'=-\sum_{i=1}^{n}\e_{t_i}$ with $\e_j=(\delta_{1j},\delta_{2j}, \cdots, \delta_{tj})$, for
$1\le j\le t$. Note that for $1\le j\le t-1$,
\begin{equation}
(\CC^*\hat{\m})_j=(\Cb\m-L\overline{\e}_1)_j
\end{equation}
where $\overline{\e}_1=(1,0\cdots,0)$ is the $t-1$-dimensional vector. 

Using $\uu'$ let us define $t-1$-dimensional vectors $\uu'_{\flat}$ and $\uu'_{\sharp}$ as
\begin{equation*}
(u'_{\flat})_j=\begin{cases} 0 & \quad \text{if $t_k<j\le t_{k+1}, k\equiv 0 \pmod {2}$};\\
                               {u}_j' & \quad \text{if $t_k<j\le t_{k+1}, k\not\equiv 0 (\mod 2)$};
                 \end{cases}                   
\end{equation*} 
and 
\begin{equation*}
(u'_{\sharp})_j=\begin{cases} u_j' & \quad \text{if $t_k<j\le t_{k+1}, k\equiv 0 \pmod {2}$};\\
                               0 & \quad \text{if $t_k<j\le t_{k+1}, k\not\equiv 0 (\mod 2)$}.
                 \end{cases}                   
\end{equation*}  
Here $\uu'=\uu'_{\flat}+\uu'_{\sharp}$.

Note that $\CC^*$ is invertible. We define $Q_i\in\{0,1\}$ for $0\le i<t$ by
\begin{equation}
(Q_0,Q_1,\cdots, Q_{t-1})\equiv (\CC^*)^{-1}\uu.
\end{equation}
Using this we define the $t-1$-dimensional parity vector $Q=(Q_1,Q_2, \cdots, Q_{t-1})$.

The fermionic formula is given by
\begin{equation} \label{eq:ferm_r=b=1}
F_{1,1}^{p,p'}(L;q)= \sum_{\hat{\m}\equiv Q(\uu)}q^{\frac{1}{4}\m^t\overline{C}\m
-\frac{1}{2}(\uu'_{\flat}+\uu'_{\sharp})\m+\frac{1}{4}\gamma(\uu')}
\prod_{j=1}^{t-1} \left[ \begin{array}{c} m_j+n_j\\m_j
\end{array}\right]_q+ F_{1,1}^{\hat{p},\hat{p}'}(L;q) 
\end{equation}      
where $\hat{\m}\equiv Q(\uu)$ means $L\equiv Q_0 \pmod 2$ and  $m_i\equiv Q_i \pmod{2}$ for
 $1\le i\le t-1$. This ensures that $m_j+n_j$ is an integer. Let us define the term $\gamma(\uu')$.
 We iteratively generate sequences
 $(\beta_t,\beta_{t-1},\cdots, \beta_0)$, $(\alpha_t, \alpha_{t-1},\cdots,\alpha_0)$ and
 $(\gamma_t,\gamma_{t-1},\cdots, \gamma_0)$ as follows. Let $\alpha_t=\beta_t=\gamma_t=0$.
 Now, for $j\in\{t,t-1,\cdots,1\}$, obtain $\alpha_{j-1},\beta_{j-1}$, and $\gamma_{j-1}$ in the following
 three steps. First let
 \begin{equation*}
 (\beta'_{j-1},\gamma'_{j-1})=(\beta_j,  -2\alpha_j u'_j).
 \end{equation*}  
Then set 
\begin{equation*}
 (\alpha''_{j-1},\gamma''_{j-1})=(\alpha_j+\beta'_{j-1}, \gamma'_{j-1} -(\beta'_{j-1})^2).
 \end{equation*}  
Finally set
\begin{equation*}
 (\alpha_{j-1},\beta_{j-1},\gamma_{j-1})=
 \begin{cases}
 (\alpha''_{j-1},\alpha_j,  -(\alpha''_{j-1})^2-\gamma''_{j-1}) &\text{ if $j=t_k+1$ with 
 $1\le k\le n_0$},\\
 (\alpha''_{j-1},\beta'_{j-1},\gamma''_{j-1}) \quad \text{ otherwise}.
 \end{cases}
 \end{equation*}  
Now define $\gamma(\uu')=\gamma_0$.

To define $\hat{p}$ and $\hat{p}'$ let us first define $y_k$ and $z_k$ for $-1\le k \le n_0+1$.
 \begin{equation*}
 \begin{split}
  y_{-1}=0, \quad y_0=1, \quad y_{k}&=y_{k-2}+c_{k-1} y_{k-1}, \quad \text{for $1\le k\le n_0+1$},\\
  z_{-1}=1, \quad z_0=0, \quad z_{k}&=z_{k-2}+c_{k-1} z_{k-1}, \quad \text{for $1\le k\le n_0+1$}.
\end{split}
\end{equation*}
Define $\{\xi_i\}_{i=0}^{2c_{n_0}-1}$ and $\{\tilde{\xi}_i\}_{i=0}^{2c_{n_0}-1}$ according to
$\xi_0=\tilde{\xi_0}=0$, $\xi_{2c_{n_0}-1}=p'$, $\tilde{\xi}_{2c_{n_0}-1}=p$ and
\begin{equation*}
\begin{split}
&\xi_{2k-1}=ky_{n_0}; \quad \quad\quad \quad \tilde{\xi}_{2k-1}=kz_{n_0};\\
&\xi_{2k}=ky_{n_0}+y_{n_0-1}; \quad \quad \tilde{\xi}_{2k}=kz_{n_0}+z_{n_0-1},
\end{split}
\end{equation*}
for $1\le k \le c_{\n_0}$. Now for $1\le a<p'$ define $\eta(a)=\ell$ where $\xi_\ell \le a<\xi_{\ell+1}$.
Similarly, for $1\le b<p$ define $\tilde{\eta}(b)=\ell$ where $\tilde{\xi}_\ell \le b<\tilde{\xi}_{\ell+1}$.
Then the values of $\hat{p}$ and $\hat{p}'$ that appeared in \eqref{eq:ferm_r=b=1} are given by
 $\hat{p}'=\xi_{\eta(b)+1}-\xi_{\eta(b)}$ and $\hat{p}=\tilde{\xi}_{\eta(b)+1}-\xi_{\eta(b)}$. Note that we have $a=b=1$.

It is easy to show that 
 \begin{equation}\label{eq:ferminverse_r=b=1}
\begin{split}
&F_{1,1}^{p,p'}(L;q^{-1})=\\
&\sum_{\hat{\m}\equiv Q(\uu)}q^{\frac{1}{4}\m^t\overline{C}\m-\frac{1}{2}Lm_1+\frac{1}{2}
(\uu'_{\flat}+\uu'_{\sharp})\m-\frac{1}{2}\overline{\uu} \m-\frac{1}{4}\gamma(\uu')}
\prod_{j=1}^{t-1} \left[ \begin{array}{c}m_j+n_j\\m_j
\end{array}\right]_q +F_{1,1}^{\hat{p},\hat{p}'}(L;q^{-1}).
\end{split}
\end{equation} 
where $\overline{\uu}=(u_1,\cdots,u_{t-1})$ is the $t-1$-dimensional vector.


\section{$N=1$ Superconformal character from $M(p,p')$} \label{sec:N1}
The $N=1$ superconformal algebra is the infinite dimensional Lie 
super algebra~\cite{NS:1971,Ra:1971,Ki:1987} with basis $L_n,G_r,\tilde{C}$ and 
(anti)-commutation relation given by
\begin{equation*}
\begin{split}
\left[L_m,L_n\right] &= (m-n)L_{m+n}+\frac{\tilde{C}}{8}(m^3-m)\delta_{m+n,0}\\
\left[L_m,G_r\right] &= (\frac{m}{2}-r)G_{m+r} \\
\{G_r,G_s\} &= 2L_{r+s}+\frac{\tilde{C}}{2}(r^2-\frac{1}{4})\delta_{r+s,0}
\end{split}
\end{equation*}
where $m,n$ are integers $\tilde{C}$ is the central charge and its eigen value is parametrized by 
$\tilde{c}=\frac{3}{2}-\frac{3(p-p')^2}{pp'}$. If $r,s$ are integers the algebra is called Neveu-Schwarz (NS)
algebra and if $r,s$ are half integers then the algrebra is called Ramond (R) algebra.
Let us denote these algebras by $SM(p,p')$.
 
The character formula of these algebras are calculated in  ~\cite{Dob:1987,GKO:1986} and are
given by,
\begin{equation}\label{eq:N1char}
\tilde{\chi}^{(p,p')}_{r,s}(q)=\tilde{\chi}^{(p,p')}_{p-r,p'-s}(q)
=\frac{(-q^{\epsilon_{r-s}})_{\infty}}{(q)_{\infty}}
\sum_{j=-\infty}^{\infty}\Big(q^{\frac{j(jpp'+rp'-sp)}{2}}-q^{\frac{(jp-r)(jp'-s)}{2}}
\Big),
\end{equation}
where $1\le r\le p-1,1\le s\le p'-1$, $p$ and $(p'-p)/2$ are relatively prime and
\begin{equation}
 \epsilon_i=
 \begin{cases} \frac{1}{2}& \text{if $i$ is even (NS-sector),}\\
            1 & \text{ if $i$ is odd (Ramond-sector).}
\end{cases}              
\end{equation}

In this section we are going to consider the specialization of the form of 
\eqref{eq:half_finite_spec} in~\eqref{eq:ra1-1} and~\eqref{eq:ra1-2}.
We will see that these give Bailey flows from the minimal model $M(p,p')$ to the
superconformal models $SM(p',2p+p')$ and $SM(p',3p'-2p)$. 

\subsection{The model $SM(p',2p+p')$}

Specializing $\rho_1\longrightarrow \infty$ and $\rho_2=-q^{\frac{b-s+1}{2}}$ with 
$x=0$ in \eqref{eq:ra1-1}  and comparing with \eqref{eq:N1char} we find for $b-s$ even (NS-sector)
\begin{equation}\label{eq:N1NS-1}
\tilde{\chi}^{(p',2p+p')}_{s,2r+b}(q)= \sum_{n\ge 0} \frac{q^{\frac{1}{2}(n^2+nb-ns)}
(-q^{\frac{1}{2}})_{n+(b-s)/2}}{(q)_{2n+b-s}}q^{-\mathcal{N}_{r,s}}F_{r,s}^{(p,p')}(2n+b-s,b;q)
\end{equation}
and for $b-s$ odd (Ramond-sector)
\begin{equation}\label{eq:N1R-1}
\tilde{\chi}^{(p',2p+p')}_{s,2r+b}(q)=\sum_{n\ge 0} \frac{q^{\frac{1}{2}(n^2+nb-ns)}
(-q)_{n+(b-s-1)/2}}{(q)_{2n+b-s}}q^{-\mathcal{N}_{r,s}}F_{r,s}^{(p,p')}(2n+b-s,b;q).
\end{equation}

Hence there is a Bailey flow from $M(p,p')$ to the
superconformal model $SM(p',2p+p')$. Let us calculate the fermionic formula 
using section~\ref{sec:fermionic}.\\

\noindent 
{\bf Fermionic formula for $SM(2p+p',p')$ with $p<p'<2p$:}
To obtain an explicit fermionic formula set $m_0=L=2n+b-s$ and insert \eqref{eq:fermi1-1}
into \eqref{eq:N1NS-1}. Then using
\begin{equation}\label{eq:-q1/2}
(-q^{\frac{1}{2}})_{\frac{m_0}{2}}=\sum_{k=0}^{\frac{m_0}{2}}q^{\frac{1}{2}(\frac{m_0}{2}-k)^2}
\left[ \begin{array}{c}\frac{m_0}{2} \\ k \end{array} \right]_q
\end{equation}
we find
\begin{equation}\label{eq:fermiN1NS-1}
\begin{split}
\tilde{\chi}^{(p',2p+p')}_{s,2r+b}(q) & = q^{-\frac{1}{8}(b-s)^2-\mathcal{N}_{r,s}+k_{b,s}}
 \sum_{\substack{m_0=0\\ \text{$m_0$ even}}}^{\infty}\sum_{k=0}^{\frac{m_0}{2}}
 \sum_{\m\equiv Q_{\uu+\vv}}q^{\frac{1}{8}m_0^2+\frac{1}{2}(\frac{m_0}{2}-k)^2}\\ 
& \times q^{\frac{1}{4}\m^tB\m-\frac{1}{2}A_{\uu,\vv}\m} 
 \times \frac{1}{(q)_{m_0}} \left[ \begin{array}{c}\frac{m_0}{2} \\ k \end{array} \right]_q
 \prod_{j=1}^{t_{n_0+1}}\left[ 
\begin{array}{c} n_j+m_j\\m_j \end{array} \right]_q^{'}.
\end{split}
\end{equation}

Setting $\p=(k,m_0,\m) \in \Z^{t_{n_0+1}+2}$  we  can write
\begin{equation}
\frac{1}{8}m_0^2+\frac{1}{2}(\frac{m_0}{2}-k)^2+\frac{1}{4}\m^tB\m=\frac{1}{4}\p^t\tilde{B}\p
\end{equation}
where
\begin{equation}
\tilde{B} = \left( \begin{array}{cc|cc}
                            2  & -1 & 0\\ 
			    -1 & 1 & 1\\\hline  
			    0 & -1 & B \\
			    \end{array}
         \right)\\
\end{equation}         
Using this the NS-sector character  \eqref{eq:fermiN1NS-1} 
can be rewritten as  
\begin{equation}\label{eq:fermiN1NS-2}
\begin{split}
\tilde{\chi}^{(p',2p+p')}_{s,2r+b}(q) &=q^{-\frac{1}{8}(b-s)^2-\mathcal{N}_{r,s}+k_{b,s}}
 \sum_{\substack{\p\in\Z^{t_{n_0+1}+2}\\p_i \equiv (\tilde{Q}_{\uu,\vv})_i,i\ge 2}}
 q^{\frac{1}{4}\p^t\tilde{B}\p-\frac{1}{2}\tilde{A}_{\uu,\vv}\p}\\
 & \times \frac{1}{(q)_{p_2}} \prod_{\substack{j=1,j\ne 2}}^{t_{n_0+1}+2}\left[ 
\begin{array}{c} \frac{1}{2}(\I_{\tilde{B}}\p+\tilde{\uu}+\tilde{\vv})_j \\ p_j \end{array} 
\right]_q^{'}
\end{split}
\end{equation}
where $\I_{\tilde{B}}=2I_{t_{n_0+1}+2}-\tilde{B}$,
\begin{equation}\label{eq:C-1}
\begin{split}
\tilde{A}_{\uu,\vv}&=(0,0,A_{\uu,\vv}),\\
\tilde{\uu}^t&=(0,0,\uu^t),\\
\tilde{\vv}^t&=(0,0,\vv^t),\\
\tilde{Q}_{\uu+\vv}^t&=(0,0,Q_{\uu+\vv}^t).
\end{split}
\end{equation}

Similarly setting $m_0=2n+b-s$ in (\ref{eq:N1R-1}) and using 
\begin{equation}\label{eq:-q}
(-q)_{\frac{m_0-1}{2}}=\frac{1}{2}\sum_{k=0}^{\frac{m_0+1}{2}}q^{\frac{1}{2}(\frac{m_0+1}{2}-k)
(\frac{m_0-1}{2}-k)}
\left[ \begin{array}{c}\frac{m_0+1}{2} \\ k \end{array} \right]_q
\end{equation}

we get the fermionic formula in the Ramond-sector,
\begin{multline}\label{eq:fermiN1R-1}
\tilde{\chi}^{(p',2p+p')}_{s,2r+b}(q)  =\frac{1}{2}
 q^{-\frac{1}{8}((b-s)^2+1)-\mathcal{N}_{r,s}+k_{b,s}}
 \sum_{\substack{\p \in \Z^{t_{n_0+1}+2}\\ p_i\equiv (\tilde{Q}_{\uu,\vv})_i, i\ge 2}}
 q^{\frac{1}{4}\p^t\tilde{B}\p-\frac{1}{2}\tilde{A}_{\uu,\vv}\p}\\
  \times \frac{1}{(q)_{p_2}}  \prod_{\substack{j=1,j\ne 2}}^{t_{n_0+1}+2}\left[ 
\begin{array}{c} \frac{1}{2}(\I_{\tilde{B}}\p+\tilde{\uu}+\tilde{\vv})_j \\ p_j \end{array} 
\right]_q^{'}
\end{multline}
where $\tilde{B},\tilde{A},\tilde{\vv}$ are as in \eqref{eq:C-1} and 
$\tilde{\uu}^t=(1,0,\uu^t)$, $\tilde{Q}_{\uu+\vv}^t=(0,1,Q_{\uu+\vv}^t)$.\\

\noindent
{\bf Fermionic formula for $SM(2p+p',p')$ with $p'>2p$:} In this section we will just state the 
formulas without showing the calculations. The calculation is very similar to the previous case.\\

\noindent
Using (\ref{eq:fermip'>2p}) in NS-sector we get,
\begin{equation}\label{eq:fermip'>2pN1NS-1}
\begin{split}
\tilde{\chi}^{(2p+p',p')}_{2r+b,s}(q) & =q^{-\frac{1}{8}(b-s)^2-\mathcal{N}_{r,s}+C_{b,s}}
 \sum_{\p\equiv \hat{\overline{Q}}_{\uu}}q^{\frac{1}{4}\p^tT\p-\frac{1}{2}\hat{A}^t\p}
 \times \frac{1}{(q)_{p_2}} \left[ \begin{array}{c}\frac{p_2}{2} \\ p_1 \end{array} \right]_q \\
 & \times \prod_{j=3}^{t_{n_0+1}+2}\left[ \begin{array}{c} \frac{1}{2}(\I_T\p+\hat{B}
\hat{\uu}_{+}+\hat{\e}_{\nu_0+1}(\hat{\uu}_+^t.\hat{V})+\hat{\uu}_{-})_j \\ p_j \end{array} 
\right]_q^{'}
\end{split}
\end{equation} 
where
\begin{eqnarray}\label{eq:T-1}
\p&=& (k,m_0,\tm) \in Z^{t_{n_o+1}+2}\nonumber\\
T & = & \left( \begin{array}{cc|cccccccccccc}
                            2 & -1 &0&.&.&.&.&.&.&.&0\\ 
			                     -1 & 1 & 2&.&.&2&1&0&.&.&0\\\hline  
			                      0 &-2 & &&&&&&&& \\
			                      . & . & &&&&&&&& \\
			                      . & . & &&&&&&&& \\
			                      . &-2 & &&&&2B&&&& \\
			                      . &-1 & &&&&&&&& \\
			                      . &0 & &&&&&&&& \\
			                      . &. & &&&&&&&& \\
			                      . &. & &&&&&&&& \\
			                      0 &0 & &&&&&&&& \\
			                      
			    \end{array}
         \right)\nonumber\\
         \I_T& = & 2I_{t_{n_o+1}+2}-T \nonumber\\
 \hat{X}^t&=& (0,0,X^t) \text{ for $X \in Z^{t_{n_o+1}}$ }
\nonumber\\
 \hat{B}&=& \left( \begin{array}{cc|cc}
                            0 & 0 &0\\ 
			                      0 & 0 &0\\\hline  
			                      0&0 &B \\
			                      \end{array}
                           \right)\nonumber\\
\p\equiv \hat{\overline{Q}}_{\uu}&:=&\begin{cases} & 0 \le k\le \frac{m_0}{2}\\
                                            & m_0 \ge 0, \text{$m_0$ is even}\\
                                           & \tm \equiv \overline{Q}_{\uu}
                              \end{cases} \nonumber
\end{eqnarray}
\noindent
Note $T$ is a $(t_{n_o+1}+2\times t_{n_o+1}+2)$ matrix, number of 2 in second row is $\nu_0$ and  
number of -2 in second column is $\nu_0$.\\

\noindent
In Ramond-sector we get,
\begin{equation}\label{eq:fermip'>2pN1NS-2}
\begin{split}
\tilde{\chi}^{(2p+p',p')}_{2r+b,s}(q) & =q^{-\frac{1}{8}((b-s)^2+1)-\mathcal{N}_{r,s}+C_{b,s}}
 \sum_{\p\equiv \hat{\overline{Q}}_{\uu}}q^{\frac{1}{4}\p^tT\p-\frac{1}{2}\hat{A}^t\p}
 \times \frac{1}{(q)_{p_2}} \left[ \begin{array}{c}\frac{p_2+1}{2} \\ p_1 \end{array} \right]_q \\
 & \times \prod_{j=3}^{t_{n_0+1}+2}\left[ \begin{array}{c} \frac{1}{2}(\I_T\p+\frac{\hat{B}}{2}
\hat{\uu}_{+}+\frac{1}{2}\hat{\e}_{\nu_0+1}(\hat{\uu}_+^t.\hat{V})+\frac{1}{2}\hat{\uu}_{-})_j 
\\ p_j \end{array} \right]_q^{'}
\end{split}
\end{equation} 
\noindent
where $T,\hat{A},\hat{\uu}, \hat{B}$ are as in (\ref{eq:T-1}) and  
\begin{equation*}
\p\equiv \hat{\overline{Q}}_{\uu}:=\begin{cases} & 0 \le k\le \frac{m_0+1}{2}\\
                                            & m_0 \ge 0, \text{$m_0$ is odd}\\
                                           & \tm \equiv \overline{Q}_{\uu}
                              \end{cases}  
\end{equation*} 

\subsection{The model $SM(p',3p'-2p)$} Similarly using the same specialization with the dual 
Bailey pair in \eqref{eq:ra1-2}  and comparing the bosonic side with \eqref{eq:N1char} with 
we find for $b-s$ even in the NS-sector
\begin{equation}\label{eq:N1NS-2}
\tilde{\chi}^{(p',3p'-2p)}_{s,3b-2r}(q)=
\sum_{n\ge 0} \frac{q^{\frac{3n}{2}(n+b-s)}(-q^{\frac{1}{2}})_{n+(b-s)/2}}{(q)_{2n+b-s}}  
              q^{\mathcal{N}_{r,s}}F_{r,s}^{(p,p')}(2n+b-s,b;q^{-1}) 
\end{equation}
and for $b-s$ odd in the Ramond-sector
\begin{equation}\label{eq:N1R-2}
\tilde{\chi}^{(p',3p'-2p)}_{s,3b-2r}(q)=
\sum_{n\ge 0} \frac{q^{\frac{3n}{2}(n+b-s)}(-q)_{n+(b-s-1)/2}}{(q)_{2n+b-s}}  
              q^{\mathcal{N}_{r,s}}F_{r,s}^{(p,p')}(2n+b-s,b;q^{-1}).
\end{equation}  \\

{\bf Fermionic formula for $SM(p',3p'-2p)$ with $p<p'<2p$:}
To obtain the fermionic formula, as before we are going to set $m_0=2n+b-s$.
Inserting \eqref{eq:-q} and \eqref{eq:fermi1-2} into \eqref{eq:N1R-2}
we get in the Ramond-sector
\begin{equation}\label{eq:fermiN1R-2}
\begin{split}
\tilde{\chi}^{(p',3p'-2p)}_{s,3b-2r}(q) &= \frac{1}{2}q^{-\frac{1}{8}(3(b-s)^2+1)+ 
\mathcal{N}_{r,s}-k_{b,s}}
 \sum_{\substack{m_0=0\\ \text{$m_0$ odd}}}^{\infty}\sum_{k=0}
^{\frac{m_0+1}{2}}\sum_{\m\equiv Q_{\uu+\vv}}\\
&\times q^{\frac{1}{2}(m_0^2+k^2-m_0k-m_0m_1)} 
  q^{\frac{1}{4}\m^tB\m-\frac{1}{2}\m^t(\uu+\vv)+\frac{1}{2}A_{\uu,\vv}\m}\\ 
& \times \frac{1}{(q)_{m_0}} \left[ \begin{array}{c}\frac{m_0+1}{2} \\ k \end{array} \right]_q
 \prod_{j=1}^{t_{n_0+1}}\left[ 
\begin{array}{c} n_j+m_j\\m_j \end{array} \right]_q^{'}.
\end{split}
\end{equation}
Define $\p=(k,m_0,\m) \in \Z^{t_{n_0+1}+2}$, so that \eqref{eq:fermiN1R-2} in the Ramond-sector
can be rewritten as
\begin{multline}\label{eq:fermiN1R-3}
\tilde{\chi}^{(p',3p'-2p)}_{s,3b-2r}(q) =\frac{1}{2}
 q^{-\frac{1}{8}(3(b-s)^2+1)+\mathcal{N}_{r,s}-k_{b,s}}
 \sum_{\substack{\p\in \Z^{t_{n_0+1}+2}\\ p_i\equiv (\tilde{Q'}_{\uu+\vv})_i, i\ge 2}}
q^{\frac{1}{4}\p^t\tilde{B}'\p+\frac{1}{2}\tilde{A}_{\uu,\vv}\p}\\
 \times \frac{1}{(q)_{p_2}} \prod_{\substack{j=1,j\ne 2}}^{t_{n_0+1}+2}\left[ 
\begin{array}{c} \frac{1}{2}(\I_{\tilde{B}'}\p+\tilde{\uu}+\tilde{\vv})_j \\ p_j \end{array} 
\right]_q^{'}
\end{multline}
where $\I_{\tilde{B}'}=2I_{t_{n_0+1}+2}-\tilde{B}'$, $\tilde{\vv}$ as in \eqref{eq:C-1}, 
$\tilde{\uu}^t =(1,0,\uu^t)$, $(\tilde{Q}_{\uu+\vv}')^t =(0,1,Q_{\uu+\vv}^t)$, and
\begin{equation}\label{eq:C'-1}
\begin{split}
\tilde{B}'&=\left( \begin{array}{cc|cc}
                            2  & -1 & 0\\ 
			    -1 & 2 & -1\\\hline  
			    0 & -1 & B \\
			    \end{array}
         \right)\\
\tilde{A}_{\uu,\vv}&=(0,0,A_{\uu,\vv}-\uu^t-\vv^t).
\end{split}
\end{equation}

Similarly, for the NS-sector it follows from \eqref{eq:N1NS-2}
\begin{multline}\label{eq:fermiN1NS-3}
\tilde{\chi}^{(p',3p'-2p)}_{s,3b-2r}(q) =q^{-\frac{3}{8}(b-s)^2+\mathcal{N}_{r,s}-k_{b,s}}
 \sum_{\substack{\p\in\Z^{t_{n_0+1}+2}\\ p_i\equiv (\tilde{Q'}_{\uu+\vv})_i, i\ge 2}}
 q^{\frac{1}{4}\p^t\tilde{B}'\p+\frac{1}{2}\tilde{A}_{\uu,\vv}\p}\\
  \times \frac{1}{(q)_{p_2}} \prod_{\substack{j=1,j\ne 2}}^{t_{n_0+1}+2}\left[ 
\begin{array}{c} \frac{1}{2}(\I_{\tilde{B}'}\p+\tilde{\uu}+\tilde{\vv})_j \\ p_j \end{array} 
\right]_q^{'}
\end{multline}
with $\tilde{B}'$ and $\tilde{A}_{\uu,\vv}$ as in \eqref{eq:C'-1},
$(\tilde{Q}_{\uu+\vv}')^t=(0,0,Q_{\uu+\vv}^t)$, $\tilde{\uu}^t =(0,0,\uu^t)$
and $\tilde{\vv}^t =(0,0,\vv^t)$.\\

{\bf Fermionic formula for $SM(p',3p'-2p)$ with $p'>2p$:} Again we state the formulas without 
showing the calculations. \\

\noindent
Using (\ref{eq:invfermip'>2p}) in NS-sector we get,
\begin{equation}\label{eq:fermip'>2pN1NS-3}
\begin{split}
\tilde{\chi}^{(p',3p'-2p)}_{s, 3b-2r}(q) & =q^{-\frac{3}{8}(b-s)^2+\mathcal{N}_{r,s}-C_{b,s}}
 \sum_{\p\equiv \hat{\overline{Q}}_{\uu}}q^{\frac{1}{4}\p^tT'\p}\\
 &\times q^{-\hat{A}^t\p-\frac{1}{2}\p^t\hat{B}\hat{\uu}_+
-\frac{1}{2}\p^t\hat{\e}_{\nu_0+1}(\hat{\uu}_+^t.\hat{V})-\frac{1}{2}\p^t\hat{\uu}_{-}}
 \times \frac{1}{(q)_{p_2}} \left[ \begin{array}{c}\frac{p_2}{2} \\ p_1 \end{array} \right]_q \\
 & \times \prod_{j=3}^{t_{n_0+1}+2}\left[ \begin{array}{c} \frac{1}{2}(\I_{T'}\p+
\frac{\hat{B}}{2}\hat{\uu}_{+}+\frac{1}{2}\hat{\e}_{\nu_0+1}(\hat{\uu}_+^t.\hat{V})+
\frac{1}{2}\hat{\uu}_{-})_j \\ p_j \end{array} \right]_q^{'}
\end{split}
\end{equation} 
where
\begin{eqnarray}\label{eq:T'-1}
\p&=& (k,m_0,\tm) \in Z^{t_{n_o+1}+2}\\
T' & = & \left( \begin{array}{cc|cccccccccccc}
                            2 & -1 &0&.&.&.&.&.&.&.&0\\ 
			   -1 & 2 & -2&.&.&-2&-1&0&.&.&0\\\hline  
	                    0 &-2 & &&&&&&&& \\
                            . & . & &&&&&&&& \\
	                    . & . & &&&&&&&& \\
			    . &-2 & &&&&2B&&&& \\
	                    . &-1 & &&&&&&&& \\
                            . &0 & &&&&&&&& \\
                            . &. & &&&&&&&& \\
                            . &. & &&&&&&&& \\
                            0 &0 & &&&&&&&& \\
			                      
			    \end{array}
         \right)\nonumber\\
         \I_{T'}& = & 2I_{t_{n_o+1}+2}-T' \nonumber
\end{eqnarray}
\noindent
In Ramond-sector we get,
\begin{equation}\label{eq:fermip'>2pN1R-1}
\begin{split}
\tilde{\chi}^{(p',3p'-2p)}_{s, 3b-2r}(q) & =q^{-\frac{3}{8}(b-s)^2+\mathcal{N}_{r,s}-C_{b,s}}
 \sum_{\p\equiv \hat{\overline{Q}}_{\uu}}q^{\frac{1}{4}\p^tT'\p}\\
 &\times q^{-\hat{A}^t\p-\frac{1}{2}\p^t\hat{B}\hat{\uu}_+
-\frac{1}{2}\p^t\hat{\e}_{\nu_0+1}(\hat{\uu}_+^t.\hat{V})-\frac{1}{2}\p^t\hat{\uu}_{-}}
 \times \frac{1}{(q)_{p_2}} \left[ \begin{array}{c}\frac{p_2+1}{2} \\ p_1 \end{array} \right]_q \\
 & \times \prod_{j=3}^{t_{n_0+1}+2}\left[ \begin{array}{c} \frac{1}{2}(\I_{T'}\p+\frac{\hat{B}}{2}
\hat{\uu}_{+}+\frac{1}{2}\hat{\e}_{\nu_0+1}(\hat{\uu}_+^t.\hat{V})+\frac{1}{2}\hat{\uu}_{-})_j 
\\ p_j \end{array} \right]_q^{'}
\end{split}
\end{equation} 
\noindent
where $T'$ is as in (\ref{eq:T'-1}) and  
\begin{equation*}
\p\equiv \hat{\overline{Q}}_{\uu}:=\begin{cases} & 0 \le k\le \frac{m_0}{2}\\
                                            & m_0 \ge 0, \text{$m_0$ is even}\\
                                           & \tm \equiv \overline{Q}_{\uu}
                              \end{cases}  
\end{equation*}     

\section{$N=2$ Character formulas} \label{sec:N2}
\subsection{$N=2$ superconformal algebra and Spectral flow}
The $N=2$ superconformal algebra $\mathcal{A}$ is the infinite dimensional Lie 
super algebra~\cite{EG:1996} with basis $L_n,T_n,G_r^{\pm},C$ and (anti)-commutation 
relation given by
\begin{equation*}
\begin{split}
\left[L_m,L_n\right] &= (m-n)L_{m+n}+\frac{C}{12}(m^3-m)\delta_{m+n,0}\\
\left[L_m,G_r^{\pm}\right] &= (\frac{1}{2}m-r)G_{m+r}^{\pm}\\
\left[L_m,T_n\right]&= -nT_{m+n}\\
\left[T_m,T_n\right] &= \frac{1}{3}cm\delta_{m+n,0}\\
\left[T_m,G_r^{\pm}\right] &= \pm G_{m+r}^{\pm}\\
\{G_r^+,G_s^-\} &= 2L_{r+s}+(r-s)T_{r+s}+\frac{C}{3}(r^2-\frac{1}{4})
\delta_{r+s,0}\\
\left[L_m,C\right]&=\left[T_n,C\right]=\left[G_r^{\pm},C\right]=0\\
\{ G_r^+,G_s^+\} &= \{ G_r^-,G_s^-\}=0
\end{split}
\end{equation*}
where $n,m\in\Z$, but $r,s$ are integers in Ramond-sector and half-integer in NS-sector.
The element $C$ is the central element and its eigenvalue $c$ is parametrized as 
$c=3(1-\frac{2p}{p'})$, where $p,p'$ are relatively prime positive integers.
Let us denote this algebra by $\mathcal{A}(p,p')$.

It was observed in~\cite{K:2003,SS:1987} that there exits a family of outer 
automorphisms $\alpha_{\eta}:\mathcal{A}\rightarrow \mathcal{A}$ which maps the 
$N=2$ superconformal algebras to itself. These are explicitly given by
\begin{equation}\label{eq:auto}
\begin{split}
\alpha_{\eta}(G_r^+)&=\hat{G}_r^+=G_{r-\eta}^+\\
\alpha_{\eta}(G_r^-)&=\hat{G}_r^-=G_{r+\eta}^-\\
\alpha_{\eta}(L_n)&=\hat{L}_n=L_n-\eta T_n +\frac{c}{6}\eta^2\delta_{n,0}\\
\alpha_{\eta}(T_n)&=\hat{T}_n=T_n-\frac{c}{3}\eta \delta_{n,0}
\end{split}
\end{equation}
This family of automorphisms is called \textbf{spectral flow} and $\eta \in \R$ is 
called the \textbf{flow parameter}. 
When $\eta \in \Z$ each sector of the algebra is mapped to itself.
When $\eta \in \Z+\frac{1}{2}$ the Neveu-Schwarz sector is mapped to the Ramond-sector 
and vice-versa. We are going to use the spectral flow $\eta=\pm \frac{1}{2}$ to map 
the NS-sector to the Ramond-sector.

\subsection{Spectral flow and characters}
We denote the Verma module generated from a highest weight state $|h,Q,c\rangle$ with 
$L_0$ eigenvalue $h$, $T_0$ eigenvalue $Q$ and central charge $c$ by $V_{h,Q}$. 
The character $\chi_{V_{h,Q}}$ of a highest weight representation $V_{h,Q}$ is 
defined as
\begin{equation*}
\chi_{V_{h,Q}}(q,z)=\mathrm{Tr}_{V_{h,Q}}(q^{L_0-c/24}z^{T_0}).
\end{equation*}
Following~\cite{K:2003} the character transforms under the spectral flow in the 
following way 
\begin{equation}\label{eq:specchar}
\mathrm{Tr}_{V_{h,Q}}(q^{\hat{L}_0-c/24}z^{\hat{T}_0})
=\mathrm{Tr}_{V_{h^{\eta},Q^{\eta}}}(q^{L_0-c/24}z^{T_0}),
\end{equation}
where $h^{\eta}$ and $Q^{\eta}$ are the eigenvalues of $\hat{L}_0$ and $\hat{T}_0$, 
respectively, as defined in \eqref{eq:auto}. This means the new character 
$\chi_{V_{h^{\eta},Q^{\eta}}}(q,z)$ which is the trace of the transformed operators
over the original representation equals the character of the representation defined 
by the eigenvalues $h^{\eta}$ and $Q^{\eta}$ of $\hat{L}_0$ and $\hat{T}_0$, 
respectively. So the new character is the character of the representation 
$V_{h^{\eta},Q^{\eta}}$.

For $\eta=\frac{1}{2}$ the spectral flow $\alpha_{\frac{1}{2}}$ takes a
NS-sector character to an Ramond-sector character. Let $\chi_{V_{h,Q}}^{NS}(q,z)$ be a 
NS-sector character corresponding to the representation  $V_{h,Q}$. Then 
by~\eqref{eq:specchar} and~\eqref{eq:auto} the new Ramond-sector character
$\chi_{V_{h^{\eta},Q^{\eta}}}^{R}(q,z)$ is derived using
\begin{equation}\label{eq:specns}
\begin{split}
\chi_{V_{h^{\eta},Q^{\eta}}}^{R}(q,z)
&=\mathrm{Tr}_{V_{h,Q}}(q^{\hat{L}_0-c/24}z^{\hat{T}_0})
=\mathrm{Tr}_{V_{h,Q}}(q^{L_0-\frac{1}{2}T_0+\frac{c}{24}-\frac{c}{24}}
z^{T_0-\frac{c}{6}})\\
&=q^{\frac{c}{24}}z^{-\frac{c}{6}}\mathrm{Tr}_{V_{h,Q}}(q^{L_0-\frac{c}{24}}{
(zq^{-\frac{1}{2}})}^{T_0})
=q^{\frac{c}{24}}z^{-\frac{c}{6}}\chi_{V_{h,Q}}^{NS}(q,zq^{-\frac{1}{2}}).
\end{split}
\end{equation}

\subsection{Ramond-sector character from NS-sector character} 
{}To simplify notation we are going to use a slightly different notation for characters. 
Since we are only dealing with the vacuum character in the NS-sector for which $h=0,Q=0$,
we write $\hat{\chi}^{NS}_{p,p'}(q,z)$. The Ramond-sector character is denoted by
$\hat{\chi}^{R}_{p,p'}(q,z)$ with the corresponding $(h,Q)$ specified separately.

Following~\cite{Dobrev:1986,D:1998,EG:1996,K:2003,K:2004} the vacuum character for the $N=2$ 
superconformal algebra with central element $c=3(1-\frac{2p}{p'})$ in the 
NS-sector is given by
\begin{multline}\label{eq:vacuum1}
\hat{\chi}_{p,p'}^{NS}(q,z)=q^{-c/24}\prod_{n=1}^{\infty}\frac{(1+zq^{n-
\frac{1}{2}})(1+z^{-1}q^{n-\frac{1}{2}})}{{(1-q^n)}^2}\\
 \times \Big(1-\sum_{n=0}^\infty q^{p'(n+1)(p(n+1)-1)}
 +\frac{zq^{pn(p'n+1)+p'n+\frac{1}{2}}}{1+zq^{pn+\frac{1}{2}}}+
 \frac{z^{-1}q^{pn(p'n+1)+p'n+\frac{1}{2}}}{1+z^{-1}q^{pn+\frac{1}{2}}}\\
 +\sum_{n=1}^\infty q^{p'n(pn+1)}+\frac{zq^{pn(p'n+1)-p'n-
\frac{1}{2}}}{1+zq^{pn-\frac{1}{2}}}+\frac{z^{-1}q^{pn(p'n+1)-p'n-
\frac{1}{2}}}{1+z^{-1}q^{pn-\frac{1}{2}}}\Big).
\end{multline}
This formula can be verified using the embedding diagram for the vacuum character as 
described in~\cite{EG:1996,K:2003} and can be rewritten as (as will be useful later)
\begin{multline}\label{eq:vacuum2}
\hat{\chi}^{NS}_{p,p'}(q,z)=q^{-c/24}\prod_{n=1}^{\infty} 
\frac{(1+zq^{n-\frac{1}{2}})(1+z^{-1}q^{n-\frac{1}{2}})}{(1-q^n)^2}\\
\times \sum_{j=-\infty}^{\infty}q^{pj(p'j+1)}\frac{1-q^{2p'j+1}}
{(1+zq^{p'j+\frac{1}{2}})(1+z^{-1}q^{p'j+\frac{1}{2}})}.
\end{multline}
The unitary case $p=1$ of these character formulas was given 
in~\cite{Dobrev:1987a,Ki:1987,M:1987,RY:1987}. 
In particular if we put $z=1$ in \eqref{eq:vacuum2} we obtain the following formula 
derived in~\cite{EG:1996}
\begin{equation}\label{eq:vacuum3}
\hat{\chi}^{NS}_{p,p'}(q)=q^{-c/24}\prod_{n=1}^{\infty} 
\frac{(1+q^{n-\frac{1}{2}})^2}{(1-q^n)^2}
\sum_{j=-\infty}^{\infty}q^{pj(p'j+1)}\frac{1-q^{p'j+\frac{1}{2}}}
{1+q^{p'j+\frac{1}{2}}}.
\end{equation}

Let us apply \eqref{eq:specns} to the NS-sector vacuum character \eqref{eq:vacuum2}
to get a Ramond-sector character. From \eqref{eq:auto} it follows that
\begin{equation*}
\begin{split}
\hat{L}_0&=L_0-\frac{1}{2}T_0+\frac{c}{24}\\
\hat{T}_0&=T_0-\frac{c}{6}.
\end{split}
\end{equation*}

For the vacuum character in the NS-sector $(h,Q)=(0,0)$, so the new eigenvalues 
are $(h^{\eta},Q^{\eta})= (\frac{c}{24},-\frac{c}{6})$ in the Ramond-sector. 
Hence the new character in the Ramond-sector corresponds to $(h^{\eta},Q^{\eta})$ and 
by~\eqref{eq:specns}
\begin{multline}\label{eq:Ramond1}
\hat{\chi}^{R}_{p,p'}(q,z)=q^{\frac{c}{24}}z^{-\frac{c}{6}}
\hat{\chi}^{NS}_{p,p'}(q,zq^{-\frac{1}{2}})\\
=z^{-\frac{c}{6}}\frac{(-z)_{\infty}(-z^{-1}q)_{\infty}}{(q)^2_{\infty}}
\sum_{j=-\infty}^{\infty}q^{pj(p'j+1)}
\frac{1-q^{2p'j+1}}{(1+zq^{p'j})(1+z^{-1}q^{p'j+1})}.
\end{multline}

\subsection{Bailey flows from the minimal model $M(p,p')$ to $N=2$ superconformal}
We will consider two set of special values for $r$ and $b$ to find Bailey flows from
the minimal model $M(p,p')$ to $N=2$ superconformal models.
 
First we use $r=0$ and $b=1$ in \eqref{eq:ra1-1} and we obtain
\begin{multline}\label{eq:ra6-1}
\sum_{n=0}^{\infty}(\rho_1)_n(\rho_2)_n(aq/\rho_1\rho_2)^n
\frac{q^{-\mathcal{N}_{0,s}}}{\left(aq\right)_{2n}}
F_{0,s}^{(p,p')}(2n+1-s+2x,1;q)\\
=\frac{(aq/\rho_1)_{\infty}(aq/\rho_2)_{\infty}}{(aq)_{\infty}(aq/\rho_1 \rho_2)_{\infty}}
\sum_{j=-\infty}^{\infty} \Bigg(\frac{(\rho_1)_{jp'-x}(\rho_2)_{jp'-x}}{(aq/\rho_1)_
{jp'-x}(aq/\rho_2)_{jp'-x}}(aq/\rho_1 \rho_2)^{jp'-x}\\
-\frac{(\rho_1)_{jp'-1-x}(\rho_2)_{jp'-1-x}}{(aq/\rho_1)_{jp'-1-x}(aq/\rho_2)_{jp'-1-x}}
(aq/\rho_1 \rho_2)^{jp'-1-x} \Bigg)q^{jp(jp'-s)}.
\end{multline}
 
Then we are going to assume $r(b)=b=1$ in \eqref{eq:ra1-2}. This gives us

\begin{equation}\label{eq:ra6-2}
\begin{split}
&\sum_{n=0}^{\infty}(\rho_1)_n(\rho_2)_n(aq/\rho_1\rho_2)^n\frac{q^{\mathcal{N}_
{1,s}}}{\left(aq\right)_{2n}}q^{n^2} F_{1,s}^{(p,p')}(2n+1-s+2x,1;q^{-1})\\
&=\frac{(aq/\rho_1)_{\infty}(aq/\rho_2)_{\infty}}{(aq)_{\infty}(aq/\rho_1 \rho_2)_{\infty}}
\sum_{j=-\infty}^{\infty} \Bigg(\frac{(\rho_1)_{jp'-x}(\rho_2)_{jp'-x}}{(aq/\rho_1)_
{jp'-x}(aq/\rho_2)_{jp'-x}}(aq/\rho_1 \rho_2)^{jp'-x}\\
& \quad -\frac{(\rho_1)_{jp'-1-x}(\rho_2)_{jp'-1-x}}{(aq/\rho_1)_
{jp'-1-x}(aq/\rho_2)_{jp'-1-x}}(aq/\rho_1 \rho_2)^{jp'-1-x} \Bigg)q^{j^2p'(p'-p)-js(p'-p)-
x(1+x-s)}
\end{split}
\end{equation}

In \eqref{eq:ra6-1} and \eqref{eq:ra6-2} we consider the specialization 
\begin{equation*}
\rho_1=\text{finite},\quad \rho_2=\text{finite}.
\end{equation*} 
so that $\frac{aq}{\rho_1\rho_2}\longrightarrow 1$.
Taking the limit $\frac{aq}{\rho_1\rho_2}\longrightarrow 1$ in \eqref{eq:ra6-1}, we find
\begin{multline}\label{eq:ra6-3}
\sum_{n=0}^{\infty}(\rho_1)_n(\rho_2)_n\frac{q^{-\mathcal{N}_{0,s}}}
{\left(aq\right)_{2n}}F_{0,s}^{(p,p')}(2n+1-s+2x,1;q)\\
=\frac{(\rho_1)_{\infty}(\rho_2)_{\infty}}{(\rho_1\rho_2)_{\infty}(q)_{\infty}} 
\sum_{j=-\infty}^{\infty}q^{jp(jp'-s)}
 \frac{\rho_1\rho_2q^{2(jp'-x-1)}-1}{(1-\rho_1q^{jp'-x-1})(1-\rho_2q^{jp'-x-1})}.
\end{multline}

Similarly taking the limit $\frac{aq}{\rho_1\rho_2}\longrightarrow 1$ in \eqref{eq:ra6-2}, 
we find
\begin{equation}\label{eq:ra6-4}
\begin{split}
&\sum_{n=0}^{\infty}(\rho_1)_n(\rho_2)_n\frac{q^{\mathcal{N}_{1,s}}}{\left(aq
\right)_{2n}}q^{n^2} F_{1,s}^{(p,p')}(2n+1-s+2x,1;q^{-1})\\
&=q^{x(s-x-1)}\frac{(\rho_1)_{\infty}(\rho_2)_{\infty}}{(\rho_1\rho_2)_{\infty}(q)_
{\infty}} \sum_{j=-\infty}^{\infty}q^{j^2p'(p'-p)-js(p'-p)}
 \frac{\rho_1\rho_2q^{2(jp'-x-1)}-1}{(1-\rho_1q^{jp'-x-1})(1-\rho_2q^{jp'-x-1})}
\end{split}
\end{equation}
Now we will consider appropriate finite values for  $\rho_1$ and $\rho_2$.

Here we use two specializations,
\begin{equation}\label{eq:spec_finite}
\rho_1=-zq^{x+\frac{1}{2}} \quad \text{and}\quad \rho_2=-z^{-1}q^{x+\frac{1}{2}},
\end{equation}
and 
\begin{equation}\label{eq:spec_finite2}
\rho_1=-zq^{x} \quad \text{and}\quad \rho_2=-z^{-1}q^{x+1}.
\end{equation}
The specialization \eqref{eq:spec_finite} gives NS-sector characters and 
the specialization \eqref{eq:spec_finite2} gives Ramond-sector characters for $N=2$
superconformal algebras.
\subsection{Fermionic formula for $p<p'<2p$ and $r=0,b=1$}

{\bf NS-sector characters:}\\
Let us use the specialization \eqref{eq:spec_finite} 
in~\eqref{eq:ra6-3}, which implies $a=q^{2x}$ and $s=1$. Making the variable
change $j\longrightarrow -j$ in \eqref{eq:ra6-3} and setting $x=0$ we obtain
\begin{multline}\label{eq:ra3-3}
\sum_{n=0}^{\infty}(-zq^{\frac{1}{2}})_n(-z^{-1}q^{\frac{1}{2}})_n
\frac{q^{-\mathcal{N}_{0,1}}}{\left(q\right)_{2n}} F_{0,1}^{(p,p')}(2n,1;q)\\
=\frac{(-zq^{\frac{1}{2}})_{\infty}(-z^{-1}q^{\frac{1}{2}})_{\infty}}
{(q)_{\infty}^2} \sum_{j=-\infty}^{\infty}q^{jp(jp'+1)}
 \frac{1-q^{2jp'+1}}{(1+zq^{jp'+\frac{1}{2}})(1+z^{-1}q^{jp'+\frac{1}{2}})}.
\end{multline}
Comparing with \eqref{eq:vacuum2}, we obtain
\begin{equation}\label{eq:ra3-5}
\hat{\chi}^{NS}_{p,p'}(q,z)=q^{-\frac{c}{24}-\mathcal{N}_{0,1}}
\sum_{n=0}^{\infty} \frac{(-zq^{\frac{1}{2}})_n(-z^{-1}q^{\frac{1}{2}})_n}{\left(q\right)_{2n}}
F_{0,1}^{(p,p')}(2n,1;q).
\end{equation}
This gives us a Bailey flow from $M(p,p')$ model to $N=2$ superconformal model in the NS
sector. Now we calculate the fermionic side for $p<p'<2p$.
  
Setting $z=1$ and inserting the fermionic formula \eqref{eq:fermi1-1}, we find
\begin{multline}\label{eq:fermi3-1}
\hat{\chi}^{NS}_{p,p'}(q)=q^{-\frac{c}{24}-\mathcal{N}_{0,1}+k_{1,1}}
\sum_{n=0}^{\infty}\Bigg( \frac{{(-q^{\frac{1}{2}})}^2_n}{(q)_{2n}}
 \sum_{\m\equiv Q_{\uu+\vv}} q^{\frac{1}{4}\m^tB\m-\frac{1}{2}A_{\uu,\vv}\m}\\
 \times \prod_{j=1}^{t_{n_0+1}}\left[ \begin{array}{c} n_j+m_j\\m_j \end{array} 
 \right]_q^{'}\Bigg).
\end{multline}
Let us set $m_0=2n$ and use \eqref{eq:-q1/2} to get
\begin{multline}\label{eq:fermi3-2}
\hat{\chi}^{NS}_{p,p'}(q) = q^{-\frac{c}{24}-\mathcal{N}_{0,1}+k_{1,1}} 
\sum_{\substack{m_0=0\\ \text{$m_0$ even}}}^{\infty}\sum_{k_1=0}^{\frac{m_0}{2}}
 \sum_{k_2=0}^{\frac{m_0}{2}}\sum_{\m\equiv Q_{\uu+\vv}}
 q^{\frac{1}{2}(\frac{m_0}{2}-k_1)^2+\frac{1}{2}(\frac{m_0}{2}-k_2)^2}\\
 \times q^{\frac{1}{4}\m^tB\m-\frac{1}{2}A_{\uu,\vv}\m}
 \frac{1}{(q)_{m_0}} \left[ \begin{array}{c}\frac{m_0}{2} \\ k_1 \end{array} 
 \right]_q
\left[ \begin{array}{c}\frac{m_0}{2} \\ k_2 \end{array} \right]_q 
\prod_{j=1}^{t_{n_0+1}}\left[ 
\begin{array}{c} n_j+m_j\\m_j \end{array} \right]_q^{'}.
\end{multline}

Define $\p=(k_1,k_2,m_0,\m)\in \Z^{t_{n_0+1}+3}$, so that \eqref{eq:fermi3-2}
can be rewritten as
\begin{multline}\label{eq:fermi3-3}
\hat{\chi}^{NS}_{p,p'}(q)=q^{-\frac{c}{24}-\mathcal{N}_{0,1}+k_{1,1}}
\sum_{\substack{\p\in\Z^{t_{n_0+1}+3}\\ p_i\equiv (\hat{Q}_{\uu,\vv})_i, i\ge 3}}
q^{\frac{1}{4}\p^tD\p -\frac{1}{2}\hat{A}_{\uu,\vv}\p}\\
 \times \frac{1}{(q)_{p_3}}  \prod_{\substack{j=1,j\ne3}}^{t_{n_0+1}+3}\left[ 
\begin{array}{c} \frac{1}{2}(\I_D\p+\hat{\uu}+\hat{\vv})_j\\p_j \end{array} \right]_q^{'},
\end{multline}
where $\I_D=2I_{t_{n_0+1}+3}-D$ and
\begin{equation}\label{I-D}
\begin{split}
D&=\left( \begin{array}{ccc|cc}
            2 & 0 & -1 & 0\\ 
	    0 & 2 & -1 & 0\\ 
	   -1 & -1 & 1 & 1\\\hline 
            0 & 0 &  -1 & B\\
		\end{array}
         \right),\\
\hat{A}_{\uu,\vv}&=(0,0,0,A_{\uu,\vv}),\\
\hat{\uu}^t&=(0,0,0,\uu^t),\\
\hat{\vv}^t&=(0,0,0,\vv^t),\\
\hat{Q}^t_{\uu,\vv}&=(0,0,0,Q_{\uu+\vv}^t).
\end{split}
\end{equation}
This gives a new fermionic expression for the NS-sector character.\\

{\bf Ramond-sector characters:}
Let us set $\rho_1=-zq^x,\rho_2=-z^{-1}q^{x+1}$ in \eqref{eq:ra6-4}, which implies 
$a=q^{2x}$ and $s=1$. Setting $x=0$ and changing $j\longrightarrow -j$ we obtain
\begin{multline}\label{eq:R2}
\sum_{n=0}^{\infty} \frac{{(-z)}_n{(-z^{-1}q)}_n}{(q)_{2n}} q^{-\mathcal{N}_{0,1}} 
F_{0,1}^{(p,p')}(2n,1;q) \\
=\frac{{(-z)}_{\infty}{(-z^{-1}q)}_{\infty}}{{(q)^2}_{\infty}}\sum_{j=-\infty}^{\infty}
q^{jp(jp'+1)}\frac{1-q^{2jp'+1}}{(1+zq^{jp'})(1+z^{-1}q^{jp'+1})}.
\end{multline}
Comparing with \eqref{eq:Ramond1} we get 
\begin{equation}\label{eq:R3}
\hat{\chi}^{R}_{p,p'}(q,z)=z^{-\frac{c}{6}}q^{-\mathcal{N}_{0,1}}\sum_{n=0}^{\infty} 
\frac{{(-z)}_n{(-z^{-1}q)}_n}{(q)_{2n}}F_{0,1}^{(p,p')}(2n,1;q). 
\end{equation}
This shows a Bailey flow from $M(p,p')$ to $N=2$ superconformal algebra in the Ramond-sector.

Again using \eqref{eq:fermi1-1} in a similar way to the NS-sector and setting $z=1$ we find
\begin{multline}\label{eq:Rfermi3-1}
\hat{\chi}^{R}_{p,p'}(q)= 2q^{-\mathcal{N}_{0,1}+k_{1,1}}
\sum_{n=0}^{\infty}\Bigg( \frac{{(-q)}_{n-1}{(-q)}_n}{(q)_{2n}}
 \sum_{\m\equiv Q_{\uu+\vv}} q^{\frac{1}{4}\m^tB\m-\frac{1}{2}A_{\uu,\vv}\m}\\
 \times \prod_{j=1}^{t_{n_0+1}}\left[ 
\begin{array}{c} n_j+m_j\\m_j \end{array} \right]_q^{'}\Bigg).
\end{multline}
Using
\begin{equation*}
(x)_n =\sum_{k=0}^n (-x)^{(n-k)}q^{\frac{1}{2}(n-k)(n-k-1)}\left[ \begin{array}{c} n\\k \end{array} 
\right]_{q}
\end{equation*}
and setting $m_0=2n$, equation \eqref{eq:Rfermi3-1} can be rewritten as
\begin{multline}\label{eq:Rfermi3-2}
\hat{\chi}^{R}_{p,p'}(q) = 2q^{-\mathcal{N}_{0,1}+k_{1,1}}
\sum_{\substack{m_0=0\\ \text{$m_0$ even}}}^{\infty}\sum_{k_1=0}
^{\frac{m_0}{2}-1}\sum_{k_2=0}^{\frac{m_0}{2}}\sum_{\m\equiv Q_{\uu+\vv}}
 q^{\frac{1}{4}(m_0^2+2k_1^2+2k_2^2-2m_0k_1-2m_0k_2)}\\
 \times q^{\frac{1}{4}\m^tB\m-\frac{1}{2}A_{\uu,\vv}\m+\frac{1}{2}(k_1-k_2)}
 \frac{1}{(q)_{m_0}} \left[ \begin{array}{c}\frac{m_0}{2}-1 \\ k_1 \end{array} \right]_q
\left[ \begin{array}{c}\frac{m_0}{2} \\ k_2 \end{array} \right]_q 
\prod_{j=1}^{t_{n_0+1}}\left[ 
\begin{array}{c} n_j+m_j\\m_j \end{array} \right]_q^{'}.
\end{multline}
Setting $\p=(k_1,k_2,m_0,\m)\in \Z^{t_{n_0+1}+3}$ this becomes
\begin{multline}\label{eq:Rfermi3-3}
\hat{\chi}^{R}_{p,p'}(q)= 2q^{-\mathcal{N}_{0,1}+k_{1,1}}
\sum_{\substack{\p\in \Z^{t_{n_0+1}+3}\\ p_i\equiv (\hat{Q}_{\uu,\vv})_i, i\ge 3}}
q^{\frac{1}{4}\p^tD\p -\frac{1}{2}\hat{A}_{\uu,\vv}\p}\\
 \times \frac{1}{(q)_{p_3}}
\prod_{\substack{j=1,j\ne 3}}^{t_{n_0+1}+3}\left[ 
\begin{array}{c} \frac{1}{2}(\I_D\p+\hat{\uu}+\hat{\vv})_j\\p_j \end{array} \right]_q^{'},
\end{multline}
with the same notations as in \eqref{I-D} except
\begin{equation*}
\begin{split}
\hat{A}_{\uu,\vv}&=(1,-1,0,A_{\uu,\vv}),\quad \hat{\uu}^t=(-1,0,0,\uu^t),\quad
\hat{\vv}^t=(-1,0,0,\vv^t).
\end{split}
\end{equation*} 
This gives a new fermionic expression of the new Ramond-sector character.
\subsection{Fermionic formula for $p<p'<2p$ and $r=b=1$}

{\bf NS-sector characters:}\\
Now we use the specialization \eqref{eq:spec_finite} in \eqref{eq:ra6-4}. 
This implies $a=q^{2x}, s=1$ and by taking $j\longrightarrow -j$ in (\ref{eq:ra6-4}) we get,
\begin{equation}\label{eq:ra3-3.5}
\begin{split}
&\sum_{n=0}^{\infty}(-zq^{x+\frac{1}{2}})_n(-z^{-1}q^{x+\frac{1}{2}})_n\frac{q^{\mathcal{N}_{1,1}}}
{\left(q^{2x+1}\right)_{2n}}q^{n^2} F_{1,1}^{(p,p')}(2n+2x;q^{-1})\\
&=q^{-x^2}\frac{(-zq^{x+\frac{1}{2}})_{\infty}(-z^{-1}q^{x+\frac{1}{2}})_{\infty}}
{(q^{2x+1})_{\infty}(q)_
{\infty}} \sum_{j=-\infty}^{\infty}q^{j(p'-p)(jp'+1)}
 \frac{1-q^{2jp'+1}}{(1+zq^{jp'+\frac{1}{2}})(1+q^{jp'+\frac{1}{2}})}.
\end{split}
\end{equation}
\noindent
Now let $x=0$ in (\ref{eq:ra3-3.5}) to get 
\begin{equation}\label{eq:ra3-4}
\begin{split}
&\sum_{n=0}^{\infty}(-zq^{\frac{1}{2}})_n(-z^{-1}q^{\frac{1}{2}})_n\frac{q^{\mathcal{N}_{1,1}}}
{\left(q\right)_{2n}}q^{n^2} F_{1,1}^{(p,p')}(2n;q^{-1})\\
&=\frac{(-zq^{\frac{1}{2}})_{\infty}(-z^{-1}q^{\frac{1}{2}})_{\infty}}{(q)^2_{\infty}} 
\sum_{j=-\infty}^{\infty}q^{j(p'-p)(jp'+1)}
 \frac{1-q^{2jp'+1}}{(1+zq^{jp'+\frac{1}{2}})(1+q^{jp'+\frac{1}{2}})}.
\end{split}
\end{equation}
Comparing with (\ref{eq:vacuum2}) we get,
\begin{equation}\label{eq:ra3-4.5}
\hat{\chi}^{NS}_{p'-p,p'}(q,z)=q^{-c(p'-p,p')+\mathcal{N}_{1,1}}\sum_{n=0}^{\infty}
\frac{(-zq^{\frac{1}{2}})_n(-z^{-1}q^{\frac{1}{2}})_n}
{\left(q\right)_{2n}}q^{n^2} F_{1,1}^{(p,p')}(2n;q^{-1}).
\end{equation}
Now putting $z=1$ and comparing with (\ref{eq:vacuum3}) we get 
\begin{equation}\label{eq:ra3-5.5}
\hat{\chi}^{NS}_{p'-p,p'}(q)=q^{-c(p'-p,p')+\mathcal{N}_{1,1}}
\sum_{n=0}^{\infty}\frac{{(-q^{\frac{1}{2}})}^2_n}{(q)_{2n}}q^{n^2}F_{1,1}^{p,p'}(2n;q^{-1}).
\end{equation}
When $p<p'<2p$ recall that $\mathcal{N}_{1,1}=0$. Hence the formula
looks like
\begin{equation}\label{eq:r=b=1_1}
\hat{\chi}^{NS}_{p'-p,p'}(q)=q^{-c(p'-p,p')}
\sum_{n=0}^{\infty}\frac{{(-q^{\frac{1}{2}})}^2_n}{(q)_{2n}}q^{n^2}F_{1,1}^{p,p'}(2n;q^{-1}).
\end{equation}
Let us calculate the fermionic side using the fermionic formula 
\eqref{eq:ferminverse_r=b=1} from section~\ref{sec:fermionic}.
We get
\begin{equation}\label{eq:r=b=1_2}
\begin{split}
\hat{\chi}^{NS}_{p'-p,p'}(q)=&q^{-c(p'-p,p')}
\sum_{n=0}^{\infty}\frac{{(-q^{\frac{1}{2}})}^2_n}{(q)_{2n}}q^{n^2}\times \sum_{\hat{\m}\equiv Q(\uu)}
\Bigg(q^{\frac{1}{4}\m^t\overline{C}\m-nm_1}\\
\times &q^{\frac{1}{2}(\uu'_{\flat}+\uu'_{\sharp})\m-\frac{1}{2}\overline{\uu}\m-\frac{1}{4}\gamma(\uu')}
\prod_{j=1}^{t-1} \left[ \begin{array}{c}m_j+n_j\\ m_j
\end{array}\right]_q\Bigg)\\
+&q^{-c(p'-p,p')+\mathcal{N}_{1,1}}
\sum_{n=0}^{\infty}\frac{{(-q^{\frac{1}{2}})}^2_n}{(q)_{2n}}q^{n^2}F_{1,1}^{\hat{p},\hat{p'}}(2n;q^{-1}).
\end{split}
\end{equation} 

Let $m_0=L=2n$ and we can rewrite \eqref{eq:r=b=1_2} as
\begin{equation}\label{eq:r=b=1_3}
\begin{split}
\hat{\chi}^{NS}_{p'-p,p'}(q) 
&= q^{-c(p'-p,p')} 
\sum_{\substack{m_0=0\\ \text{$m_0$ even}}}^{\infty}\sum_{k_1=0}
^{\frac{m_0}{2}}\sum_{k_2=0}^{\frac{m_0}{2}}\sum_{\m\equiv Q(\uu)}\\ 
& \times \Bigg(q^{\frac{1}{2}(\frac{m_0}{2}-k_1)^2+\frac{1}{2}(\frac{m_0}{2}-k_2)^2+\frac{m_0^2}{4}
-\frac{1}{2}m_0m_1+\frac{1}{2}(\uu'_{\flat}+\uu'_{\sharp})\m-\frac{1}{2}\overline{\uu}\m-\frac{1}{4}\gamma(\uu')}\\
&\times \frac{1}{(q)_{m_0}} \left[ \begin{array}{c}\frac{m_0}{2} \\ k_1 \end{array} \right]_q
\left[ \begin{array}{c}\frac{m_0}{2} \\ k_2 \end{array} \right]_q \prod_{j=1}^{t-1} 
\left[ \begin{array}{c}m_j+n_j\\ m_j
\end{array}\right]_q\Bigg) \\
&+q^{-c(p'-p,p')}
\sum_{n=0}^{\infty}\frac{{(-q^{\frac{1}{2}})}^2_n}{(q)_{2n}}q^{n^2}F_{1,1}^{\hat{p},\hat{p'}}(2n;q^{-1}).
\end{split}
\end{equation}

Define $\p=(k_1,k_2,m_0,\m)$, we can rewrite \eqref{eq:r=b=1_3} as 
\begin{equation}\label{eq:r=b=1_4}
\begin{split}
\hat{\chi}^{NS}_{p'-p,p'}(q)=&q^{-c(p'-p,p')}
\sum_{\p\equiv \hat{Q}(\uu)}\Bigg(q^{\frac{1}{4}\p^t\tilde{D}\p
+\frac{1}{2}(\hat{\uu}'_{\flat}+\hat{\uu}'_{\sharp})\p-\frac{1}{4}\hat{\gamma}(\uu')}\\
& \times \frac{1}{(q)_{p_3}} 
\prod_{j=1,j\ne 3}^{t+2} \left[ \begin{array}{c}\frac{1}{2}(\I_{\tilde{D}}\p+\hat{\uu})_j\\p_j
\end{array}\right]_q \Bigg)\\
+&q^{-c(p'-p,p')+\mathcal{N}_{1,1}}
\sum_{n=0}^{\infty}\frac{{(-q^{\frac{1}{2}})}^2_n}{(q)_{2n}}q^{n^2}F_{1,1}^{\hat{p},\hat{p'}}(2n;q^{-1})
\end{split}
\end{equation} 
where\begin{equation}\label{I-Dtilde}
\begin{split}
\tilde{D}& =  \left( \begin{array}{ccc|cc}
                            2 & 0 & -1 & 0\\ 
			    0 & 2 & -1 & 0\\ 
			    -1 & -1 & 2 & -1\\\hline 
			    0 & 0 & -1 & \overline{C}\\
			    \end{array}
         \right)\\
\I_{\tilde{D}} & = 2I_{t+2}-\tilde{D}\\
\hat{\uu} &=(0,0,0,\overline{\uu}) \\
\hat{\uu}'_\flat+\hat{\uu}'_\sharp &=(0,0,0,\uu_\flat'+\uu_\sharp'-\overline{\uu}) \\
\p\equiv \hat{Q}(\uu):&=\begin{cases} & 0 \le k_1\le \frac{m_0}{2}\\
                                           & 0 \le k_2\le \frac{m_0}{2}\\
                                           & m_0 \ge 0, \text{$m_0$ is even}\\
                                           & \m \equiv Q(\overline{\uu}),
                              \end{cases} \nonumber
\end{split}
\end{equation}
Note that 
\begin{equation}
\hat{\chi}^{NS}_{\hat{p'}-\hat{p},\hat{p'}}(q)=
q^{-c(\hat{p'}-\hat{p},\hat{p'})}
\sum_{n=0}^{\infty}\frac{{(-q^{\frac{1}{2}})}^2_n}{(q)_{2n}}q^{n^2}F_{1,1}^{\hat{p},\hat{p'}}(2n;q^{-1}).
\end{equation}
Hence \eqref{eq:ra3-5} can be written as
\begin{equation}
\begin{split}
\hat{\chi}^{NS}_{p'-p,p'}(q)=&q^{-c(p'-p,p')}
\sum_{\p\equiv \hat{Q}(\uu)}\Bigg(q^{\frac{1}{4}\p^t\tilde{D}\p
+\frac{1}{2}(\hat{\uu}'_\flat+\hat{\uu}'_\sharp)\p-\frac{1}{4}\hat{\gamma}(\uu')}\\
& \times \frac{1}{(q)_{p_3}} 
\prod_{j=1,j\ne 3}^{t+2} \left[ \begin{array}{c}\frac{1}{2}(\I_{\tilde{D}}\p+\hat{\uu})_j\\p_j
\end{array}\right]_q\Bigg)\\
+&q^{-c(p'-p,p')+c(\hat{p'}-\hat{p},\hat{p'})}\hat{\chi}^{NS}_{\hat{p'}-\hat{p},\hat{p'}}(q).
\end{split}
\end{equation}

{\bf Ramond-sector characters:}\\
Let us use the specialization \eqref{eq:spec_finite2} in \eqref{eq:ra6-4} with $x=0$ to get
\begin{equation}\label{eq:r=b=1_R1}
\begin{split}
\sum_{n=0}^{\infty} & \frac{{(-z)}_n{(-z^{-1}q)}_n}{(q)_{2n}} q^{\mathcal{N}_{1,1}} 
q^{n^2}F_{1,1}^{(p,p')}(2n;q^{-1}) \\
& =\frac{{(-z)}_{\infty}{(-z^{-1}q)}_{\infty}}{{(q)^2}_{\infty}}\sum_{j=-\infty}^{\infty}
q^{j(p'-p)(jp'+1)}\frac{1-q^{2jp'+1}}{(1+zq^{jp'})(1+z^{-1}q^{jp'+1})}.
\end{split}
\end{equation}
Comparing with \eqref{eq:Ramond1} we get 
\begin{equation}\label{eq:r=b=1_R2}
\hat{\chi}^{R}_{p'-p,p'}(q,z)=z^{-\frac{c(p,p')}{6}}q^{\mathcal{N}_{1,1}}\sum_{n=0}^{\infty} 
\frac{{(-z)}_n
{(-z^{-1}q)}_n}{(q)_{2n}} q^{n^2}F_{1,1}^{(p,p')}(2n;q^{-1}).
\end{equation}
Let us put $z=1$  and calculate the fermionic side using 
\eqref{eq:ferminverse_r=b=1}. The calculations are similar to the NS-sector case.
Hence we just state the result here.
We get
\begin{equation}\label{eq:r=b=1_R3}
\begin{split}
\hat{\chi}^{R}_{p'-p,p'}(q)= & 2
\sum_{\p\equiv \hat{Q}(\uu)}\Bigg( q^{\frac{1}{4}\p^t\tilde{D}\p +
\frac{1}{2}(\hat{\uu}'_\flat+\hat{\uu}'_\sharp)\p-\frac{1}{4}\hat{\gamma}(\uu')}\\
 & \times \frac{1}{(q)_{p_3}} \prod_{j=1,j\ne 3}^{t+2}\left[ 
\begin{array}{c} \frac{1}{2}(\I_{\tilde{D}}\p+\hat{\uu})_j\\p_j \end{array} \right]_q^{'}\Bigg)
+\hat{\chi}^{R}_{\hat{p}'-\hat{p},\hat{p}'}(q),
\end{split}
\end{equation}
\noindent
with the same notations as in (\ref{I-Dtilde}) except
\begin{equation}
\p\equiv \hat{Q}(\uu):=\begin{cases} & 0 \le k_1\le \frac{m_0}{2}-1\\
                                           & 0 \le k_2\le \frac{m_0}{2}\\
                                           & m_0 \ge 0, \text{$m_0$ is even}\\
                                           & \m \equiv Q(\overline{\uu}),
                              \end{cases}
\end{equation} 
$\hat{\uu}'_\flat+\hat{\uu}'_\sharp=(1,-1,0,\uu_\flat'+\uu_\sharp')$ and 
$\hat{\uu}=(-2, 0,0,\overline{\uu})$.

\subsection{Fermionic formula for $p'>2p$} 
In this case we only give the result for $r=b=1$ case. The fermionic formula for
the case $r=0,b=1$ and $p'>2p$ can be found similarly using the formula given in
\cite{W:2002}.   

{\bf NS-sector character:}\\

All the calculations are done in a similar 
way as in the previous section,  but using the fermionic formula for $p'>2p$. 
Hence we just state the result here. 
  
The explicit fermionic formula in this case looks like
\begin{equation}\label{eq:fermip'>2pN2NS}
\begin{split}
\hat{\chi}^{NS}_{p'-p,p'}(q)= & q^{-c(p'-p,p')+\mathcal{N}_{1,1}-C_{1,1}}
\sum_{\p\equiv \hat{\overline{Q}}}\Bigg( q^{\frac{1}{4}\p^tD'\p -\hat{A'}\p}\\
&\times q^{-\frac{1}{2}\p^t\hat{B'}\hat{\uu_+}-\frac{1}{2}\p^t\hat{\e}_{\nu_0+1}
(\hat{\uu}_+^t.\hat{V})-\frac{1}{2}\p^t\hat{\uu_{-}}}
\times \frac{1}{(q)_{p_3}} \left[ \begin{array}{c}\frac{p_3}{2} \\ p_1 \end{array} \right]_q
\left[ \begin{array}{c}\frac{p_3}{2} \\ p_2 \end{array} \right]_q \\
&\times \prod_{j=4}^{t_{n_0+1}+3}\left[ 
\begin{array}{c} \frac{1}{2}(\I_D'\p+\hat{B'}\hat{\uu}_{+}+\hat{\e'}_{\nu_0+1}(\hat{\uu'}_+^t.
\hat{V'})+\hat{\uu'}_{-})_j\\p_j \end{array} \right]_q^{'}\Bigg)
\end{split}
\end{equation}
\noindent
where
\begin{eqnarray}\label{I-D'}
D' & = & \left( \begin{array}{ccc|cccccccccccc}
                            2 & 0 &-1 &0 &.&.&.&.&.&.&.&0\\ 
			                      0 & 2 &-1 &0 &.&.&.&.&.&.&.&0\\
			                     -1 & -1 & 2 &-2&.&.&-2&-1&0&.&.&0\\\hline 
			                      0 & 0 &-2 & &&&&&&&& \\
			                      . & . &.  & &&&&&&&& \\
			                      . & . &.  & &&&&&&&& \\
			                      . & . &-2 & &&&&2B&&&& \\
			                      . & . &-1 & &&&&&&&& \\
			                      . & . & 0 & &&&&&&&& \\
			                      . & . & . & &&&&&&&& \\
			                      . & . & . & &&&&&&&& \\
			                      0 & 0 & 0 & &&&&&&&& \\
			                      \end{array}\right)\nonumber\\
\I_{D'} & = & 2I_{t_{n_0+1}+3}-D' \nonumber\\
\hat{A'}&=& (0,0,0,A)\\
\hat{X}^t&=& (0,0,0,X^t)\nonumber\\
\p\equiv \hat{\overline{Q'}}_{\uu}&:=&\begin{cases} & 0 \le k_1\le \frac{m_0}{2}\\
                                           & 0 \le k_2\le \frac{m_0}{2}\\
                                           & m_0 \ge 0, \text{$m_0$ is even}\\
                                           & \tm \equiv \overline{Q}_{\uu}
                              \end{cases} \nonumber
\end{eqnarray}
\noindent
Note that $D'$ is a $(t_{n_0+1}+3)\times(t_{n_0+1}+3)$  matrix and  number of -2 in the 
third row and column is $\nu_0$.  

{\bf Ramond-sector character:} Using the specialization \eqref{eq:spec_finite2} we get the
Ramond-sector characters. All the calculations are similar to the
previous case. Hence we state the formula only. We get
\begin{equation}\label{eq:fermip'>2pN2R}
\begin{split}
\hat{\chi}^{R}_{p'-p,p'}(q)= & 2q^{\mathcal{N}_{1,1}-C_{1,1}}
\sum_{\p\equiv \hat{\overline{Q}}}\Bigg( q^{\frac{1}{4}\p^tD'\p -\hat{A'}^t\p}\\
&\times q^{-\frac{1}{2}\p^t\hat{B'}\hat{\uu_+}-\frac{1}{2}\p^t\hat{\e}_{\nu_0+1}
(\hat{\uu}_+^t.\hat{V})-\frac{1}{2}\p^t\hat{\uu_{-}}}
\times \frac{1}{(q)_{p_3}} \left[ \begin{array}{c}\frac{p_3}{2}-1 \\ p_1 \end{array} \right]_q
\left[ \begin{array}{c}\frac{p_3}{2} \\ p_2 \end{array} \right]_q \\
&\times \prod_{j=4}^{t_{n_0+1}+3}\left[ 
\begin{array}{c} \frac{1}{2}(\I_{D'}\p+\hat{B'}\hat{\uu}_{+}+\hat{\e'}_{\nu_0+1}
(\hat{\uu'}_+^t.\hat{V'})+\hat{\uu'}_{-})_j\\p_j \end{array} \right]_q^{'}\Bigg)
\end{split}
\end{equation}
with the same notations as in (\ref{I-D'}) except
\begin{equation}
\p\equiv \hat{\overline{Q'}}_{\uu}:=\begin{cases} & 0 \le k_1\le \frac{m_0}{2}-1\\
                                           & 0 \le k_2\le \frac{m_0}{2}\\
                                           & m_0 \ge 0, \text{$m_0$ is even}\\
                                           & \m \equiv \overline{Q}_{\uu}
                              \end{cases} \nonumber
\end{equation} 
\noindent 
(\ref{eq:Rfermi3-3})and (\ref{eq:fermip'>2pN2R}) gives us new fermionic expressions of the new 
Ramond-sector character.

\section{Conclusion} \label{sec:conclusion}
In this dissertation we only considered the vacuum character for the $N=2$ superconformal algebra
with central charge $c=3(1-\frac{2p}{p'})$ with $p<p'$ in the NS-sector and the 
Ramond-sector character derived from the vacuum character. We believe that 
similar Bailey flows exist for the general $N=2$ superconformal characters, but explicit
formulas are not yet available in the literature.

    %
    %
     \newchap{Summary of implementation and user manual}
      \label{chap:program}    
      
In this chapter we describe the programs that were used to
verify the conjectures for our results on unrestricted Kostka 
polynomials presented in Chapter~\ref{chap:kostka}. The bijection $\Phi$ in 
Chapter~\ref{chap:kostka} has been implemented as a program written in C++. Several 
different versions of this program have been used to carry out calculations 
regarding the unrestricted rigged configurations. We used six different
programs to verify conjectures regarding the lower bound 
conditions, the convexity
property of the unrestricted rigged configurations, and the fact that the bijection $\Phi$
preserves the statistics. We describe three of these programs in this chapter.
 The programs presented here can be used by anyone studying unrestricted Kostka polynomials. 
The code for these programs are provided in the appendix and can be downloaded
at {\tt http://math.ucdavis.edu/\~{}deka}.

The programs have also been incorporated into MuPAD-Combinat as a dynamic
module by Francois Descouens~\cite{MuPAD:2005}.
For example, the command
\begin{multline*}
\text{{\tt riggedConfigurations::RcPathsEnergy::}}\\
\text{{\tt fromOnePath([[[3]],[[2],[1]],[[4,5,6],[1,2,3]]])}}
\end{multline*}
calculates $\Phi(b)$ with $b$ as in Example~\ref{ex:b}.
 
The programs have been compiled and tested using the gnu c++ compiler, (g++), 
version 3.2.3.

\section{ Program: $allpaths\_bijection.c$}\label{prog:allpaths}

The program named $allpaths\_bijection.c$ performs the bijection $\Phi$ from
the set of unrestricted paths $\Path(B,\la)$ to the set of  
unrestricted rigged configurations $\RC(L,\la)$ for a fixed value of
$\la$ and $L$. Let us recall that $\la$ is the weight vector for the unrestricted 
paths in  $\Path(B,\la)$ and
$L$ is the multiplicity matrix for the shape of the tensor product
$B=B^{r_k,s_k}\otimes\cdots \otimes B^{r_1,s_1}$. Let us denote the shape of 
$B$ by a sequence of rectangular partitions $\mu$. 
This program first calculates the set $\Path(B,\la)$. Then, for each $b \in \Path(B,\la)$
it calculates $\Phi(b)$, thus calculates the set $\RC(L,\la)$. The program sorts the elements
of $\RC(L,\la)$ so that all the rigged configurations with the same shape appear
together. It also calculates the statistics for each pair $(b, \Phi(b))$
and finally prints out the unrestricted Kostka polynomial $X(\la,B)$.\\

{\bf Input:}
Here we explain how to input data for the program.
The input file for this program is called  $input\_allpaths$.
The input data for this program are:\\

$n$ = The rank of the Lie algebra of type $A_n$ and we input it in the 1st line of
  the input file. \\

$\la$ = The fixed weight of the unrestricted paths. $\la$ is an  $n+1$ 
        tuple of non-negative numbers. We enter this in the 2nd line of 
        the input file with exactly $n+1$ parts including $0$ if necessary.\\

$\mu$ = The fixed shape of the paths. $\mu$ is a sequence of rectangular 
    partitions. We enter a rectangular partition column-wise.
    For example, 2 2 2 0 represents a rectangular partition with 2 boxes in
    the 1st column, 2 boxes in the 2nd column and 2 boxes in the 3rd  
    column. The 0 at the end indicates the end of that partition. We input 
    each component of $\mu$ in a new line.\\
  
We illustrate how to enter $n$, $\la$ and $\mu$ to the program using a small example
given below.
\begin{equation*}
n=5, \quad \la=(1, 2, 2, 2,1 0),\quad \mu=((1,1),(3,3)).  
\end{equation*}
The input file is:
\begin{verbatim}
5
1 2 2 2 1 0
2 0
2 2 2 0
\end{verbatim}
\noindent
{\bf WARNING:} Do not leave any extra blank space at the end of a line. Do not
forget to include 0 if necessary to make $\la$ an $n+1$ tuple. Do not
forget to put $0$ at the end of a part of $\mu$. The program
will read the input data incorrectly if you forget any of these and 
you will get a wrong answer.

\begin{remark} The maximum size of the rank $n$ of the algebra $A_n$ is limited by
 `RIGSIZE' which is defined to be 20. For a larger $n$, the `RIGSIZE' needs to be 
increased accordingly in the beginning of the program. For a very large $n$, 
the program might take longer to compile and run. 
\end{remark}

{\bf Output:} Let us consider the input data:
$n=3, \la=(0, 1, 1, 1)$ and $\mu=((1),(1,1))$. Input file is
\begin{verbatim}
3
0 1 1 1
1 0
2 0
\end{verbatim}

The output of the program for this example is shown below. 
\begin{verbatim}
n = 3 
Lambda is:  0 1 1 1 
mu is :
..........
 1 
..........
 2 
---------------------------------------
There are 3 unrestricted paths.
---------------------------------------

Path (1): 
---------------------------------------
 2 
---------------------------------------
 3 
 4 

Corresponding rigged configuration is:
---------------------------------------
(1)
 ___ ___
|   | -1| -1
 --- ---
---------------------------------------
(2)
 ___
| -1| -1
 ---
| -1| -1
 ---
---------------------------------------
(3)
 ___
|  0| 0
 ---
Statistic is = 0   
************************************************

Path (2): 
---------------------------------------
 3 
---------------------------------------
 2 
 4 

Corresponding rigged configuration is:
---------------------------------------
(1)
 ___ ___
|   | -1| -1
 --- ---
---------------------------------------
(2)
 ___ ___
|   |  0| 0
 --- ---
---------------------------------------
(3)
 ___
| -1| -1
 ---
Statistic is = 0   
************************************************

Path (3): 
---------------------------------------
 4 
---------------------------------------
 2 
 3 

Corresponding rigged configuration is:
---------------------------------------
(1)
 ___
| -1| -1
 ---
| -1| -1
 ---
---------------------------------------
(2)
 ___
|  0| 0
 ---
|  0| 0
 ---
---------------------------------------
(3)
 ___
|  0| 0
 ---
Statistic is = 1   
************************************************

Unrestricted Kostka polynomial is:  2q^0  + 1q^1 
\end{verbatim}
The output  means there are 3 unrestricted paths. Each of the path
and the corresponding rigged configuration are printed together along with the
statistic. Finally the unrestricted Kostka polynomial,
$X(B,\la)=2+q$. 

\section{Program: $one\_path\_bij.c$}\label{prog:onepath}

The program called $one\_path\_bij.c$ calculates the image of
an unrestricted path $b\in \Path(B,\la)$ under the map $\Phi$.
The output is an unrestricted rigged configuration in $\RC(L,\la)$.
The program also calculates the corresponding statistic. 

{\bf Input:}
 The input file for this program is called  \textit{inputpath}.
We explain how to enter a path to the program using an example.
Suppose we want to input the following path:

$b=\young(234) \otimes \young(123,234,356)\otimes \young( 2,4)\otimes
\young(245,356)$   for the algebra of type $A_5$.
          
This path has 4 parts and the rank $n=5$. 
The input for this example will be
\begin{verbatim}
5 4            [First entry is n=5, 2nd is number of parts=4]
2 3 4         [This is the 1st part of the path]
0               [0 separates the parts]
1 2 3
2 3 4          [2nd part]
3 5 6
0               [0 separates the parts]
2               [3rd part]
4
0               [0 separates the parts]
2 4 5          [4th part of the path]
3 5 6
0               [0 to indicate the end of the path]
\end{verbatim}

{\bf WARNING:} Do not leave any extra blank space at the end of a line. The program 
will read the input incorrectly in that case and you will get a wrong answer.

{\bf Output:}
Now the output for our example is
\begin{verbatim}
n=5
Given path is:
---------------------------------------
 2  3  4
---------------------------------------
 1  2  3
 2  3  4
 3  5  6
---------------------------------------
 2
 4
---------------------------------------
 2  4  5
 3  5  6
---------------------------------------

Corresponding rigged configuration is :
---------------------------------------
(1)
 ___ ___ ___ ___ ___ ___ ___
|   |   |   |   |   |   | -4| -4
 --- --- --- --- --- --- ---
|   | -1| 0
 --- ---
---------------------------------------
(2)
 ___ ___ ___ ___ ___ ___
|   |   |   |   |   | -1| -1
 --- --- --- --- --- ---
|   |   |  0| 1
 --- --- ---
|   |  1| 1
 --- ---
---------------------------------------
(3)
 ___ ___ ___ ___
|   |   |   | -1| -1
 --- --- --- ---
|   |   |  0| 0
 --- --- ---
|   | -1| 0
 --- ---
---------------------------------------
(4)
 ___ ___ ___
|   |   | -1| 0
 --- --- ---
|   | -1| 0
 --- ---
---------------------------------------
(5)
 ___ ___
|   |  0| 0
 --- ---
---------------------------------------
Statistic = 5
\end{verbatim}

For the path, the dotted lines separates the parts of the path.
For the rigged configuration the program gives the component number
for the rigged partition and separates the different components with
dotted lines. In the end, the program gives the statistics corresponding 
to the path and the rigged configuration. We proved in 
Chapter~\ref{chap:kostka} that the statistics for
the path and the rigged configuration are preserved under the bijection.

\section{Program: $inverse\_bijection.c$}\label{prog:inverse_bijection}

The program called $inverse\_bijection.c$ computes the inverse map of $\Phi$.
It takes an unrestricted rigged configuration $(\nu,J)\in \RC(L,\la)$ as an input
and finds the image under the inverse bijection. The output is a path in
$\Path(B,\la)$.

{\bf Input:}
Let us explain the input file with an example. The input file is called $inputrigged$
for this program. Suppose we want to find the image of the rigged configuration
$$(\nu,J)=\yngrc(4,-1,3,0) \quad \yngrc(4,0,1,-1,1,-1) \quad \yngrc(3,-1,1,0)
\quad \yngrc(2,-1) $$
with $n=5$ and $\mu=((1),(1),(3,3),(2,2,2))$. 
To enter components of  $mu$ we use the dimension of the rectangular box. For 
example, the third component of $\mu$ in our example is entered as 2 3 indicating
a 2 by 3 rectangle.

The input file for this example is
\begin{verbatim}
5 4        [n=5, 4 is the number of components in mu]
1 1        [first component of mu] 
1 1        [2nd component of mu]
2 3        [3rd component of mu]
3 2        [4th component of mu]
4 1        [first rigged partition]
-1 0       [riggings for respective parts right below]
4 3        [2nd rigged partition]
-1 0       [riggings for respective parts]
4 1 1     [3rd rigged partition]
0 -1 -1   [riggings for respective parts]
3 1        [4th rigged partition]
-1 0       [riggings for respective parts]
2           [5th rigged partition]
-1          [riggings for respective parts]
\end{verbatim}
{\bf WARNING:} As in the previous cases, do not leave any extra blank 
space at the end of a line. The program will get confused and  
will give a wrong result.

{\bf Output:}
The output for our example is:
\begin{verbatim}
n = 5  L= 4
mu 
1 1 
1 1 
2 3 
3 2 

Given rigged configuration is:
---------------------------------------
(1)
 ___ ___ ___ ___
|   |   |   | -1| -1
 --- --- --- ---
|  0| 0
 ---
---------------------------------------
(2)
 ___ ___ ___ ___
|   |   |   | -1| 0
 --- --- --- ---
|   |   |  0| 0
 --- --- ---
---------------------------------------
(3)
 ___ ___ ___ ___
|   |   |   |  0| 1
 --- --- --- ---
| -1| -1
 ---
| -1| -1
 ---
---------------------------------------
(4)
 ___ ___ ___
|   |   | -1| -1
 --- --- ---
|  0| 0
 ---
---------------------------------------
(5)
 ___ ___
|   | -1| -1
 --- ---
---------------------------------------
The corresponding path is:
-------------------------
 3
-------------------------
 2
-------------------------
 1 3 5
 2 4 6
-------------------------
 1 2
 3 4
 5 6
-------------------------
\end{verbatim}

Note that the program first prints out the input data and then prints the image of the 
rigged configuration
under the inverse map which is an unrestricted path. The different parts
of the path are separated by dotted lines.
The output of the above example is the path
$$\young(3)\otimes \young(2) \otimes \young(135,246)\otimes \young(12,34,56).$$

\begin{subappendices}
\singlespacing
\section{Code for $allpaths\_bijection.c$}\label{allpaths}

\textit{This program computes the bijection for all unrestricted paths, finds all the
unrestricted rigged configurations for a fixed $\la$ and $\mu$. It also calculates
the statistics and the unrestricted Kostka polynomial corresponding to the $\la$ and
$\mu$.}

\begin{verbatim}
#include <stdio.h>
#define UNUSED 9999
#define RIGSIZE 20
int n, l, num_shapes;
int lambda[100];
int tab_shape[100];
int tableau[100][100];
int r, tab_indx, num_rc_lb_tab;
int *cum_lambda;
int *new_lambda;
int bigL [RIGSIZE][RIGSIZE];
int curL [RIGSIZE][RIGSIZE];
int path_index;
int tblu_index;
int num_paths;
int exp[1000];

FILE *fp;
class shape_class;

class tblu_row {
public:
    int *col;
    int num_col;
    tblu_row(int c);
    void print_row();
};
tblu_row::tblu_row(int c):num_col(c) {
    col = new int[c];
    for (int i=0; i < c; i++) col[i] = UNUSED;
}
void tblu_row::print_row() {
    if (col[0] == UNUSED) return;
    int i = 0;
    while (col[i] != UNUSED && i < num_col) {
       fprintf(stderr, "%2d ", col[i]);
       i++;
    }
    fprintf (stderr, "\n");
}
\end{verbatim}
\textit{ A doubly linked list of objects of type $tblu\_class$ makes up a path.
 Each object of type $tblu\_class$ represents a tableau which is a part 
of a path.}
\begin{verbatim}
class tblu_class {
public:
    int tblu_id;
    tblu_row* row;
    int* tab_lambda;
    int num_row;
    tblu_class* next;
    tblu_class* prev;
    shape_class* tblu_shape; // pointer to the mu
         // from which we got the shape
    tblu_class(int r, int c);
    void print_tblu();
};
tblu_class::tblu_class(int r, int c):num_row(r) {
    row = new tblu_row [r](c);
    tab_lambda = new int[n+1];
    for (int i=0; i<=n; i++) tab_lambda[i] = 0;
    next = NULL;
    prev = NULL;
    tblu_shape = NULL;
}
\end{verbatim}
\textit{Prints a tableau}
\begin{verbatim}
void tblu_class::print_tblu(){
    fprintf (stderr,"------------------------\n");
    for (int i=0; i < num_row; i++) {
       row[i].print_row();
    }
}
tblu_class *tblu_list;
tblu_class *tblu_list_end;
typedef tblu_class* tblu_class_ptr;
tblu_class_ptr *tblu_array;
class shape_class {
public:
    int* shape;
    int num_col;
    int num_row;
    shape_class* prev;
    shape_class* next;
    tblu_class* first_tblu;
    shape_class(int ncol, int nrow);
};
shape_class::shape_class (int ncol, int nrow){
    num_col = ncol;
    num_row = nrow;
    shape = new int[ncol];
    for (int i=0; i<ncol; i++) shape[i]=nrow;
    first_tblu = NULL;
    prev = NULL;
    next = NULL;
}
shape_class* shape_list;
shape_class* shape_list_end;
\end{verbatim}
\textit{
An object of type path$\_$class represents a path. A doubly linked 
list of objects of type path$\_$class has all the unrestricted paths and 
the corresponding rigged configuration for each path. 
}
\begin{verbatim}
class path_class {
public:
    int path_len;
    int index;
    int rigged[RIGSIZE][5][RIGSIZE]; 
              // the rigged set for this path
    tblu_class_ptr *path;    
              // the array of pointers to tableaux
    path_class* next;
    path_class* prev;
    int cocharge;
    int Energy;
    path_class(int path_len);
    void print_path();
    void path_class::reset_flags(int pathi);
    void path_class::print_rigged_for_this_path();
    int path_class::
        find_largest_inside_outside_others 
        (int index,int old_largest_index,int pathi);
    int path_class::
        find_largest_inside_outside_first 
        (int index, int pathi) ;
    void path_class::add_new_col
             (int index, int pathi);
    void path_class::init_unused_rigged
        (int index, int pathi);
    void path_class::add_to_rigged 
         (int index, int column_index, int pathi) ;
    int path_class::num_box_1k_col 
        (unsigned int i, int k, int pathi);
    int path_class::add_box_to_rigged 
        (int index, int begin, 
        int old_largest_index, int pathi) ;
    int path_class::second_func(int part_size, 
        int rig_num);
    void path_class::calc_outer_label 
         (int rig_num, int pathi);
    void path_class::calc_inner_label 
         (int i, int pathi);
    void path_class::insert_element_to_rigged 
        (int num, int row_indx, int col_indx,
         int  nrow, int ncol);
    void path_class::insert_tableau_to_rigged 
        (tblu_class* cur_tblu);
    void path_class::build_rigged_for_path () ;
    void path_class::calculate_cocharge();
    int path_class::alpha(int k, int i);                                                                          
};
path_class::path_class(int len):path_len(len) {
    path = new tblu_class_ptr [len];
    for (int i=0; i<len; i++) path[i] = NULL;
    next = NULL;
    prev = NULL;
    index=0;
    for (int i=0; i < RIGSIZE; i++) {
       for (int j = 0; j < 3; j++)
          for (int k = 0;  k < RIGSIZE; k++)
             rigged [i][j][k] = UNUSED;
       for (int j = 3; j < 5; j++)
          for (int k = 0;  k < RIGSIZE; k++)
             rigged [i][j][k] = 0;
    }
};
\end{verbatim}
\textit{Prints a path.}
\begin{verbatim}
void path_class::print_path() {
    for (int i=0; i <path_len; i++) {
        if (path[i] != NULL) path[i]->print_tblu();
    }
}
typedef path_class* path_class_ptr;
path_class_ptr tmp_path;
path_class_ptr path_list;
path_class_ptr path_list_end;
path_class_ptr *path_array;
int *print_order;

void reset_tableau() {
    for (int i=0;i<RIGSIZE; i++){
        for (int j=0; j<RIGSIZE;j++){
            tableau[i][j]=UNUSED;
        }
    }    
}
void initialize_lambda() {
    int i, j, k,m;
    for (i=0; i < RIGSIZE; i++) {
        lambda[i] = -1;
       tab_shape[i]=-1;
    }
    for (i=0; i<1000; i++);
        exp[i]=0;
    reset_tableau(); 
}
\end{verbatim}
\textit{Reads the input file.}
\begin{verbatim}
void read_input(){
    int i, tmp;
    tmp = UNUSED;
    i = 0;
    l = 0;
    fp = fopen ("input_allpaths","rw");
    fscanf(fp,"%d\n", &n);
    while (i< n+1) {
        fscanf (fp, "%d", &tmp);
        lambda[i] = tmp;
        l = l + tmp;
        i++;
    }
}
\end{verbatim}
\textit{Prints the input data: $n$, $\lambda$ and $\mu$.}
\begin{verbatim}
void print_input() {
    int i;
    fprintf (stderr, "n = %d \n", n);
    fprintf (stderr, "Lambda is:  ");
    for (i=0; i <=n; i++) {
    if (lambda[i] == -1) break;
        fprintf (stderr, "%d ", lambda[i]);
    }
    fprintf (stderr, "\n");
}
\end{verbatim}
\textit{Copies the tableau constructed to the tableau list.} 
\begin{verbatim}
void print_and_copy_tableau(int k, int nrow, 
       int ncol, shape_class*  shape_obj ) {
    int i, j;
    i=0;
    tblu_class *my_tblu = 
              new tblu_class(nrow,ncol);
    my_tblu->tblu_id = tblu_index;
    tblu_index += 1;
    while (tableau[i][0] != UNUSED){
       j=0;
       while ( tableau[i][j]!=UNUSED){
          my_tblu->row[i].col[j] = tableau[i][j];
          my_tblu->tab_lambda[tableau[i][j] - 1] = 
          my_tblu->tab_lambda[tableau[i][j] - 1] 
                                             + 1;
          j++;
       }
       i++;
    }
    if (shape_obj->first_tblu == NULL) 
       shape_obj->first_tblu = my_tblu;
    my_tblu->tblu_shape = shape_obj;
    if (tblu_list == NULL){ 
        tblu_list = my_tblu;
        tblu_list_end = my_tblu;
    } else {
        tblu_list_end->next = my_tblu;
        my_tblu->prev = tblu_list_end;
        tblu_list_end = my_tblu;
    }         
}
\end{verbatim}  
\textit{Builds a  tableau recursively.}
\begin{verbatim}
void build_tableau(int m,int row,int col, int nrow, 
       int ncol, shape_class* shape_obj){
    int k,p,q,m1,h,h1,valid,num,a,b,c ; 

    if ( (col == 0 && row == 0) ||
       ((col == 0) && (m > tableau[row - 1][col]) )
       || (row == 0 && m >= tableau[row][col-1]) || 
       (row > 0 && col > 0 && m > 
          tableau[row-1][col] && m >= 
          tableau[row][col-1]))
    { 
        tableau[row][col]=m;
        if ( row == tab_shape[col]-1) {
           if ( tab_shape[col+1] < 1) {
              print_and_copy_tableau(tab_indx, 
                nrow, ncol, shape_obj);
              tab_indx=tab_indx+1;
              return;
           }
           else{ 
               for ( k=1; k <= n+1; k++) {
                  build_tableau(k,0,
                   col+1,nrow,ncol, shape_obj);
               }     
        }
        } else{ 
          for ( k=m+1; k<=n+1; k++) {
             build_tableau(k,row+1,col,nrow,
                    ncol, shape_obj);
           } 
        }
          
    }     
}
\end{verbatim}
\textit{ Finds all possible  tableaux of a given shape.}
\begin{verbatim}
void find_tableau(){
  int i,j,k,h,tmp;
  int nrow, ncol;
  tab_indx=0;num_rc_lb_tab=0;
  num_shapes = 0;
  for (i=0;i<RIGSIZE;i++){
    tab_shape[i]=0;
  } 
  j = 0; 
  tblu_list = NULL;
  tblu_list_end = NULL;
  shape_list = NULL;
  shape_list_end = NULL;
  tblu_index = 0;
  fprintf (stderr, "mu is :\n");
  while (fscanf (fp, "%d", &tmp) != EOF){
     if (tmp == 0)  {
        k = 0;
        shape_class* my_shape_obj 
          = new shape_class(k,tab_shape[0]);
        if (shape_list == NULL){
           shape_list = my_shape_obj; 
           shape_list_end = my_shape_obj; 
        } else {
           my_shape_obj->prev = shape_list_end;
           shape_list_end->next = my_shape_obj;
           shape_list_end = my_shape_obj;
        }
             
        fprintf (stderr, "..........\n");
  
        while (tab_shape[k] != 0) {
           fprintf (stderr, "%2d ", tab_shape[k]);
           k++;
        }
        num_shapes++;
        fprintf (stderr, "\n");
        nrow = tab_shape[0];
        ncol = k;
        reset_tableau();
        for (int i = 1; i <= n+1; i++){
           build_tableau(i, 0, 0, 
             nrow, ncol, my_shape_obj);
        }
        for (int i = 0; i < RIGSIZE; i++) 
           tab_shape[i] = 0;
        j = 0;
     } else {
        tab_shape[j] = tmp;
        j++;
     }
  }
}
\end{verbatim}
\textit{Finds the possible parts (which are tableaux) of a path with the given $\mu$ and 
$\lambda$.}
\begin{verbatim} 
void find_path_element(int position, shape_class* 
       cur_shape, tblu_class* cur_tblu) {

    if (cur_tblu == NULL) {
        return;
    }

    int *this_lambda = cur_tblu->tab_lambda;
    bool satisfy = true;
    for (int i=0; i<= n; i++) {
        new_lambda[i] = 
               cum_lambda[i] + this_lambda[i];
        if (new_lambda[i] > lambda[i]) {
           satisfy = false;
           break;
        }
    }
    if (satisfy == false) {
       return;
    } else {
       for (int i=0; i<=n; i++) 
          cum_lambda[i] = new_lambda[i];
       tmp_path->path[position] = cur_tblu;
       if (position == num_shapes - 1) {
          bool found= true;
          for (int i=0; i<=n; i++)
            if (cum_lambda[i] != lambda[i]) 
               found = false;
          if (found) { 
               path_class* tmp_path_list = 
                 new path_class (num_shapes);
               path_index=path_index+1;
               tmp_path_list -> index=path_index;
               for (int i=0; i<num_shapes; i++) {
                   tmp_path_list->path[i] = 
                     tmp_path->path[i];
               }
               if (path_list == NULL) {
                  path_list = tmp_path_list;
                  path_list_end = tmp_path_list;
               } else {
                  path_list_end->next = 
                                  tmp_path_list;
                  tmp_path_list->prev = 
                                  path_list_end;
                  path_list_end = tmp_path_list;
               }
          }
       } else {
          tblu_class* this_tblu = 
               cur_shape->next->first_tblu;
          shape_class* this_shape =cur_shape->next;
          while (this_tblu != NULL && 
             this_tblu->tblu_shape == this_shape){
             find_path_element(position +1,
                this_shape, this_tblu);
             this_tblu = this_tblu->next;
          }
       }
       for (int i=0; i<=n; i++) cum_lambda[i] =
          cum_lambda[i] - this_lambda[i];
       tmp_path->path[position] = NULL;
    }
}
\end{verbatim}
\textit{Builds paths of shape $\mu$ and weight $\lambda$.}
\begin{verbatim}
void build_paths() {
    cum_lambda = new int [n+1];
    new_lambda = new int [n+1];
    path_index = 0;
    for (int i = 0; i <=n; i++) {
       cum_lambda[i] = 0;
       new_lambda[i] = 0;
    }
    tmp_path = new path_class (num_shapes);
    path_list = NULL;
    path_list_end = NULL;
    tblu_class* this_tblu = shape_list->first_tblu;
    while (this_tblu != NULL && 
           this_tblu->tblu_shape == shape_list) {
       find_path_element (0, shape_list, this_tblu);
       this_tblu = this_tblu->next;
    }

    path_class* tmp_path_list = path_list;
    while (tmp_path_list != NULL) {

       tmp_path_list = tmp_path_list->next;
    }

}


void path_class::reset_flags (int pathi) {
   int i, k;
   for (i=0; i < RIGSIZE; i++) {
      for (k=0; k < RIGSIZE; k++) {
          rigged [i][3][k] = 0;
      }
   }
}
\end{verbatim}
\textit{ This finds the cocharge for a given rigged configuration.}
\begin{verbatim}
int path_class::alpha(int k, int i){
       int num_coln,j;
       num_coln=0;
       if (k>=n) num_coln=0; 
     else{
       for (j=0; j<RIGSIZE;j++){
          if (rigged[k][0][0]==UNUSED){
               num_coln=0 ;
               break;
          }
          else{
             if (rigged[k][0][j]!=UNUSED){ 
                  if (rigged[k][0][j]>=i+1){
                     num_coln=num_coln+1;
             }
       }
     }   
        }
     }
       return num_coln;
}
void path_class::calculate_cocharge(){
  int k,j,i,sum,cosum;
  sum=0;
  for (k=0;k<=n-1;k++){
    if (rigged[k][0][0]!=UNUSED){
      for (i=0;i< rigged[k][0][0];i++){
      sum=sum+alpha(k,i)*(alpha(k,i)-alpha(k+1,i));
      }  
    }
  }     
  cosum=sum;
  for (k=0;k<=n-1;k++){
    for (j=0;j<RIGSIZE;j++){
      if (rigged[k][2][j]==UNUSED) break;
      if ( rigged[k][2][j]!=UNUSED){
        cosum=cosum + rigged[k][2][j];
      }
    }
  }
  cocharge=cosum;
  fprintf (stderr, "Statistic is = %d   \n", cosum );
  exp[cosum]=exp[cosum]+1;  
}
\end{verbatim}
\textit{Prints the rigged configuration corresponding to a path.}
\begin{verbatim}
void path_class::print_rigged_for_this_path() {
   int i, j, k, l,a,b;
   
   for (i = 0; i < n ; i++) {
      if (rigged[i][0][0]==UNUSED) {
         fprintf (stderr,"-------------------\n");
         fprintf(stderr,"(%d)   Empty\n", i+1);
      } else {   
         if (rigged[i][0][0] != UNUSED) {
            fprintf (stderr,"----------------\n");
            fprintf(stderr, "(%d)\n", i+1);

            j=0;
            if (rigged[i][0][j] != UNUSED)
               for (k=0; k <rigged[i][0][j]; k++)
                   fprintf (stderr, " ___");
            fprintf (stderr,"\n");
            while (rigged[i][0][j] !=UNUSED){
                k=rigged[i][0][j];
                for (l=0; l<k-1;l++){
                   fprintf (stderr, "|   ");
                } 
                if (l==k-1) fprintf (stderr, 
                   "| %2d",rigged[i][2][j]); 
   
                fprintf (stderr, 
                     "| %d\n",rigged[i][1][j]);
                for (l=0; l<k; l++){
                   fprintf (stderr, " ---");
                }
                fprintf (stderr,"\n");
                j++;
            }  
         }
      }
 
   }
   calculate_cocharge();
   fprintf (stderr, "********************\n");
   fprintf (stderr,"\n");          
}
\end{verbatim}
\textit{Finds the largest singular string in a rigged partition other than the first one
which is bigger or equal to the string selected in the previous rigged partition.}
\begin{verbatim}
int path_class::find_largest_inside_outside_others
         (int index, int old_largest_index, 
          int pathi) {
   int i = 0;
   int largest_index = UNUSED;
   while ((rigged[index][0][i] != UNUSED) && 
           i < RIGSIZE) {
      if ((rigged[index][1][i] == 
                     rigged[index][2][i]) &&
          (rigged[index][1][i] != UNUSED) &&
          (rigged[index][0][i] 
                          <= old_largest_index)) {
         largest_index = i;
         break;
      }
      i++;
   }
   return largest_index;
}
\end{verbatim}
\textit{Finds the largest singular string in the first rigged partition}
\begin{verbatim}
int path_class::find_largest_inside_outside_first
                    (int index, int pathi) {
   int i = 0;
   int largest_index = UNUSED;
   while ((rigged[index][0][i] != UNUSED) && 
                                i < RIGSIZE) {
      if ((rigged[index][1][i] == 
                  rigged[index][2][i]) &&
          (rigged[index][1][i] != UNUSED)) {
         largest_index = i;
         break;
      }
      i++;
   }
   return largest_index;
}
void path_class::add_new_col(int index, int pathi){
   int i=0;
   while (rigged[index][0][i] != UNUSED) i++;
   rigged[index][0][i] = 1;
   rigged [index][3][i] = 1;
}
void path_class::init_unused_rigged 
                     (int index, int pathi) {
   rigged [index][0][0] = 1;
   rigged [index][3][0] = 1;
}
void path_class::add_to_rigged (int index, 
               int column_index, int pathi) {
   rigged [index][0][column_index] += 1;
   rigged [index][3][column_index] = 1;
}
\end{verbatim}
\textit{Calculates the number of boxes in the first $k$ columns of a rigged partition.}
\begin{verbatim}
int path_class::num_box_1k_col 
          (unsigned int i, int k, int pathi){
   int j, l;
   int num_boxes = 0;
   for (l=1; l <= k; l++) {
      for (j = 0; j < RIGSIZE; j++){
         if (rigged[i][0][j] == UNUSED) break;
         if (rigged[i][0][j] >= l){
             num_boxes += 1;
         }
      }
   }
   return num_boxes;
}
\end{verbatim}
\textit{This adds a box to a rigged partition while doing the bijection.}
\begin{verbatim}
int path_class::add_box_to_rigged (int index, 
     int begin, int old_largest_index, int pathi){
   int largest_index;
   if (index == begin) {
      largest_index = 
         find_largest_inside_outside_first(index, 
                                           pathi);
   } else {
      largest_index =
         find_largest_inside_outside_others(index, 
                old_largest_index, pathi);
   }
   if (largest_index == UNUSED) {
      if (rigged[index][0][0] == UNUSED)
         init_unused_rigged (index, pathi);
      else 
         add_new_col (index, pathi);
      return 0;
   } else {
      add_to_rigged (index, largest_index, pathi);
      return (rigged [index][0][largest_index] - 1);
   }
}
\end{verbatim}
\textit{This calculates the second function in the definition of vacancy numbers.}
\begin{verbatim}
int path_class::second_func
        (int part_size, int rig_num) {
   int sum = 0;
   for (int i=1; i<RIGSIZE; i++) {
      if (curL[rig_num+1][i] != 0) {
         int minimum = part_size;
         if (i < part_size) minimum = i;
         sum =sum + minimum * curL[rig_num+1][i];
      }
   } 
   return sum;
}
\end{verbatim}
\textit{This calculates the vacancy numbers.}
\begin{verbatim}
void path_class::calc_outer_label (int rig_num, 
                                   int pathi){
   int part_num=0;
   int part_size, p;

   if (rig_num == 0)  {
      for (int part_num=0; part_num < RIGSIZE; 
                                part_num++) {
         part_size = rigged[0][0][part_num];
         if(part_size == UNUSED) break;
         p = (-2*num_box_1k_col (0, part_size, 
                                    pathi)) +
           (num_box_1k_col (1, part_size, pathi)) +
           second_func(part_size, rig_num);
         rigged [0][1][part_num] = p;
      }
   } else {
      for (part_num=0; part_num < RIGSIZE; 
                                       rt_num++) {
         part_size = rigged[rig_num][0][part_num];
         if(part_size == UNUSED) break;
         p = -2*num_box_1k_col (rig_num, part_size, 
                                           pathi) +
             num_box_1k_col (rig_num-1, part_size, 
                                           pathi) +
             num_box_1k_col (rig_num+1, part_size, 
                                           pathi) +
             second_func(part_size, rig_num);
         rigged [rig_num][1][part_num] = p;
      }

   }
}
\end{verbatim}
\textit{This calculates the labels or the riggings.}
\begin{verbatim}
void path_class::calc_inner_label 
            (int i, int pathi) {
   int j,k;
   int tmp;
   for (j = 0; j < RIGSIZE; j++) {
      if (rigged[i][0][j] == UNUSED) break;
      if (rigged[i][3][j] == 1) {
         rigged [i][2][j] = rigged [i][1][j];
      }
   }
   for (j = 0; j < RIGSIZE; j++) {
      for (k = 1; k < RIGSIZE; k++) {
         if (rigged[i][0][k] == UNUSED) break;
         if (rigged[i][0][k] == rigged[i][0][k-1] 
           && rigged[i][1][k] == rigged[i][1][k-1]
           && rigged[i][2][k-1] < rigged[i][2][k]) 
           {
            tmp = rigged[i][2][k-1];
            rigged[i][2][k-1] = rigged[i][2][k];
            rigged[i][2][k] = tmp;
         }
      }
   }
}
\end{verbatim}
\textit{This inserts each element of a part in the path into the bijection}
\begin{verbatim}
void path_class::insert_element_to_rigged 
         (int num, int row_indx, int col_indx,
                         int  nrow, int ncol) {
   int old_largest_index;
   reset_flags(0);
   for (int i = (num-2); i >= row_indx; i--) {
      old_largest_index = 
      add_box_to_rigged(i, num-2, old_largest_index, 0);
   }
   for (int i=0; i < RIGSIZE; i++) 
      for (int j = 0; j < RIGSIZE; j++)
         curL[i][j] = bigL[i][j];
   if (row_indx == nrow - 1) {
      curL [nrow][1] += 1;
      if (col_indx != ncol - 1) {
         curL [nrow][ncol - col_indx - 1] += 1;
      }
   } else {
      curL [nrow][ncol - col_indx - 1] += 1;
      curL [row_indx + 1][1] += 1;
   }
   for (int i = 0; i < n; i++) {
      calc_outer_label (i, 0);
   }
   for (int i = 0; i < n; i++) {
      calc_inner_label (i, 0);
   }     
   if (row_indx == nrow - 1) {
      for (int i=0; i < RIGSIZE; i++) 
         for (int j = 0; j < RIGSIZE; j++)
            curL[i][j] = bigL[i][j];

      // update curL - this is different from above
      curL [nrow][ncol - col_indx] += 1;
      for (int i = 0; i < n; i++) {
         calc_outer_label (i, 0);
      }
   }
}
\end{verbatim}
\textit{This inserts each part of a path into the bijection}
\begin{verbatim}
void path_class::insert_tableau_to_rigged 
         (tblu_class* cur_tblu) {
   int nrow = cur_tblu->num_row;
   int ncol = cur_tblu->row[0].num_col;
   for (int j = ncol-1; j >=0; j--) {
     for (int i = 0; i < nrow; i++) {
        insert_element_to_rigged
         (cur_tblu->row[i].col[j], i,j,nrow, ncol);
     }
   }
   bigL [nrow][ncol] = bigL [nrow][ncol] + 1;
   for (int i = 0; i < RIGSIZE; i++)
      for (int j = 0; j < RIGSIZE; j++) 
         curL[i][j] = bigL[i][j];
}
\end{verbatim}
\textit{ This insert all the parts of a path to the bijection}
\begin{verbatim}
void path_class::build_rigged_for_path () {
   int j;
   for (int i=0; i < RIGSIZE; i++) {
      for (j =0; j < RIGSIZE; j++)
         bigL[i][j] = 0;
      curL[i][j] = 0;
   }
   for (int i = path_len - 1; i >= 0; i--) {
      insert_tableau_to_rigged (path[i]);
   }  
}
\end{verbatim}
\textit{Sorts the rigged configurations in a order so that the 
configurations with the same shape appears together.}
\begin{verbatim}
void sort_rigged(){
  int i,j,k,a,b,p,T,m,h,a1,b1;
  a=0;p=0;a1=0;

  for (int a=0; a < num_paths; a++) {
     if  (path_array[a]->rigged[0][4][0] != 0) 
              continue;
     print_order[p] = a;
     p = p + 1;
     for (int i=a+1; i < num_paths; i++){
        if (path_array[i]->rigged[0][4][0] != 0)
           continue;

        T = 0;
        for (b=0;b<=n-1;b++){
        // pick the b-th rig-element of 
        // a-th path with i-th path
           k=0;
           while(path_array[a]->rigged[b][0][k]
                                      !=UNUSED) {
             if (path_array[i]->rigged[b][0][k] !=
                path_array[a]->rigged[b][0][k]){
                  T=1;break;
             }
             if (path_array[i]->rigged[b][4][k]
                                           !=0) {
                T=1;break;
             }
             k++;
           }
           // Make sure if we exited while 
           // loop because a[..] unused
           if (path_array[a]->rigged[b][0][k] 
               != path_array[i]->rigged[b][0][k]) 
              T = 1;
           // if unequal quit searching 
           if (T == 1) break; 
        } 
        // if rig-element comparison failed 
        // quit path comparison
        if (T==1) { 
           continue;
        }
        print_order[p]=i;
        path_array[i]->rigged[0][4][0]=1;
        p=p+1;
     }   // inner for loop - i
  }       // outer for loop - a
}  
\end{verbatim}
\textit{ Calculates the rigged configurations corresponding to each of the
possible paths.}
\begin{verbatim}
void build_rigged() {
    path_class* tmp_path_list = path_list;
    num_paths = 0;
    while (tmp_path_list != NULL) {
       tmp_path_list->build_rigged_for_path ();
       tmp_path_list = tmp_path_list->next;
       num_paths++;
    }
    path_array = new path_class_ptr [num_paths];
    print_order = new int [num_paths];
    for (int i=0; i<num_paths; i++) {
       print_order[i] = i;
    }
    tmp_path_list = path_list;
    while (tmp_path_list != NULL) {
       path_array[tmp_path_list->index - 1] = 
                             tmp_path_list;
       tmp_path_list = tmp_path_list->next;
    }
    sort_rigged();
    fprintf(stderr,"-----------------------\n");
    fprintf(stderr,"There are %d unrestricted 
                      paths.\n", num_paths);
    fprintf(stderr,"-----------------------\n");
    fprintf(stderr,"\n");
    for (int i=0; i<num_paths; i++) {
       fprintf (stderr, "Path (%d): \n",i+1);
       path_array[print_order[i]]->print_path();
       fprintf (stderr, "\nCorresponding rigged 
                 configuration is:\n");
       path_array[print_order[i]]->
             print_rigged_for_this_path();

    }
    fprintf(stderr,
           "Unrestricted Kostka polynomial is: ");   
    int begin=1;   
    for (int i=0; i<1000; i++) {
       if (exp[i]!= 0) {
          if (!begin)
             fprintf(stderr, " + ");
          begin=0;
          fprintf(stderr, "%dq^%d ", exp[i],i);
       }
    }
    fprintf(stderr,"\n");    
}
\end{verbatim}
\textit{Main program.}
\begin{verbatim}
int main() {
   int i;
   initialize_lambda();  // initialization.
   read_input();          // this reads the input.
   print_input();       // this prints the input.
   find_tableau();// this finds all the tableaux of 
                  // all the shapes from mu
   build_paths(); // this finds all the paths for 
                  // given lambda and mu
   build_rigged();// this calculates all the rigged 
                // configurations via the bijection 
                // and sorts them in the order so 
                // that all the configurations 
                // with the same shape appear together. 
}
\end{verbatim}

\section{Code for the program $one\_path\_bij.tex.c$}

\textit{
This program does the bijection from the set of paths to the set of 
rigged configurations. 
Input data is a single path and the program calculates the corresponding
rigged configuration using the bijection.
It also calcutales the statistics.
}

\begin{verbatim}
#include <stdio.h>
#define UNUSED 9999 
#define RIGSIZE 50
int n,l;
int tab_shape[100];
int tableau[100][100];
int r;
int rigged[RIGSIZE][5][RIGSIZE];
int bigL [RIGSIZE][RIGSIZE];
int curL [RIGSIZE][RIGSIZE];
int path_index;
int tblu_index;
FILE *fp;

\end{verbatim}
\textit{ 
A doubly linked list of objects of type tblu$\_$class makes up a path.
Each object of type tblu$\_$class represents a tableau which is a part 
of a path.
}
\begin{verbatim}
class tblu_class {
public:
    int tblu_id;
    int** tb; // 2-dimensional array of integers 
                // holding the tableau
    int* tab_lambda;
    int num_row;
    int num_col;
    tblu_class* next;
    tblu_class* prev;
    tblu_class(int r, int c);
    void print_tblu();
};

tblu_class::tblu_class(int r, int c)
              :num_row(r),num_col(c) {
   tb = new int* [r];
   for (int i=0; i < r; i++) {
      tb[i] = new int [c];
      for (int j = 0; j < c; j++)
         tb [i][j] = UNUSED;
   }
   tab_lambda = new int[n+1];
   for (int i=0; i<=n; i++) tab_lambda[i] = 0;
   next = NULL;
   prev = NULL;
}
\end{verbatim}
\textit{Prints a tableau.}
\begin{verbatim}
void tblu_class::print_tblu(){
   fprintf (stderr, "------------------------\n");
   for (int i=0; i < num_row; i++) {
      for (int j=0; j < num_col; j++)
         fprintf (stderr, "%2d ", tb[i][j]);
      fprintf (stderr, "\n");
   }
}


typedef tblu_class* tblu_class_ptr;
tblu_class_ptr *tblu_array;

\end{verbatim}
\textit{
An object of type path$\_$class represents a path. In this
program it has only one object and the corresponding
rigged configuration. 
}
\begin{verbatim}
class path_class {
public:
    int path_len;
    int index;
    int rigged[RIGSIZE][5][RIGSIZE]; 
      // this is the rigged set for this path
    tblu_class_ptr *path;  // this is the array 
       // of pointers to tableaux
    path_class* next;
    path_class* prev;
    int cocharge;
    path_class(int path_len);
    void print_path();
    void path_class::reset_flags(int pathi);
    void path_class::print_rigged_for_this_path();
    int path_class::find_largest_inside_outside_others
        (int index, int old_largest_index, int pathi);
    int path_class::find_largest_inside_outside_first 
                               (int index, int pathi) ;
    void path_class::add_new_col
           (int index, int pathi);
    void path_class::init_unused_rigged
           (int index, int pathi);
    void path_class::add_to_rigged (int index, 
                      int column_index, int pathi) ;
    int path_class::num_box_1k_col (unsigned int i, 
                                  int k, int pathi);
    int path_class::add_box_to_rigged (int index, 
           int begin, int old_largest_index, 
           int pathi) ;
    int path_class::second_func
           (int part_size, int rig_num);
    void path_class::calc_outer_label 
           (int rig_num, int pathi);
    void path_class::calc_inner_label 
           (int i, int pathi);
    void path_class::insert_element_to_rigged 
           (int num, int row_indx, int col_indx,
            int  nrow, int ncol);
    void path_class::insert_tableau_to_rigged 
           (tblu_class* cur_tblu);
    void path_class::build_rigged_for_path () ;
    void path_class::calculate_cocharge();
    int path_class::alpha(int k, int i);

};


path_class::path_class(int len):path_len(len) {
    path = new tblu_class_ptr [len];
    for (int i=0; i<len; i++) path[i] = NULL;
        next = NULL;
        prev = NULL;
        index=0;
        for (int i=0; i < RIGSIZE; i++) {
            for (int j = 0; j < 3; j++)
                for (int k = 0;  k < RIGSIZE; k++)
                    rigged [i][j][k] = UNUSED;
            for (int j = 3; j < 5; j++)
                for (int k = 0;  k < RIGSIZE; k++)
                    rigged [i][j][k] = 0;
       }
};


void path_class::print_path() {
    fprintf (stderr, "Given path is:\n");
    for (int i=0; i <path_len; i++) {
        if (path[i] != NULL) path[i]->print_tblu();
    }
}


typedef path_class* path_class_ptr;
path_class_ptr input_path;


void reset_tableau() {
    for (int i=0;i<RIGSIZE; i++){
        for (int j=0; j<RIGSIZE;j++){
            tableau[i][j]=UNUSED;
        }
    }    
}  

\end{verbatim}
\textit{
This reads the input file.
1st number we input is "$n$" which is the number of nodes in the Dynkin
diagram of type A$\_${n}.
2nd number we input is the number of parts in the path, called path length.
Then we input each part of the path which is a tableau,
we seperate the parts by putting a "0".
At the end of the last part we put a "0" to ensure the end of the path.
}
\begin{verbatim}
void read_input(){

   int i, j, k, tmp, c;
   int path_len, path_index;
   int col;
   tmp = UNUSED;
   path_index = 0;
   i = 0;
   j = 0;
   k = 0;
   col = 0;
   l = 0;
   fp = fopen ("inputpath","rw");
   fscanf(fp,"%d %d\n", &n, &path_len);
   input_path = new path_class (path_len);
   reset_tableau();
   while (fscanf (fp, "%d", &tmp) != EOF) {
      if (tmp == 0) {
         tblu_class *my_tblu = 
            new tblu_class(i,col);
         i = 0;
         my_tblu->tblu_id = tblu_index;
         tblu_index += 1;
         while (tableau[i][0] != UNUSED){
               
            j=0;
	    while ( tableau[i][j]!=UNUSED){
	       my_tblu->tb[i][j] = tableau[i][j];
	       my_tblu->tab_lambda[tableau[i][j] - 1] = 
                my_tblu->tab_lambda[tableau[i][j] - 1] + 1;
               j++;
	    }
	    i++;
	 }
         input_path->path[path_index] = my_tblu;
         reset_tableau ();
         path_index += 1;
	 i = 0; j = 0;
         continue;
      }
      tableau[i][j] = tmp;
      j += 1;
      c = fgetc(fp);
      if (c == '\n') {
         i += 1;
         col = j;
         j = 0;
      }
   }
}



void path_class::reset_flags (int pathi) {
   int i, k;
   for (i=0; i < RIGSIZE; i++) {
      for (k=0; k < RIGSIZE; k++) {
         rigged [i][3][k] = 0;
      }
   }
}

\end{verbatim}
\textit{
This prints the rigged configuration obtained from the bijection and prints the 
corresponding statistic.}
\begin{verbatim}
void path_class::print_rigged_for_this_path() {
   int i, j, k, l,a,b;
   for (i = 0; i < n; i++) {
      if (rigged[i][0][0]==UNUSED) {
         fprintf (stderr,"--------------------\n");
         fprintf(stderr,"(%d)   Empty\n", i+1);
      } 
      else {   
        if (rigged[i][0][0] != UNUSED) {
           fprintf (stderr,"------------------\n");
           fprintf(stderr, "(%d)\n", i+1);
           j=0;
           if (rigged[i][0][j] != UNUSED)
           for (k=0; k <rigged[i][0][j]; k++)
              fprintf (stderr, " ___");
           fprintf (stderr,"\n");
           while (rigged[i][0][j] !=UNUSED){
               k=rigged[i][0][j];
               for (l=0; l<k-1;l++){
                  fprintf (stderr, "|   ");
               } 
               if (l==k-1) fprintf 
                  (stderr, "| %2d",rigged[i][2][j]); 
   
               fprintf (stderr, 
                  "| %d\n",rigged[i][1][j]);
               for (l=0; l<k; l++){
                  fprintf (stderr, " ---");
               }
               fprintf (stderr,"\n");
               j++;
            }  
         }
       } 
    }
    fprintf (stderr,"--------------------------\n");
    calculate_cocharge();
}

\end{verbatim}
\textit{
Finds the largest singular string in a rigged partition other than the first one
which is bigger or equal to the string selected in the previous rigged partition
by $\delta$.}
\begin{verbatim}
int path_class::find_largest_inside_outside_others 
   (int index, int old_largest_index, int pathi) {
   int i = 0;
   int largest_index = UNUSED;
   while ((rigged[index][0][i] != UNUSED) && 
      i < RIGSIZE) {
      if ((rigged[index][1][i] == 
                        rigged[index][2][i]) 
            && (rigged[index][1][i] != UNUSED) &&
         (rigged[index][0][i] <= old_largest_index)) {
         largest_index = i;
         break;
      }
      i++;
   }
   return largest_index;
}

\end{verbatim}
\textit{
Finds the largest singular string in the first rigged partition.
}
\begin{verbatim}
int path_class::find_largest_inside_outside_first 
      (int index, int pathi) {
    int i = 0;
    int largest_index = UNUSED;
    while ((rigged[index][0][i] != UNUSED) 
             && i < RIGSIZE) {
       if ((rigged[index][1][i] == 
            rigged[index][2][i]) &&
           (rigged[index][1][i] != UNUSED)) {
          largest_index = i;
          break;
       }
       i++;
    }
    return largest_index;
}

void path_class::add_new_col (int index, int pathi)
{
    int i=0;
    while (rigged[index][0][i] != UNUSED) i++;
    rigged[index][0][i] = 1;
    rigged [index][3][i] = 1;
}

void path_class::init_unused_rigged (int index, 
      int pathi) {
    rigged [index][0][0] = 1;
    rigged [index][3][0] = 1;
}

void path_class::add_to_rigged (int index, 
     int column_index, int pathi) {
    rigged [index][0][column_index] += 1;
    rigged [index][3][column_index] = 1;
}

\end{verbatim}
\textit{Calculates the number of boxes in the first $k$ columns of a rigged partition.}
\begin{verbatim}
int path_class::num_box_1k_col (unsigned int i, 
     int k, int pathi){
    int j, l;
    int num_boxes = 0;
    for (l=1; l <= k; l++) {
       for (j = 0; j < RIGSIZE; j++){
          if (rigged[i][0][j] == UNUSED) break;
          if (rigged[i][0][j] >= l){
              num_boxes += 1;
          }
       }
    }
    return num_boxes;
}

\end{verbatim}
\textit{This adds a box to a rigged partition while doing the bijection.}
\begin{verbatim}
int path_class::add_box_to_rigged (int index, 
       int begin, int old_largest_index, 
       int pathi) {
    int largest_index;
    if (index == begin) {
       largest_index = 
         find_largest_inside_outside_first
                     (index, pathi);
    } else {
       largest_index = 
          find_largest_inside_outside_others 
            (index, old_largest_index, pathi);
    }
    if (largest_index == UNUSED) {
       if (rigged[index][0][0] == UNUSED)
          init_unused_rigged (index, pathi);
       else 
          add_new_col (index, pathi);
       return 0;
    } else {
       add_to_rigged (index, largest_index, pathi);
       return 
          (rigged [index][0][largest_index] - 1);
    }
}

\end{verbatim}
\textit{This calculates the second function in the definition of vacancy numbers.}
\begin{verbatim}
int path_class::second_func(int part_size, 
                            int rig_num) {
    int sum = 0;
    for (int i=1; i<RIGSIZE; i++) {
        if (curL[rig_num+1][i] != 0) {
        int minimum = part_size;
        if (i < part_size) minimum = i;
            sum =sum + minimum * curL[rig_num+1][i];
        }
    }
    return sum;
}

\end{verbatim}
\textit{This calculates the vacancy numbers.}
\begin{verbatim}
void path_class::calc_outer_label 
       (int rig_num, int pathi){
   int part_num=0;
   int part_size, p;

   if (rig_num == 0)  {
      for (int part_num=0; part_num < RIGSIZE; 
            part_num++) {
         part_size = rigged[0][0][part_num];
         if(part_size == UNUSED) break;
         p = (-2*num_box_1k_col (0, part_size, 
                                         pathi)) 
             + (num_box_1k_col (1, part_size,
                                         pathi)) 
             + second_func(part_size, rig_num);
         rigged [0][1][part_num] = p;
      }
   } else {
      for (part_num=0; part_num < RIGSIZE; 
           part_num++) {
         part_size = rigged[rig_num][0][part_num];
         if(part_size == UNUSED) break;
         p = -2*num_box_1k_col (rig_num, part_size, 
             pathi) + num_box_1k_col (rig_num-1, 
             part_size, pathi) + num_box_1k_col 
             (rig_num+1, part_size, pathi) +
             second_func(part_size, rig_num);
         rigged [rig_num][1][part_num] = p;
      }

   }
}

\end{verbatim}
\textit{This calculates the labels or the riggings. }
\begin{verbatim}
void path_class::calc_inner_label 
      (int i, int pathi) {
   int j,k;
   int tmp;
   for (j = 0; j < RIGSIZE; j++) {
      if (rigged[i][0][j] == UNUSED) break;
      if (rigged[i][3][j] == 1) {
         rigged [i][2][j] = rigged [i][1][j];
      }
   }
   for (j = 0; j < RIGSIZE; j++) {
      for (k = 1; k < RIGSIZE; k++) {
         if (rigged[i][0][k] == UNUSED) break;
         if (rigged[i][0][k] == rigged[i][0][k-1] 
            && rigged[i][1][k] == rigged[i][1][k-1] 
            && rigged[i][2][k-1] < rigged[i][2][k]){
            tmp = rigged[i][2][k-1];
            rigged[i][2][k-1] = rigged[i][2][k];
            rigged[i][2][k] = tmp;
         }
      }
   }
}

\end{verbatim}
\textit{
This inserts each element of a part in the path into the
bijection.
}
\begin{verbatim}
void path_class::insert_element_to_rigged (int num,
      int row_indx, int col_indx,int nrow, int ncol) 
{
   int old_largest_index;
   reset_flags(0);
   // add the new element - num - to rigged and 
   // add box if necessary
   for (int i = (num-2); i >= row_indx; i--) {
      old_largest_index = 
        add_box_to_rigged(i, num-2, 
          old_largest_index, 0);
   }
       
   // initialize curL to bigL
   for (int i=0; i < RIGSIZE; i++) 
      for (int j = 0; j < RIGSIZE; j++)
         curL[i][j] = bigL[i][j];

       // update curL to include the part of the
       //  tableau seen so far

       // we just finished a column
       if (row_indx == nrow - 1) {
          curL [nrow][1] += 1;
          if (col_indx != ncol - 1) {
             curL [nrow][ncol - col_indx - 1] += 1;
          }
          // we are in the middle of a column
          } else {
             curL [nrow][ncol - col_indx - 1] += 1;
             curL [row_indx + 1][1] += 1;
          }
          // calculate outer and inner labels 
          // based on curL
          for (int i = 0; i < n; i++) {
             calc_outer_label (i, 0);
          }
          for (int i = 0; i < n; i++) {
             calc_inner_label (i, 0);
          }
     
         // only if we are at the end of a column   
         // initialize curL to bigL - we'll update
         //  curL differently now in a FUSED way 
         // and recompute outer labels
         if (row_indx == nrow - 1) {
            for (int i=0; i < RIGSIZE; i++) 
               for (int j = 0; j < RIGSIZE; j++)
                 curL[i][j] = bigL[i][j];

            // update curL - 
            // this is different from above
            curL [nrow][ncol - col_indx] += 1;
            for (int i = 0; i < n; i++) {
               calc_outer_label (i, 0);
            }
         }
}

\end{verbatim}
\textit{ This inserts each part (which is a tableau) of a path to the bijection.}
\begin{verbatim}
void path_class::insert_tableau_to_rigged 
       (tblu_class* cur_tblu) {
    int nrow = cur_tblu->num_row;
    int ncol = cur_tblu->num_col;
    for (int j = ncol-1; j >=0; j--) {
       for (int i = 0; i < nrow; i++) {
          insert_element_to_rigged(
            cur_tblu->tb[i][j], i, j, nrow, ncol);
       }
    }
    bigL [nrow][ncol] = bigL [nrow][ncol] + 1;
    for (int i = 0; i < RIGSIZE; i++)
       for (int j = 0; j < RIGSIZE; j++) 
          curL[i][j] = bigL[i][j];
}

\end{verbatim}
\textit{This inserts all the parts of a path to the bijection.}
\begin{verbatim}
void path_class::build_rigged_for_path () {
    int j;
    for (int i=0; i < RIGSIZE; i++) {
        for (j =0; j < RIGSIZE; j++)
            bigL[i][j] = 0;
        curL[i][j] = 0;
    }
    for (int i = path_len - 1; i >= 0; i--) {
        insert_tableau_to_rigged (path[i]);
    }
}

\end{verbatim}
\textit{Calculates the $\alpha$ function used in the definition of cocharge. }
\begin{verbatim}
int path_class::alpha(int k, int i){
   int num_coln,j;
   num_coln=0;
   if (k>=n) num_coln=0; 
   else{
       for (j=0; j<RIGSIZE;j++){
         if (rigged[k][0][0]==UNUSED){
             num_coln=0 ;
             break;
          } else{
             if (rigged[k][0][j]!=UNUSED){ 
                if (rigged[k][0][j]>=i+1){
                num_coln=num_coln+1;
                }
             }
          }   
       }
   }
   return num_coln;
}

\end{verbatim}
\textit{Calculates the cocharge for the rigged configuration corresponding.}
\begin{verbatim}
void path_class::calculate_cocharge(){
  int k,j,i,sum,cosum;
  sum=0;
  for (k=0;k<=n-1;k++){
    if (rigged[k][0][0]!=UNUSED){
      for (i=0;i< rigged[k][0][0];i++){
      sum=sum+alpha(k,i)*(alpha(k,i)-alpha(k+1,i));
      }  
    }
  }     
  cosum=sum;
  
  for (k=0;k<=n-1;k++){
    for (j=0;j<RIGSIZE;j++){
      if (rigged[k][2][j]==UNUSED) break;
      if ( rigged[k][2][j]!=UNUSED){
        cosum=cosum + rigged[k][2][j];
      }
    }
  }
  cocharge=cosum;
  
  fprintf (stderr, "Statistic = %d \n", cosum );
}
 
\end{verbatim}
\textit{This is the main program.}
\begin{verbatim}
int main() {
    int i;
    read_input();   
            // this reads the input file.
    fprintf (stderr, "n=%d\n", n);
    input_path->print_path();              
            // this prints the input path
    input_path->build_rigged_for_path();   
            //finds the corresponding rigged 
           //configuration via the bijection.

    fprintf (stderr, "-------------------------\n");
    fprintf (stderr, "\n");
 
    fprintf (stderr, 
      "Corresponding rigged configuration is : \n");
        
    input_path->print_rigged_for_this_path();
       // prints the resulting rigged configuration. 
       
}
\end{verbatim}    
    
\section{Code for the program $inverse\_bijection.c$}

\textit{
This program does the inverse bijection from rigged configuration (RC) 
to path. Given a rigged 
configuration, n, path length and the shape of the path it calculates the 
corresponding path via the bijection. 
}

\begin{verbatim}
#include <stdio.h>
#define UNUSED 9999 
#define RIGSIZE 20
int n, l, num_shapes;
int lambda[100];
int path_shape[100][100];
int tableau_list[40000][10][10];
int tableau[100][100][100];
int pick[100];
int r, tab_indx, num_rc_lb_tab;
int rigged[RIGSIZE][5][RIGSIZE];
int bigL [RIGSIZE][RIGSIZE];
int curL [RIGSIZE][RIGSIZE];
int tblu_index , path_len;
int num_paths;
FILE *fp;

void initialize() {
   int i, j, k,m;
   for (i=0; i<RIGSIZE;i++){
      for (j=0; j<5; j++){
         for (k=0; k<RIGSIZE; k++){
            rigged[i][j][k]=UNUSED;
            curL[i][k]=0;
            bigL[i][k]=0;

         }
      }
   }
   for (i=0; i<100;i++){
      for (j=0; j<100; j++){
         for (k=0; k<100; k++){
            tableau[i][j][k]=UNUSED;
         }
      }
   }
   for (i=0;i<100; i++){
      for (j=0; j<100; j++){
         path_shape[i][j]=UNUSED;
      }
   } 
}
\end{verbatim}
\textit{Reads the input from file called "inputrigged".}
\begin{verbatim}
void read_input(){
    int i,j,k,tmp,tmu ;
    char c,c1;
    fp = fopen ("inputrigged","rw");
    fscanf(fp,"%d %d\n", &n, &path_len);

    for (i=0;i<=path_len-1;i++){
        for (j=0; j<2; j++){
            fscanf(fp,"%d",&tmu);
            path_shape[i][j]=tmu;
        }  
    }
    i = 0;
    j = 0;
    k = 0;
    while (fscanf (fp, "%d", &tmp)!= EOF){
      rigged[i][j][k] = tmp;
        k++;
        c=fgetc(fp);
        if (c=='\n'){
           k=0;
           if(j<1) j=j+2;
           else {j=0; i++;}
        }
    }
}
\end{verbatim}
\textit{Prints the input.}
\begin{verbatim}
void print_input() {
    int i,j;
    fprintf (stderr, "n = %d  L= %d\n", n, path_len);
    fprintf (stderr, "mu \n");
    for (i=0; i <=path_len-1; i++) {
       j=0;
       while (path_shape[i][j] != UNUSED){
           fprintf (stderr, "%d ", path_shape[i][j]);
           j++;
       } 
       fprintf(stderr,  "\n");
    }

    fprintf (stderr, "\n");
}

void reset_flags () {
    int i, k;
    for (i=0; i < RIGSIZE; i++) {
        for (k=0; k < RIGSIZE; k++) {
            rigged [i][3][k] = 0;
        }
    }
}
\end{verbatim}
\textit{
Prints the RC.
}
\begin{verbatim}
void print_rigged() {
   int i, j, k, l,a,b;
   fprintf(stderr, "Given rigged configuration is:\n");

   for (i = 0; i < n; i++) {
      if (rigged[i][0][0]==UNUSED) {
         fprintf (stderr,"------------------------\n");
         fprintf(stderr,"(%d)   Empty\n", i+1);
      }
      else {
         if (rigged[i][0][0] != UNUSED) {
            fprintf (stderr,"---------------------\n");
            fprintf(stderr, "(%d)\n", i+1);

            j=0;
            if (rigged[i][0][j] != UNUSED)
               for (k=0; k <rigged[i][0][j]; k++)
                  fprintf (stderr, " ___");
            fprintf (stderr,"\n");
            while (rigged[i][0][j] !=UNUSED){
               k=rigged[i][0][j];
               for (l=0; l<k-1;l++){
                   fprintf (stderr, "|   ");
               }
               if (l==k-1) 
               fprintf (stderr, "| %2d",rigged[i][2][j]);

               fprintf (stderr, 
                 "| %d\n",rigged[i][1][j]);
               for (l=0; l<k; l++){
                  fprintf (stderr, " ---");
               }
               fprintf (stderr,"\n");
               j++;
            }
          }
        }
    }
    fprintf (stderr,"------------------------------\n");
}
\end{verbatim}
\textit{Finds the smallest singular string in a middle RC.}
\begin{verbatim}
int find_smallest_inside_outside_others 
      (int index, int old_index) {
   int i = 0;
   int smallest_index = UNUSED;
   while ((rigged[index][0][i] != UNUSED) 
          && i < RIGSIZE) {
      if ((rigged[index][0][i] >= old_index) &&
          (rigged[index][1][i] == rigged[index][2][i]) 
          && (rigged[index][1][i] != UNUSED)) {
             smallest_index = i;
      }
      i++;
   }
   return smallest_index;
}
\end{verbatim}
\textit{Finds the smallest singular string in the starting RC.}
\begin{verbatim}
int find_smallest_inside_outside_first (int index) {
   int i = 0;
   int smallest_index = UNUSED;
   while ((rigged[index][0][i] != UNUSED) && 
           i < RIGSIZE) {
      if ((rigged[index][1][i] == rigged[index][2][i]) 
          && (rigged[index][1][i] != UNUSED)) {
             smallest_index = i;
      }
      i++;
   }
   return smallest_index;
}

\end{verbatim}
\textit{Calculates the new shape of the RC after removing a box.}
\begin{verbatim}
void remove_box_from_this_rigged 
       (int index, int column_index) {
   if (rigged [index][0][column_index]==1) {
      rigged [index][0][column_index]=UNUSED;
      rigged [index][3][column_index]=1;
   }
   else {
      rigged [index][0][column_index] -= 1;
      rigged [index][3][column_index] = 1;
   }
}

int num_box_1k_col (unsigned int i, int k){
   int j, l;
   int num_boxes = 0;
   for (l=1; l <= k; l++) {
      for (j = 0; j < RIGSIZE; j++){
        if (rigged[i][0][j] == UNUSED) break;
        if (rigged[i][0][j] >= l){
            num_boxes += 1;
        }
      }
   }
   return num_boxes;
}
\end{verbatim}
\textit{Finds the selected singular string and remove boxes from those parts.}
\begin{verbatim}
int remove_box_from_rigged 
      (int index, int begin, int old_smallest_index) {
   int smallest_index;
   if (index == begin) {
      smallest_index = 
        find_smallest_inside_outside_first(index);
   } else {
      smallest_index = 
        find_smallest_inside_outside_others (index, 
           old_smallest_index);
   }
   if (smallest_index == UNUSED)  return (-1);    
   else {
      remove_box_from_this_rigged (index, 
          smallest_index); 
         //this removes a box from selected part
      if (rigged[index][0][smallest_index]
           ==UNUSED) return 1;
      else return (rigged [index][0][smallest_index]+1); 
      // returns the length of the selected part       
   }
}
\end{verbatim}
\textit{
This calculates the extra term in the calculation of the vacancy num, 
which is the contribution from the shape of the path.}
\begin{verbatim}
int second_func(int part_size, int rig_num) {
    int sum = 0;
    for (int i=1; i<RIGSIZE; i++) {
        if (curL[rig_num+1][i] != 0) {
           int minimum = part_size;
           if (i < part_size) minimum = i;
           sum =sum + minimum * curL[rig_num+1][i];
        }
    }
    return sum;
}
\end{verbatim}
\textit{Calculates the vacancy numbers for each part of a rigged partition.}
\begin{verbatim}
void calc_outer_label(int rig_num){
   int part_num=0;
   int part_size, p;
   if (rig_num == 0)  {
      for (int part_num=0; part_num < RIGSIZE; 
            part_num++) {
         part_size = rigged[0][0][part_num];
         if(part_size == UNUSED) break;
         p = (-2*num_box_1k_col (0, part_size)) +
             (num_box_1k_col (1, part_size)) +
             second_func(part_size, rig_num);
         rigged [0][1][part_num] = p;
      }
   } else {
      for (part_num=0; part_num < RIGSIZE; part_num++) {
         part_size = rigged[rig_num][0][part_num];
         if(part_size == UNUSED) break;
         p = -2*num_box_1k_col (rig_num, part_size) +
             num_box_1k_col (rig_num-1, part_size) +
             num_box_1k_col (rig_num+1, part_size) +
             second_func(part_size, rig_num);
         rigged [rig_num][1][part_num] = p;
      }
   }
}
\end{verbatim}
\textit{Calculates the riggings for each part of a rigged partition.}
\begin{verbatim}
void calc_inner_label (int i) {
   int j,k;
   int tmp;
   for (j = 0; j < RIGSIZE; j++) {
      if (rigged[i][0][j] == UNUSED) break;
      if (rigged[i][3][j] == 1) {
         rigged [i][2][j] = rigged [i][1][j];
      }
   }
   for (j = 0; j < RIGSIZE; j++) {
      for (k = 1; k < RIGSIZE; k++) {
         if (rigged[i][0][k] == UNUSED) break;
         if (rigged[i][0][k] == rigged[i][0][k-1] 
            && rigged[i][1][k] == rigged[i][1][k-1] 
            && rigged[i][2][k-1] > rigged[i][2][k]) {
            tmp = rigged[i][2][k-1];
            rigged[i][2][k-1] = rigged[i][2][k];
            rigged[i][2][k] = tmp;
         }
      }
   }
}
\end{verbatim}
\textit{Calculates each element of a tableau in a path. }
\begin{verbatim}
void get_element_from_rc( int row_indx, 
       int col_indx, int pathi) {
   int old_smallest_index;
   reset_flags();
     
   // only if we are starting  a new column 
   // we'll update curL by spliting and 
   // recompute outer labels first
   if ((row_indx == path_shape[pathi][0]-1) 
       && (col_indx < path_shape[pathi][1]-1)) {
      curL[row_indx+1][1] +=1;
      curL [row_indx+1][path_shape[pathi][1]-col_indx] 
          -= 1;
      curL[row_indx+1][path_shape[pathi][1]-col_indx-1] 
          +=1;
      for (int i = 0; i < n; i++) {
         calc_outer_label (i);
      } 
    }

    // get a new element r from rigged and remove 
    // box if necessary
    int end=0;
    for (int i = row_indx ; i < n; i++) {
      old_smallest_index = 
         remove_box_from_rigged(i, row_indx, 
                        old_smallest_index);
      if (old_smallest_index == -1) {
         end=1; 
         tableau[pathi][row_indx][col_indx]=i+1;    
         break;
      }
    }
    if (end==0) {
        tableau[pathi][row_indx][col_indx]=n+1;
    }
    // update curL to exclude the part of the 
    // tableau seen so far

    // we just finished a column
    if (row_indx == 0) curL [1][1] -= 1; 

    // we are in the middle of a column
    else {
       curL [ row_indx+1 ][ 1 ] -= 1;
       curL [ row_indx ][ 1 ] += 1;
    }

    // calculate outer and inner labels based 
    // on updated curL
    for (int i = 0; i < n; i++) {
       calc_outer_label (i);
    }
    if (tableau[pathi][row_indx][col_indx]
        !=(row_indx+1)) {    
       for (int i = 0; i < n; i++) {
          calc_inner_label (i);
       }
    }
}
\end{verbatim}
\textit{
Calculates the rigged configuration for one tableau in the path.
}
\begin{verbatim}
void get_tableau_from_rc(int pathi) {
   int nrow = path_shape[pathi][0];
   int ncol = path_shape[pathi][1];
   for (int j = 0; j <= ncol-1; j++) {
      for (int i = nrow-1; i >=0; i--) {
         get_element_from_rc(i,j,pathi);
      }
   }
   bigL [nrow][ncol] = bigL [nrow][ncol] - 1;
   for (int i = 0; i < RIGSIZE; i++)
   for (int j = 0; j < RIGSIZE; j++) 
      curL[i][j] = bigL[i][j];
   fprintf(stderr, "-------------------------\n");
   for (int i=0; i <= nrow-1; i++) {
      for (int j=0; j<= ncol-1; j++) {
        fprintf(stderr, "%2d", tableau[pathi][i][j]);
      }
      fprintf(stderr, "\n");
   }
}
\end{verbatim}
\textit{Calculates the rigged configuration for a given path.}
\begin{verbatim}
void build_path_for_rc() {
   int i,j,k,tmp;
   for (i=0; i < RIGSIZE; i++) {
      for (j =0; j < RIGSIZE; j++) {
         if ((j==0) && (path_shape[i][j] != UNUSED)) {
            bigL[path_shape[i][j]][path_shape[i][j+1]] += 1;
            curL[path_shape[i][j]][path_shape[i][j+1]] += 1; 
         } 

      }
   }
   for (int i = 0; i < n; i++) {
      calc_outer_label (i);
   }
   for (i=0;i<n; i++) { 
      for (j = 0; j < RIGSIZE; j++) {
         for (k = 1; k < RIGSIZE; k++) {
            if (rigged[i][0][k] == UNUSED) break;
            if (rigged[i][0][k] == rigged[i][0][k-1] 
                && rigged[i][1][k] == rigged[i][1][k-1] 
                && rigged[i][2][k-1] > rigged[i][2][k]) {
               tmp = rigged[i][2][k-1];
               rigged[i][2][k-1] = rigged[i][2][k];
               rigged[i][2][k] = tmp;
            }
         }
      }
   }
   
   print_rigged();
    
   fprintf(stderr,"The corresponding path is:\n");
   for (i = 0; i < path_len; i++) {
      get_tableau_from_rc (i);
   }
   fprintf(stderr, "-------------------------\n");
}

\end{verbatim}
\textit{Main program.}
\begin{verbatim}
int main(){
    initialize();
    read_input();
    print_input();
    build_path_for_rc();
    
}
\end{verbatim}
\end{subappendices}

    %
    %

\end{document}